\documentclass[11pt,reqno,twoside]{amsart}
\usepackage{amscd}
\usepackage{amsfonts}
\usepackage{amsmath}
\usepackage{amssymb}
\usepackage{amsthm}
\usepackage{amsaddr}
\usepackage{fancyhdr}
\usepackage{latexsym}
\usepackage[colorlinks=true, pdfstartview=FitV, linkcolor=blue, citecolor=blue, urlcolor=blue]{hyperref}
\usepackage{enumitem}      
\usepackage{mathtools}            
\usepackage{indentfirst} 
\usepackage{color}
\usepackage{caption}          
\usepackage[normalem]{ulem}
\usepackage{subfig}
\usepackage{floatrow}
\newcommand{\nn}{\nonumber}
\newcommand{\p}{\partial}

\usepackage{thmtools}
\declaretheoremstyle[headfont=\kpfonts]{normalhead}

\newtheorem{theorem}{Theorem}[section]

\newtheorem{corollary}{Corollary}[section]

\theoremstyle{definition}

\newtheorem{remark}{Remark}[section]

\allowdisplaybreaks
\numberwithin{equation}{section}
\numberwithin{figure}{section}
\usepackage{geometry}
\makeatletter
\@namedef{subjclassname@2020}{%
  \textup{2020} Mathematics Subject Classification}
\makeatother

\begin{document}

\title{Long-time asymptotics and the radiation condition 
\\
for linear evolution equations on the half-line \\ with time-periodic boundary conditions}

\author{Yifeng Mao$^{*}$, Dionyssios Mantzavinos$^{\dagger}$, and Mark A. Hoefer$^{*}$}

\address{\normalfont $^{*}$Department of Applied Mathematics, University of Colorado, Boulder, CO, 80309, USA
\\[1mm]
$^{\dagger}$Department of Mathematics, University of Kansas, Lawrence, KS, 66045, USA}
 \email{mantzavinos@ku.edu}
\thanks{\textit{Acknowledgements.} The authors would like to thank the Isaac Newton Institute for Mathematical Sciences, Cambridge, U.K. for support and hospitality during the program ``Dispersive Hydrodynamics'', when work on this
paper was undertaken (EPSRC Grant Number EP/R014604/1).
D.M. gratefully acknowledges support from the U.S. National
Science Foundation (NSF-DMS 2206270).  Y.M.~and M.A.H.~acknowledge support from NSF-DMS 1816934.}
\date{July 27, 2023}

\begin{abstract}
  The large time $t$ asymptotics for scalar, constant coefficient,
  linear, third order, dispersive equations are obtained for
  asymptotically time-periodic Dirichlet boundary data and zero
  initial data on the half-line modeling a wavemaker acting upon an
  initially quiescent medium.  The asymptotic Dirichlet-to-Neumann
  (D-N) map is constructed by expanding upon the recently developed
  $Q$-equation method.  The D-N map is proven to be unique if and only if the
  radiation condition that selects the unique wavenumber branch of the
  dispersion relation for a sinusoidal, time-dependent boundary
  condition holds: (i) for frequencies in a finite interval, the wavenumber
  is real and corresponds to positive group velocity, (ii) for
  frequencies outside the interval, the wavenumber is complex with
  positive imaginary part.  For fixed spatial location $x$, the
  corresponding asymptotic solution is (i) a traveling wave or (ii) a
  spatially decaying, time-periodic wave.  Uniform-in-$x$
  asymptotic solutions for the physical cases of the linearized
  Korteweg-de Vries and Benjamin-Bona-Mahony (BBM) equations are obtained via integral asymptotics.  The linearized BBM asymptotics
  are found to quantitatively agree with viscous core-annular fluid
  experiments.
\end{abstract}

\maketitle
\markboth
{Yifeng Mao, Dionyssios Mantzavinos, and Mark A. Hoefer}
{Long-time asymptotics and the radiation condition for linear evolution equations on the half-line \ldots 
}

\section{Introduction}

The forced generation of waves from a boundary is a fundamental
problem in applied mathematics with broad applications in fluid
dynamics \cite{buhler_waves_2014}, material science \cite{wu_nonlinear_2010},
seismology \cite{aki_quantitative_2002}, and astronomy
\cite{cassini_imaging_team_ssi_jpl_esa_nasa_astronomy_2017} just to
name a few.  The fundamental mathematical problem can be posed as
follows.  Consider a linear evolutionary equation subject to
time-harmonic forcing at a localized source or boundary
$\propto e^{-i \omega_0 t}$ and decay at infinity.  We would expect
the solution to exhibit spatial dependence via the wavevector
$\mathbf{k}_0$ consistent with the system's dispersion relation such
that
$u(\mathbf{x},t) \propto e^{i(\mathbf{k}_0 \cdot \mathbf{x} - \omega_0
  t)}$ near the source or boundary.  However, for any evolutionary
equation whose spatial order is higher than one, the wavevector
$\mathbf{k}_0$ is not unique.  In the context of the time-independent
Helmholtz equation $(\nabla^2+k_0^2)u = -f(\mathbf{x})$, Sommerfeld's
eponymous \textit{radiation condition} in three-dimensional space,
$\lim_{r \to \infty} r(\frac{\partial u}{\partial r} - i k_0 u) = 0$,
was formulated to select the unique wave that propagates away from a
source or boundary to ensure well-posedness on an infinite domain
\cite{sommerfeld_partial_1949}.  Another, closely related approach to
this problem is called the limiting absorption principle for the
resolvent $R(k_0) = (\Delta + k_0^2)^{-1}$ \cite{sveshnikov_1950}.  The unique solution is
selected by taking a limit of the resolvent whose spectrum is
perturbed off the real, essential spectrum into the appropriate half
of the complex plane.  Ultimately, this boils down to selecting the
appropriate branch of the square root.

Both the radiation condition and the limiting absorption principle are
mathematical workarounds to the fundamental problem that the Helmholtz
equation exhibits two linearly independent, decaying solutions (in two
or more dimensions).  A physically-inspired approach to the
problem is to consider a wavemaker in which a boundary or source
time-harmonically forces an initially quiescent
($u(\mathbf{x},0) = 0$) medium.  In the context of the forced wave equation
$u_{tt} - \Delta u = - e^{-i \omega_0 t}f(\mathbf{x})$, solving this
initial-boundary value problem and taking $t \to \infty$ is known as
the limiting amplitude principle, which selects the same outgoing wave
as the Sommerfeld radiation condition and corresponds to the causal solution 
by perturbing the essential 
spectrum of the Helmholtz resolvent into the upper half of the complex plane
\cite{tikhonov_radiation_1948,morawetz_limiting_1962}.

A simple, physically inspired calculation \cite{buhler_waves_2014}
demonstrates the limiting amplitude principle.  Consider a linear, first order in time
evolutionary equation for $u(x,t)$ on the half-line $x > 0$ with
real-valued dispersion relation $\omega = \omega(k)$, the asymptotic
Dirichlet boundary condition $u(0,t) \sim e^{-i \omega_0 t}$,
$t \to \infty$ with $\omega_0 \in \mathbb{R}$, and decay
$\lim_{x \to \infty} u(x,t) = 0$ for each $t \in \mathbb{R}$.  Anticipating the wave
\begin{equation}
  \label{eq:asymptotic_traveling_wave}
  u(x,t) \sim e^{i(k_0 x - \omega_0 t)}, \quad t \to \infty, \quad x = \mathcal{O}(1),
\end{equation}
we seek the correct wavenumber branch $k_0 = k(\omega_0)$ selected
from the multivalued inverse $k(\omega_0) = \omega^{-1}(\omega_0)$.
The limiting amplitude principle states that, initially, the solution
is zero, $u(x,0) = 0$ for each $x > 0$.  Consequently, the
\textit{causal} solution can be modeled by perturbing the frequency
into the upper half of the complex plane
$\omega_0 \to \omega_0 + i \epsilon$, $0 < \epsilon \ll |\omega_0|$ so
that the boundary condition becomes
$u(0,t) = e^{-i \omega_0 t} e^{\epsilon t}$ with initially slow
amplitude growth to eventually attain the asymptotically periodic
boundary condition.  Assuming that $k(\omega)$ is analytic in a
neighborhood of $\omega_0$, we expand
\begin{equation}
  \label{eq:1}
  k(\omega_0 + i \epsilon) = k(\omega_0) + i \epsilon \frac{\mathrm{d}
    k}{\mathrm{d} \omega}(\omega_0) + \cdots = k_0 + i
  \frac{\epsilon}{c_g(\omega_0)} + \cdots ,
\end{equation}
where $c_g \equiv \mathrm{d} \omega/\mathrm{d} k$ is the group
velocity.  Then, the anticipated solution takes the form
\begin{equation}
  \label{eq:3}
  u(x,t) \sim e^{\epsilon t} e^{i(k(\omega_0)x - \omega_0 t)} =
  e^{\epsilon [t - x/c_g(\omega_0)] + \cdots} e^{i(k_0x - \omega_0 t)} .
\end{equation}
Since $u(x,t) \sim e^{(ik_0 - \epsilon /c_g)x - (i \omega_0 - \epsilon) t}$, in order to ensure $\lim_{x \to \infty} u(x,t) = 0$ in \eqref{eq:3} as $\epsilon \to 0^+$, we require
\begin{equation*}
  k_0 \in \mathbb{R}~ \mathrm{and} ~c_g(\omega_0) > 0, \quad \mathrm{or} \quad \mathrm{Im}
  \,k_0 > 0 . 
\end{equation*}
These inequalities are the
radiation condition for the
wavemaker problem.  The
positive group velocity for real wavenumbers corresponds to waves being emitted from the
boundary into the domain.  For complex $k$, the positive imaginary part implies spatial decay.  The boundary case $c_g(\omega_0) = 0$ is not included in this argument because $k(\omega)$ has a branch point at $\omega_0$; 
therefore it is not analytic.  Zero group velocity corresponds to a critical value of the frequency that represents a transition between propagating and spatially decaying waves.  We will prove the general radiation condition  
\begin{equation}
  \label{eq:general_radiation_condition}
  k_0 \in \mathbb{R}~ \mathrm{and} ~c_g(\omega_0) \ge 0, \quad \mathrm{or} \quad \mathrm{Im}
  \,k_0 > 0 . 
\end{equation}
Although we have given a physically
inspired, plausible argument for the radiation condition
\eqref{eq:general_radiation_condition}, its proof is
more involved.

In this paper, we prove that the radiation condition
\eqref{eq:general_radiation_condition} uniquely determines the
asymptotic form of the solution \eqref{eq:asymptotic_traveling_wave}
by constructing the Dirichlet-to-Neumann (D-N) map for the
wavemaker initial-boundary value problem
\begin{subequations}\label{gen-ibvp}
\begin{align}
&\big(1-A_{-2} \p_x^2\big) \, u_t = i A_0 u + A_1 u_x + i A_2 u_{xx} + A_3 u_{xxx},
\quad x, t > 0,
\label{gen-mod}
\\
&u(x, 0) = 0, \quad x>0,
\\
&u(0, t) = g_0(t), \quad t>0,
\label{gen-bc}
\end{align}
\end{subequations}
where $A_{-2}\geq 0$ to ensure well-posedness,
$A_0, A_1, A_2 \in \mathbb R$, and $A_3 \leq 0$ so that only one
boundary condition is required at $x=0$. By scaling $t$ and $x$, we
can restrict $A_3 \in \{-1,0\}$ and $A_{-2} \in \{0,1\}$ without loss
of generality.
For real wavenumbers, these restrictions guarantee the real-valued
dispersion relation
\begin{equation}
  \label{eq:5}
  \omega(k) = \frac{-A_0-A_1k+A_2 k^2+ A_3 k^3}{1+A_{-2}k^2} .
\end{equation}
The choice of parameters $\left\{A_{-2}=A_0=A_2=0, A_1=A_3=-1\right\}$
corresponds to the linear Korteweg-de Vries (KdV) equation, also known
as the first Stokes equation, while for the parameters
$\left\{A_{-2}=1, A_1=-1, A_0=A_2=A_3=0\right\}$, equation
\eqref{gen-mod} becomes the linear Benjamin-Bona-Mahoney (BBM)
equation.

In 1997, Fokas introduced a new approach for the analytical solution of initial-boundary value problems \cite{fokas1997unified}. The main motivation for Fokas's method, otherwise known as the unified transform, was twofold: 
\begin{enumerate}[label=(\roman*), topsep=1.5mm, itemsep=1mm, leftmargin=7.6mm]
\item to provide an analogue of the inverse scattering transform method, which is used for solving the initial value problem of completely integrable systems like the KdV and NLS equations, when such systems are considered in domains that involve a boundary, e.g. the half-line or the finite interval, thereby giving rise to nonlinear initial-boundary value problems;
\item
to provide a natural analogue of the Fourier transform, which yields the solution to linear evolution equations in the initial value problem setting, when these linear equations are instead formulated in the setting of initial-boundary value problems.
\end{enumerate}
While the difficulty associated with task (i) can be easily appreciated due to the presence of nonlinearity, task (ii)  comes with significant challenges as well, despite the fact that it concerns linear equations. This becomes evident through the realization that, for general nonzero boundary conditions, no classical spatial transform is available for solving initial-boundary value problems that involve linear evolution equations of spatial order higher than two. Moreover, even for second-order linear equations, classical techniques fail when the boundary data are non-separable.

The Fokas method effectively addresses both of the above objectives as a unified way of treating initial-boundary value problems for a plethora of linear and nonlinear equations, formulated in various domains and supplemented with different kinds of boundary conditions. Since its inception, the method has been further developed within the completely integrable systems and applied analysis communities in various directions such as the analytical solution of nonlinear equations via the Riemann-Hilbert/dbar problem approach, the solution of linear evolution and elliptic equations, asymptotics, and numerical studies; see the books \cite{fokas2008unified,fbook2015} as well as the review articles \cite{fs2012,dtv2014} for a detailed list of references on the above topics and beyond.

It should be emphasized that the solution formulae obtained via the unified transform have significant advantages even when the relevant problems can be solved via classical approaches. Importantly, they involve integrals with exponentially decaying integrands, due to the fact that integration takes place along complex contours in the spectral plane. As a result, they are uniformly convergent up to the boundary of the domain --- a feature that turns out particularly helpful for our purposes here. 

In this work, we employ the Fokas method to prove the existence and
uniqueness of the Dirichlet-to-Neumann (D-N) map for the wavemaker
problem \eqref{gen-ibvp} and use it to construct a series
representation of the solution for fixed $x>0$ as $t\to\infty$,
thereby providing a constructive proof of the radiation condition
\eqref{eq:general_radiation_condition}. We also obtain the large $t$
asymptotics along rays $x/t = \text{constant}$ of the wavemaker
problem for the physically important linear KdV and linear BBM
equations. Finally, we report experiments on a two-fluid system that
agree quantitatively with solutions of the linear BBM wavemaker
problem.

We are now ready to precisely state the above results. We begin with
the former result, which is actually valid for generic, nonzero
initial data provided these have sufficient smoothness and decay
(e.g. of Schwartz type).
\begin{theorem}[Dirichlet-to-Neumann map]
\label{dn-t}
Consider the evolution equation \eqref{gen-mod} formulated on the half-line $\{x>0\}$ with a sufficiently smooth and decaying initial datum and the Dirichlet boundary datum \eqref{gen-bc} that is infinitely smooth and asymptotically periodic in time, i.e. $g_0 \in C^\infty(0, \infty)$ and
$$
g_0(t) \sim \sum_{n\in\mathbb Z} a_n e^{in\omega_0 t}, \quad t\to\infty, \quad \omega_0\geq 0, \ a_n \in \mathbb C
$$
with the infinite sum being absolutely convergent.  Assume further that the initial datum and the boundary datum $g_0(t)$ satisfy the necessary compatibility conditions at the origin.
Suppose that either $\{A_{-2}=0, A_3<0\}$ (local, third-order) or $\{A_{-2}>0, A_3=0\}$ (nonlocal, second-order) with $A_2/(\omega_0 A_{-2})\notin \mathbb Z$ in the latter case.
Then, under the assumption that the solution $u(x, t)$ and its  boundary values $\p_x^j u(0, t)$, $j\in\mathbb N$, are asymptotically time-periodic, the associated generalized D-N map satisfies to leading order the asymptotic formula
\begin{equation}\label{gdn-map}
\p_x^j u(0, t) \sim i^j \sum_{n\in\mathbb Z} [k_0(n)]^j a_n e^{in\omega_0t}, \quad t\to\infty, \  j\in\mathbb N\cup\{0\},
\end{equation}
where $k_0(n)$ is the unique solution of the equation
$
A_3 k^3 + (A_2+n\omega_0 A_{-2}) k^2 -A_1 k - A_0 + n\omega_0 = 0
$
that lies along the positively oriented boundary $\p D^+$ of the region
$D^+ = \big\{k\in\mathbb C: \textnormal{Im}(k) > 0 \textnormal{ and } \textnormal{Re}\big( \frac{i(A_3 k^3 + A_2 k^2 -A_1 k - A_0)}{ 1+A_{-2}k^2}\big) < 0 \big\}$.
\end{theorem}

Theorem \ref{dn-t} is established in Section \ref{D2N-s}, where the
overall framework for the general, five-parameter
equation~\eqref{gen-mod} is first laid out and the analysis specific
to the cases $\{A_{-2}=0, A_3<0\}$ (local, third-order) and
$\{A_{-2}>0, A_3=0\}$ (nonlocal, second-order) addressed by Theorem
\ref{dn-t} is then carried out in detail.  

\begin{corollary}[Radiation condition]
  \label{rad-cond}
  Under the assumptions of Theorem \ref{dn-t}, $k_0(n)$ is the unique
  branch of the dispersion relation \eqref{eq:5} satisfying
  $n\omega_0 = -\omega(k_0(n))$ if and only if the radiation condition
  \eqref{eq:general_radiation_condition} holds.  Furthermore, the asymptotic
  solution of the wavemaker problem \eqref{gen-ibvp} takes the form
  \begin{equation}
    \label{eq:asymptotic_soln}
    u(x,t) \sim \sum_{n \in \mathbb{Z}} a_n e^{i(k_0(n)x + n \omega_0 t)},
    \quad t \to \infty, \quad x = \mathcal{O}(1) .
  \end{equation}
  Finally, there exist real numbers $\omega_{cr}^{\pm}$ such that
  $\mathrm{Im}\, k_0(1) = 0$ if and only if
  $\omega_{cr}^- \le \omega_0 \le \omega_{cr}^+$.
\end{corollary}

A few remarks are in order.
\begin{enumerate}[label=\arabic*.,topsep=1.5mm, itemsep=1mm, leftmargin=5.6mm]
\item In the nonlocal, second-order case $\{A_{-2}>0, A_3=0\}$, the scenario $A_2/(\omega_0 A_{-2}) \in \mathbb Z$, which is not covered by Theorem \ref{dn-t}, is discussed separately at the end of Subsection \ref{2-nonlocal}.
\item The assumption $g_0 \in C^\infty(0, \infty)$ implies that the
  Fourier series coefficients $a_n$ decay faster than any power of
  $n \in \mathbb Z$, thus the infinite sum in \eqref{gdn-map} is
  absolutely convergent. To see this, note that $k_0(n)$ grows like
  $n^{\frac 13}$ and $n^{\frac 12}$ in the third-order and
  second-order cases, respectively.
\item The uniqueness of the D-N map \eqref{gdn-map} follows from the
  requirement that $k_0(n)$ is the unique root of the cubic
$
A_3 k^3 + (A_2+n\omega_0 A_{-2}) k^2 -A_1 k - A_0 + n\omega_0
$
that lies along the boundary of the domain $D^+$.
As explained in Section~\ref{D2N-s}, this requirement is easy to
justify for $n\omega_0$ outside of the interval
$[\omega_{cr}^-,\omega_{cr}^+]$ (supercritical regime). However,
within the interval (subcritical regime), things become more subtle
and uniqueness follows by studying the large $t$ limit of the Fokas
solution formulae for fixed $x>0$.
This is done in Section \ref{D2N-s} for the linear KdV and linear BBM
equations, where, by taking advantage of the analyticity and
exponential decay properties of the integrands involved in the
solution formulae, we perform suitable contour deformations that allow
us to quickly derive the leading-order, long-time asymptotics of the
corresponding solutions for each fixed $x>0$.  We emphasize that this
approach is not valid when $t\to\infty$ along rays
$x/t = \text{constant}$, which is addressed in Theorems \ref{lkdv-t}
and \ref{lbbm-t}.
Importantly, these calculations show that the choice of $k_0(n)$ as
the unique root of the aforementioned cubic that lies on $\p D^+$ is
directly related to the fact that $D^+$ corresponds to the region in
the upper half of the complex $k$-plane where the imaginary part of
the dispersion relation $k(\omega)$ of equation \eqref{gen-mod} is
positive, as in the radiation condition
\eqref{eq:general_radiation_condition}. Moreover, in the general case
covered by Theorem \ref{dn-t} that goes beyond the linear KdV and BBM
equations, the same analysis can be carried out to reach the same
criterion for $k_0(n)$. This analysis addresses a relevant uniqueness
issue that was identified in \cite{fv2021}.
\item In the case of real roots of the aforementioned cubic, the root
  $k_0(n)$ is the only one that corresponds to a non-negative group
  velocity.
Indeed, the group velocity associated with equation \eqref{gen-mod} is 
$$
\frac{d}{dk} \left(\frac{A_3k^3+A_2k^2-A_1k-A_0}{1+A_{-2}k^2}\right)
=
\frac{A_{-2}A_3k^4+(3A_3+A_{-2}A_1)k^2+2(A_2+A_{-2}A_0)k-A_1}{\left(1+A_{-2}k^2\right)^2}.
$$
Thus, for $k\in\mathbb R$, the condition that $k \in \p D^+$ is
equivalent to  the group velocity being non-negative at that solution
as in the radiation condition \eqref{eq:general_radiation_condition}. 
In the case of complex roots, the root $k_0(n)$ is the unique root
with $\mathrm{Im} \, k_0(n) > 0$ that leads to a spatially decaying
solution.
\item
The calculations performed in Section \ref{D2N-s} are under the assumption of smooth data that lead to a smooth solution with sufficient decay as $x\to\infty$. Nevertheless, depending on the precise choice of
coefficients in \ref{gen-mod}, these conditions can be relaxed significantly --- see, for example, the works \cite{fokas2016korteweg, himonas2019korteweg} where the Fokas method is combined with harmonic analysis techniques for the study of well-posedness of the KdV equation on the half-line and the finite interval with data in appropriate Sobolev spaces. Such an investigation of optimal regularity properties is not our task here.

\item 
A major assumption in the hypothesis of Theorem \ref{dn-t} is that, for zero initial data and asymptotically periodic Dirichlet data, the solutions of the models corresponding to \eqref{gen-mod} with the choices $\{A_{-2}=0, A_3<0\}$ and $\{A_{-2}>0, A_3=0\}$, as well as the various boundary values of these solutions, are asymptotically periodic in time.
The explicit formulae obtained via the Fokas method for the wavemaker
problem of the linear KdV and the linear BBM equations allow us to
verify this fact and hence remove the time-periodicity assumption from
the hypothesis of Theorem \ref{dn-t} in the case of these two physical
models.  More specifically, in Section 2 it is shown that by taking
the asymptotic limit $x$ fixed and $t\to\infty$ in the solution
formulae, we obtain the asymptotically time-periodic solution
\eqref{eq:asymptotic_soln} and therefore verify the periodicity
assumption in Theorem \ref{dn-t}. Then, in Sections 3 and 4, we
combine the solution formulae with the steepest descent method to
derive the leading-order asymptotics of these solutions as
$t\to\infty$ and $x/t=\mathcal{O}(1)$, i.e. the behavior of these
solutions along rays in the $x$-$t$ plane. This type of limit provides
detailed information about the propagation of waves away from the
origin in the wavemaker problem.  We also obtain a uniform asymptotic
approximation by demonstrating that the asymptotic results from each
region $x = \mathcal{O}(1)$ and $x = \mathcal{O}(t)$ match as
$t \to \infty$.  The precise asymptotic solution statements to the
wavemaker problems are as follows.
\end{enumerate}

\begin{theorem}[Linear KdV wavemaker problem]
\label{lkdv-t}
Consider the linear KdV wavemaker problem with input frequency $\omega_0$,
\begin{align*}
    u_t + u_x + u_{xxx}  = 0, \quad & 0<x<\infty, \quad t>0, \\
    u(x,0) = 0, \quad  & 0 \le x<\infty, \\
    u(0,t) = \sin(-\omega_0 t), \quad & t \ge 0.
\end{align*}
There exists a critical frequency $\omega_{cr}:=2/(3\sqrt{3})$.   Assuming that the solution is sufficiently smooth, its leading-order asymptotic behavior as  $t\to\infty$ and $\xi := x/t=\mathcal{O}(1)$ is:
\begin{enumerate}[label=(\Roman*), leftmargin=7.6mm]
\item A time-periodic sinusoidal function, when
  $0<\omega_0<\omega_{cr}$ and $0<\xi<c_g:=1-3k_0^2<1$ with $k_0 =
  k_3(\omega_0)$ given by  \eqref{eq:KdV-poles-sub}, namely 
  \begin{equation}\label{eq:LKdV_soln_plane-intro} 
    u(\xi t,t) \sim \sin \left[ (k_0 \xi -\omega_0) t \right].
  \end{equation}
  Moreover, this is a uniform asymptotic approximation that is also
  valid for $\xi t = \mathcal{O}(1)$, $t \to \infty$.
\item A sinusoidal function in $t$ with an algebraically decaying
  amplitude of $\mathcal{O}(t^{-1/2})$, when $0<\omega_0<\omega_{cr}$
  and $c_g<\xi<1$, or when $\omega_0>\omega_{cr}$ and $0<\xi<1$ with
  the condition $l_1(\xi,\omega_0)<0$ where $l_1$ is defined in
  \eqref{eq:LKdV_l1}, i.e.
  \begin{equation}
    u(\xi t,t) \sim \frac{27 \xi  \omega_0 \cos\left( \frac{2}{9} t
        \sqrt{3(1-\xi)} (1-\xi) \right) }{-(\xi +3) \xi ^2-27
      \omega_0^2+4}  \sqrt{\frac{1}{2\pi
        t\sqrt{3(1-\xi)}}}.\label{eq:LKdV_xi<1_decay_soln-intro} 
  \end{equation}
\item A sinusoidal function in $t$ with an algebraically decaying
  amplitude of $\mathcal{O}(t^{-1/2})$ associated with an
  exponentially decaying sinusoidal function when
  $\omega_0>\omega_{cr}$ and $0<\xi<1$ with the condition
  $l_1(\xi,\omega_0)>0$ where $l_1$ is defined in \eqref{eq:LKdV_l1},
  that is
  \begin{equation}
    u(\xi t,t) \sim  e^{-\mathrm{Im}(k_0) \xi t} \sin \left[
      (\mathrm{Re}(k_0) \xi - \omega_0) t \right] + \frac{27 \xi
      \omega_0 \cos\left( \frac{2}{9} t \sqrt{3(1-\xi)} (1-\xi)
      \right) }{-(\xi +3) \xi ^2-27 \omega_0^2+4}  \sqrt{\frac{1}{2\pi
        t\sqrt{3(1-\xi)}}},  \label{eq:LKdV_xi<1_super_soln-intro} 
  \end{equation}
  with $k_0 = k_3(\omega_0)$ now given by \eqref{eq:KdV-poles-sup}.
  While the first term is beyond all orders in $t$ when
  $\xi = \mathcal{O}(1)$, this is a uniform asymptotic approximation
  that is also valid for $\xi t = \mathcal{O}(1)$, $t \to \infty$,
  where the first term is the dominant contribution.
\item An exponential, decaying in time function when $\xi>1$, i.e.
  \begin{equation}
    u(\xi t,t) 
    \sim \frac{27 \xi  \omega_0 e^{-\frac{2}{9} \sqrt{3\xi-3} (\xi-1)
        t} }{-2(\xi +3) \xi ^2-54 \omega_0^2+8}  \sqrt{\frac{1}{t \pi
        \sqrt{3(\xi-1)}}}. \label{eq:KdV_xi>1_soln_2-intro}
  \end{equation}
\end{enumerate}

\end{theorem}

\begin{theorem}[Linear BBM wavemaker problem]
\label{lbbm-t}
Consider the linear BBM wavemaker problem with input frequency $\omega_0$,
\begin{align*}
    u_t + u_x - u_{xxt}  = 0, \quad & 0<x<\infty, \quad t>0, \\
    u(x,0) = 0, \quad  & 0 \le x<\infty, \\
    u(0,t) = \sin(-\omega_0 t), \quad & t \ge 0.
\end{align*}
There exists a critical frequency $\omega_{cr}:=1/2$.   Assuming that the solution is sufficiently smooth, its leading-order asymptotic behavior as  $t\to\infty$ and $\xi := x/t=\mathcal{O}(1)$ is described by:
\begin{enumerate}[label=(\Roman*), leftmargin=7.6mm]
\item A time-periodic sinusoidal function, when
  $0<\omega_0<\omega_{cr}$ and
  $0<\xi<c_g:=\frac{1-4\omega_0^2+\sqrt{1-4\omega_0^2}}{2} <1$, namely
  \begin{equation}\label{eq:LBBM_soln_plane-intro} 
    u(\xi t,t) \sim \sin \left[ \left(
        \frac{1-\sqrt{1-4\omega_0^2}}{2\omega_0} \xi - \omega_0
      \right) t   \right].
  \end{equation}
  Moreover, this is a uniform asymptotic approximation that is also
  valid for $\xi t = \mathcal{O}(1)$, $t \to \infty$.
\item A sinusoidal function in $t$ with an algebraically decaying
  amplitude of $\mathcal{O}(t^{-1/2})$, when $0<\omega_0<\omega_{cr}$
  and $c_g<\xi<1$, or when $\omega_0>\omega_{cr}$ and $0<\xi<1$ with
  the condition $l_3(\xi,\omega_0)<0$ where $l_3$ is defined in
  \eqref{eq:BBM_kappa_pos_2}, i.e.
  \begin{small}
  \begin{equation}
 \begin{aligned}
    &u(\xi t,t) := u_{s1}(\xi t,t) 
    \\
    &\sim \frac{ \omega _0 \left(-4 \xi
        +\sqrt{8\xi +1}-1\right)  \left(\frac{\xi^{1/2} 
          \left(\sqrt{8 \xi +1}-1\right)^3}{ t \left(8 \xi + 1 
            -\sqrt{8 \xi +1}\right) \sqrt{-2 \xi -1 +\sqrt{8 \xi 
              +1} }    }\right)^{1/2} \cos\left[\frac{2\sqrt{2\xi}
        }{\left(\sqrt{8 \xi +1}-1\right)^2} \left(-2 \xi +\sqrt{8 \xi
            +1}-1 \right)^{3/2}  t  \right] }{2^{5/4} \sqrt{\pi } 
      \left(\left(4 \xi -\sqrt{8 \xi +1}+1\right) \omega_0^2+\xi
        \left(2 \xi -\sqrt{8 \xi+1}+1\right)\right)}. 
        \end{aligned}
    \label{eq:LBBM_xi<1_decay_soln-intro}
  \end{equation}
  \end{small}
\item A sinusoidal function in $t$ with an algebraically decaying
  amplitude of $\mathcal{O}(t^{-1/2})$ associated with an
  exponentially decaying sinusoidal function when
  $\omega_0>\omega_{cr}$ and $0<\xi<1$ with the condition
  $l_3(\xi,\omega_0)>0$ where $l_3$ is defined in
  \eqref{eq:BBM_kappa_pos_2}, that is
  \begin{equation}
    u(\xi t,t) \sim  \exp \left[-\frac{t \xi
        \sqrt{4\omega_0^2-1}}{2\omega_0}\right] \sin\left[\left(
        \frac{\xi}{2\omega_0}-\omega_0 \right) t \right] + u_{s1}(\xi
    t,t).  \label{eq:LBBM_xi<1_super_soln-intro} 
  \end{equation}
  While the first term is beyond all orders in $t$ when
  $\xi = \mathcal{O}(1)$, this is a uniform asymptotic approximation
  that is also valid for $\xi t = \mathcal{O}(1)$, $t \to \infty$,
  where the first term is the dominant contribution.
\item An exponentially decaying in time function when $\xi>1$, i.e.
\begin{small}
  \begin{equation}
    u(\xi t,t) 
    \sim -\frac{ \omega _0
      \left(\frac{ \xi^{1/2} \left(\sqrt{8 \xi +1}-1\right)^3  }{ t
          \sqrt{2 \xi -\sqrt{8 \xi +1}+1} \left(8 \xi -\sqrt{8 \xi
              +1}+1\right)   } \right)^{1/2} 
      \left(4 \xi +\sqrt{8 \xi +1}+8 \omega _0^2-1\right) \exp
      \left(- \frac{ \sqrt{2 \xi^2 - \xi \sqrt{8 \xi +1} +
            \xi} \left(\sqrt{8 \xi +1}-3\right)  }{ \sqrt{2}
          \left(\sqrt{8 \xi +1}-1\right)} t \right)}{4\ 2^{3/4}
      \sqrt{\pi } \left((4 \xi -1) \omega _0^2+(\xi -1) \xi +4
        \omega_0^4\right)}. 
    \label{eq:LBBM_xi>1_soln_2-intro}
  \end{equation}
  \end{small}
\end{enumerate}
\end{theorem}

The positive group velocity $c_g$ in Theorems \ref{lkdv-t} and
\ref{lbbm-t} is defined to be $c_g := \omega'(k_0)$, where $k_0$ is
the unique root of the dispersion relation $\omega(k)$ that lies on
$\p D^+$ defined in Theorem \ref{dn-t}. The ray $x/t=c_g$ for large
time $t \to \infty$ represents the leading edge of the wave's
transition from temporally periodic to temporally damped.  Theorems
\ref{lkdv-t} and \ref{lbbm-t} are established in Sections
\ref{kdv-asymp-s} and \ref{bbm-asymp-s}, respectively. Moreover, in
Section \ref{exp-s}, the analytical results of Section
\ref{bbm-asymp-s} are found to be in excellent agreement with
experimental observations of a viscous fluid conduit modeled by the
linear BBM equation in the long wavelength and small amplitude
regime. A lighter, less viscous fluid is injected into a tall
reservoir of fluid to establish a vertical fluid conduit rising
buoyantly. The spontaneous introduction of a sinusoidal injection rate
amounts to the wavemaker problem. As in previous experiments
\cite{mao2023experimental}, we observe a critical frequency below
which waves propagate and above which waves spatially decay. New
measurements include spatiotemporal wave patterns along rays
$x/t = \text{constant}$ that agree quantitatively with our large $t$
predictions.  These measurements include non-decaying waves, algebraically decaying
waves, and exponentially decaying waves. The transition between
non-decaying and algebraically decaying waves in the subcritical
regime is observed to occur when $x/t$ coincides with the group
velocity, in agreement with asymptotic predictions.

The universal nature of the Fokas method allows us to formulate a conjecture about the D-N map of the general nonlocal, third-order model \eqref{gen-mod}  in the case of asymptotically periodic Dirichlet data on the half-line. Although we do not prove this conjecture in the case of general coefficients $A_{-2}, A_0, A_1, A_2, A_3$, we anticipate that one should be able to establish it for any specific choice of these coefficients by suitably adapting the arguments leading to the proofs of Theorems \ref{dn-t}-\ref{lbbm-t}.
\\[3mm]
\textbf{Conjecture} \  
{\itshape Consider the general initial-boundary value problem
  \eqref{gen-ibvp} with an asymptotically periodic Dirichlet boundary
  condition
  $u(0, t) = g_0(t) \sim \sum_{n\in\mathbb Z} a_n e^{in\omega_0 t}$,
  $t\to\infty$, where $\omega_0\geq 0$, $a_n \in \mathbb C$. Then,
  under the assumption that the solution $u(x, t)$ along with its
  boundary values $\p_x^j u(0, t)$, $j\in\mathbb N$, are
  asymptotically time-periodic, the associated generalized D-N map
  satisfies the asymptotic formula of Theorem \ref{dn-t}. Furthermore,
  the uniqueness of the D-N map is equivalent to the radiation
  condition \eqref{eq:general_radiation_condition}. }

The paper is organized as follows.  Theorem \ref{dn-t} is proven,
under various assumptions, in Sec.~\ref{D2N-s}.  In
Sec.~\ref{sec:general-framework}, we construct the asymptotic D-N map
\eqref{gdn-map} for the general model \eqref{gen-mod}, which is used
to obtain the series representation \eqref{eq:asymptotic_soln} of the
asymptotic solution.  The asymptotic D-N map involves the root
$k_0(n)$ of the cubic polynomial $\omega(k_0(n)) = n\omega_0$ that is
proved to be the unique root lying on $\partial D^+$ obtained in the
restricted settings when $A_{-2} = 0$ (Sec.~\ref{third-ss}) or
$A_3 = 0$ (Sec.~\ref{2-nonlocal}), provided the frequency $n\omega_0$
lies outside an interval.  The proof for all real frequencies
$n\omega_0$ is obtained for the linear KdV (Sec.~\ref{third-ss}) and
linear BBM (Sec.~\ref{2-nonlocal}) equations where we also prove the
asymptotic periodicity of the solution.  Consequently, the radiation
condition \eqref{eq:general_radiation_condition} is proven for the
important physical cases of the linear KdV and linear BBM equations.
Theorems \ref{lkdv-t} and \ref{lbbm-t} are proven in
Secs.~\ref{kdv-asymp-s} and \ref{bbm-asymp-s}, respectively.  Viscous
core-annular fluid experiments that quantitatively realize Theorems
\ref{dn-t} and \ref{lbbm-t} are presented in Sec.~\ref{exp-s}.  We
conclude with a discussion of the implications of this work in
Sec.~\ref{sec:conclusions}.

\section{The D-N Map via the Fokas method}
\label{D2N-s}

As noted in the introduction, beyond its main purpose as a universal approach for the solution of linear evolution equations in the initial-boundary value problems setting, the Fokas method provides an effective and efficient way of determining the (generalized) D-N map associated with such problems when the boundary conditions are asymptotically periodic in time. This aspect of the method has been the subject of the recent articles by Fokas and van der Weele \cite{fv2021} and Fokas, Pelloni and Smith \cite{fps2022} for problems on the half-line and the finite interval, respectively. 

Here, we consider the Dirichlet problem on the half-line for the general third-order model \eqref{gen-mod}, which includes the linear KdV equation discussed in \cite{fv2021} but additionally involves a mixed spatiotemporal derivative like the one present in the linear BBM equation. Through our analysis, we demonstrate that the Fokas method for the characterization of the D-N map in the case of asymptotically time-periodic boundary conditions remains effective for this general third-order model. We note that the  analysis of this section is valid for generic initial data and not just for the zero initial condition stated in the wavemaker problem \eqref{gen-ibvp}.

Importantly, we address a uniqueness issue that arises for temporal frequencies below a certain threshold (the \textit{subcritical regime}) and manifests itself through the presence of solely real singularities among which only one specific choice gives rise to the correct D-N map, despite the fact that, at first sight, all of them qualify as removable.
In particular, by deriving the leading-order long-time asymptotics directly from the explicit solution formulae obtained via the Fokas method, we rigorously justify the empirical ``rule of thumb'' used previously in the literature, according to which, in the subcritical regime, the singularities that must be made removable are the ones along the boundary of the region $D^-$ defined by \eqref{dpm-gen-def} below.

\subsection{General framework}
\label{sec:general-framework}

Define the half-line Fourier transform of $u$ with respect to $x$ by
\begin{equation}\label{ft-hl}
\widehat u(k, t) := \int_{x=0}^\infty e^{-ikx} \, u(x, t) dx, \quad \text{Im}(k) \leq  0.
\end{equation}
We emphasize that, unlike the traditional whole-line Fourier transform, the half-line transform \eqref{ft-hl} is valid not only for $k\in\mathbb R$ but also for $\text{Im}(k) < 0$, due to the fact that $x \geq 0$.
Upon taking the transform \eqref{ft-hl} and integrating by parts, equation \eqref{gen-mod} yields
\begin{equation}\label{gen-gr-0}
\begin{aligned}
&\quad
\big(1 + A_{-2} \, k^2 \big) \, \widehat u_t (k, t)
+
i\big(A_3 k^3 + A_2 k^2 -A_1 k - A_0 \big) \, \widehat u(k, t)
\\
&= 
- A_3 \, g_2(t) -A_{-2} \, g_1'(t) - i \left(A_3 k + A_2\right)  g_1(t)   
- i A_{-2} k \, g_0'(t)  + \big(A_3 k^2 + A_2 k - A_1\big) \, g_0(t),
\end{aligned}
\end{equation}
where we have introduced the following notation for the various (known or unknown) boundary values of $u$ at $x=0$:
\begin{equation}\label{gj-def}
\p_x^j u(0, t) =:g_j(t), \quad j=0, 1, 2.
\end{equation}

Consider now the case of an asymptotically periodic Dirichlet boundary datum,  i.e. suppose that $g_0$ is given and satisfies
\begin{equation}\label{g0-per}
g_0(t) \sim \sum_{n\in\mathbb Z} a_n e^{in\omega_0 t}, \quad t\to\infty.
\end{equation}
Furthermore, suppose that the remaining (unknown) boundary values
$g_1$ and $g_2$, as well as the solution $u$, are also asymptotically
time-periodic with the same periodicity as $g_0$\footnote{This
  assumption can be justified by computing the long-time asymptotics
  of the solution obtained from integral asymptotics of the Fokas
  method solution --- see Sections \ref{lkdv-ss} and \ref{lbbm-ss} for  linear KdV and linear BBM respectively.}, i.e.
\begin{equation}
g_1(t) \sim \sum_{n\in\mathbb Z} b_n e^{in\omega_0 t},
\quad
g_2(t) \sim \sum_{n\in\mathbb Z} c_n e^{in\omega_0 t},
\quad
u(x, t) \sim \sum_{n\in\mathbb Z} U_n(x) e^{in\omega_0 t},
 \quad t\to\infty.
\end{equation}
Under these assumptions, the characterization of the D-N map for equation \eqref{gen-mod} on the half-line amounts to expressing the unknown coefficients $b_n$ and $c_n$ in terms of the  coefficients $a_n$, which are known since they correspond to the prescribed boundary datum.

Under the assumption that the half-line Fourier transform in $x$ can be interchanged with the infinite series over $n$ to infer asymptotic time periodicity for $\widehat u$, namely
\begin{equation}
\label{eq:4}  
\widehat u(k, t) \sim \sum_{n\in\mathbb Z} \widehat U_n(k) e^{in\omega_0 t}, \quad t\to\infty,
\end{equation}
and also that the various infinite series can be differentiated term by term,  equation \eqref{gen-gr-0} yields
\begin{equation}\label{gr-gen2}
\begin{aligned}
&\quad
\sum_{n\in\mathbb Z} i\left[\Omega + n\omega_0\big(1+A_{-2} \, k^2\big)\right] \widehat U_n(k) e^{in\omega_0 t}
\\
&\sim
\sum_{n\in\mathbb Z} 
\left\{
- A_3  c_n   - i \left(A_3 k + A_2 + n \omega_0 A_{-2}\right) b_n 
+ \big[A_3 k^2 + \left(A_2 + n \omega_0 A_{-2}\right) k - A_1 \big]  a_n  \right\} e^{in\omega_0 t},
\quad
t\to\infty,
\end{aligned}
 \end{equation}
where we have introduced the notation 
\begin{equation}\label{om-def-cubic}
\Omega = \Omega(k) := A_3 k^3 + A_2 k^2 -A_1 k - A_0.
\end{equation}
Consequently, $\omega(k) = \Omega(k)/(1+A_{-2}k^2)$ is the dispersion
relation \eqref{eq:5}.
We deduce from \eqref{gr-gen2}
\begin{equation}\label{Unk-gen}
\widehat U_n(k) =
\frac{ 
- A_3  c_n   - i \left(A_3 k + A_2 + n \omega_0 A_{-2}\right) b_n 
+ \big[A_3 k^2 + \left(A_2 + n \omega_0 A_{-2}\right) k - A_1 \big]  a_n}{i\left[\Omega+n\omega_0 \left(1 + A_{-2} \, k^2\right)\right]}.
\end{equation}
Zeros of the denominator correspond to roots of
$\omega(k) = -n\omega_0$.  Note that the negative sign is due to the
fact that the dispersion relation \eqref{eq:5} was defined for waves
in the form \eqref{eq:asymptotic_traveling_wave}, i.e., opposite that
of the frequency defined in \eqref{eq:4}.

As mentioned above, $\widehat u(k, t)$ is analytic for $\text{Im}(k) < 0$ and well-defined in the $L^2$-sense for $\text{Im}(k) \leq 0$ (e.g. see Theorem 7.2.4 in \cite{s1994}). Hence, it would be reasonable to expect that the expression \eqref{Unk-gen} for $\widehat U_n(k)$ is well-defined for $\text{Im}(k) \leq 0$. In turn, this would imply that any zeros of the quantity $\Omega+n\omega_0 \big(1 + A_{-2} \, k^2\big)$ that arise in $\left\{\text{Im}(k) \leq 0\right\}$ are actually removable singularities of \eqref{Unk-gen}. %
We shall see later, however, that this is true only for the zeros of $\Omega+n\omega_0 \big(1 + A_{-2} \, k^2\big)$ that lie on the \textit{boundary} $\p D^-$  of the domain $D^-$, where
\begin{equation}\label{dpm-gen-def}
D^\pm = \left\{k\in\mathbb C: \text{Im}(k) \gtrless 0 \text{ and } \text{Re}\left(\frac{i\Omega}{ 1+A_{-2} k^2}\right) < 0 \right\}.
\end{equation}
Observe that, since $n\omega_0 \in \mathbb R$, the zeros of $\Omega + n \omega_0\big(1+A_{-2} k^2\big)$ always satisfy the equation $\text{Re}\left(\frac{i\Omega}{ 1+A_{-2} k^2}\right) = 0$ and hence lie on $\p D^+ \cup \p D^-$, as $\mathbb R \subset \left(\p D^+ \cup \p D^-\right)$.
To see this, note that
\begin{equation}\label{re-iw}
\begin{aligned}
\text{Re}\left(\frac{i\Omega}{ 1+A_{-2} k^2}\right) 
&=
\frac{\text{Im}(k)}{\left|1 + A_{-2} \, k^2\right|^2}
\Big[
A_1 - 2 \left(A_2 + A_0 A_{-2}\right) \text{Re}(k) - \left(3 A_3 + A_1 A_{-2}\right) \text{Re}(k)^2 - A_{-2} A_3 \text{Re}(k)^4 
\\
&\hskip 2cm
- \left(A_1 A_{-2} - A_3 + 2 A_{-2} A_3 \text{Re}(k)^2\right) \text{Im}(k)^2 - A_{-2} A_3 \text{Im}(k)^4
\Big]
\end{aligned}
\end{equation}
and so $\text{Re}\left(\frac{i\Omega}{ 1+A_{-2} k^2}\right)$ vanishes for $k\in\mathbb R$ and   is odd with respect to $\text{Im}(k)$. Thus,  $\mathbb R \subset \left(\p D^+ \cup \p D^-\right)$ since  $\left(\p D^+ \cup \p D^-\right) \cap \mathbb R \neq \emptyset$ and  the conjugate of any point in $D^+$ belongs to $\left\{\text{Im}(k) < 0\right\} \setminus D^-$ and the conjugate of any point in $D^-$ belongs to $\left\{\text{Im}(k) > 0\right\} \setminus D^+$. 

\begin{remark}
The departure from the original intuition that the expression \eqref{Unk-gen} for $\widehat U_n(k)$ should make sense for all $\text{Im}(k) \leq 0$  is due to the fact that, although related to each other, the quantities $\widehat U_n(k)$ and $\widehat u(k, t)$ are not quite the same; it seems that the  manipulations between infinite series and integrals required in order to derive \eqref{Unk-gen} do not fully preserve  analyticity in $k$ from $\widehat u(k, t)$ to $\widehat U_n(k)$. 
\end{remark}

Explicitly, $\Omega + n \omega_0 \big(1+A_{-2} k^2\big)=0$ corresponds to the cubic equation
\begin{equation}\label{poles-eq}
A_3 k^3 + \left(A_2 + n\omega_0 A_{-2}\right) k^2 -A_1 k - A_0 + n\omega_0  = 0.
\end{equation}
The requirement that the zeros of \eqref{poles-eq} along $\p D^-$ are removable singularities of the right-hand side of \eqref{Unk-gen} implies that these zeros must also satisfy  the equation
\begin{equation}\label{dn-gen-0}
- A_3  c_n   - i \left(A_3 k + A_2 + n \omega_0 A_{-2}\right) b_n 
+ \big[A_3 k^2 + \left(A_2 + n \omega_0 A_{-2}\right) k - A_1 \big]  a_n = 0.
\end{equation}
Then, for $A_3 \neq 0$, \textit{assuming that precisely two roots $k_1^- = k_1^-(n)$ and $k_2^- = k_2^-(n)$ of the cubic equation \eqref{poles-eq} lie on $\p D^-$ while the third root $k_0 = k_0(n)$ lies on $\p D^+$}, we have the following system for the coefficients $b_n$ and $c_n$:
\begin{equation}\label{dn-gen-sys}
\begin{aligned}
&- A_3  c_n   - i \left(A_3 k_1^- + A_2 + n \omega_0 A_{-2}\right) b_n 
+ \big[A_3 (k_1^-)^2 + \left(A_2 + n \omega_0 A_{-2}\right) k_1^- - A_1 \big] a_n = 0,
\\
&- A_3  c_n   - i \left(A_3 k_2^- + A_2 + n \omega_0 A_{-2}\right) b_n 
+ \big[A_3 (k_2^-)^2 + \left(A_2 + n \omega_0 A_{-2}\right) k_2^- - A_1 \big]  a_n = 0.
\end{aligned}
\end{equation}
Solving this system while noting that, by Vieta's formulae, 
$k_1^- + k_2^- + k_0 = - \left(A_2 + n\omega_0 A_{-2}\right)/A_3$ and $k_1^- k_2^- + k_2^- k_0 + k_0 k_1^- = -A_1/A_3$, we obtain the D-N map for the third-order model \eqref{gen-mod} as
\begin{equation}\label{dn-gen}
b_n = i k_0(n) \, a_n,
\quad
c_n = -\left[k_0(n)\right]^2 \, a_n, 
\quad
n \in \mathbb Z.
\end{equation}
Similarly, for $A_3=0$, which corresponds to the special case where  \eqref{gen-mod} reduces to a second-order model, assuming that precisely one root $k_1^-=k_1^-(n)$  of the now quadratic equation \eqref{poles-eq} lies on $\p D^-$ while the second root $k_0=k_0(n)$ lies on $\p D^+$, we obtain the first of the relations \eqref{dn-gen}, which amounts to the D-N map for this second-order reduction.

\begin{remark} \textbf{(Series representation of the leading-order asymptotics via the D-N map).} 
Under the assumption of sufficient smoothness, the leading-order long-time asymptotics of the solution of the general model \eqref{gen-mod} on the half-line can be directly constructed via the D-N map \eqref{dn-gen} by means of the Maclaurin series  $u(x,t) = \sum_{j=0}^\infty \frac{\p_x^j u(0,t)}{j!} x^j$. 
Indeed, for $A_3\neq0$ \eqref{dn-gen} yields 
$$
u_x(0,t)\sim i \sum_{n\in \mathbb{Z}} k_0(n) a_n e^{in\omega_0 t}, \quad u_{xx}(0,t) \sim -\sum_{n\in \mathbb{Z}} [k_0(n)]^2 a_n e^{in\omega_0 t},
$$
which can be combined with \eqref{gen-mod} evaluated at the boundary $x=0$ and the fact that $k_0(n)$ satisfies \eqref{poles-eq} to obtain
$$
u_{xxx}(0,t) \sim -i \sum_{n\in \mathbb{Z}} [k_0(n)]^3 a_n e^{in\omega_0 t}.
$$
Moreover, along the same lines and using induction while assuming that the  $x$-derivatives of the solution are asymptotically time-periodic at $x=0$, we have the general formula
\begin{equation}
\p_x^j u(0,t) \sim i^j \sum_{n\in \mathbb{Z}} [k_0(n)]^j a_n e^{in\omega_0t}, \quad j\in\mathbb N \cup \{0\}. \label{eq:series-um}
\end{equation} 
Similarly, for $A_3=0$ we can use the first of the relations \eqref{dn-gen} to establish \eqref{eq:series-um} via induction. 
Therefore, in view of the Maclaurin series of $u$ with respect to $x$, and assuming that we can interchange the order of summation, we deduce
\begin{equation}
u(x,t) \sim \sum_{j=0}^\infty \frac{i^j x^j}{j!} \sum_{n\in \mathbb{Z}} [k_0(n)]^j a_n e^{in\omega_0t} 
=
\sum_{n\in \mathbb{Z}} a_n e^{i k_0(n)x + i n\omega_0t}, 
\quad t\to\infty,
\label{eq:LKdV_series_soln}
\end{equation}
which shows that, for asymptotically periodic boundary data that are
otherwise arbitrary, the leading-order asymptotics of the solution of
\eqref{gen-mod} is a superposition of waves that either propagate (if
$k_0(n) \in \mathbb{R}$) or spatially decay (if
$\mathrm{Im}\, k_0(n) > 0$).
\end{remark}

In order to rigorously justify the D-N map \eqref{dn-gen}, one must
study the location of the roots of the cubic equation~\eqref{poles-eq}
with respect to the curves $\p D^-$ and $\p D^+$. For general
parameters $A_j$, this task involves the study of the quartic
polynomial in $\text{Re}(k)$ and $\text{Im}(k)$ that appears in the
expression \eqref{re-iw} for
$\text{Re}\left(\frac{i\Omega}{ 1+A_{-2} k^2}\right)$. Hence, below we
focus on two particular sets of parameters, namely
$\{A_{-2} = 0, A_3<0\}$ and $\{A_3 = 0, A_{-2}>0\}$, which still
correspond to generalizations of the linear KdV and the linear BBM
equations, respectively, that are more tractable due to the fact that
the aforementioned quartic reduces to a quadratic.

\subsection{A general local third-order model: $A_{-2} = 0$ and $A_3<0$}
\label{third-ss}

For $A_{-2} = 0$ and $A_3 < 0$, equation \eqref{gen-mod} no longer involves a mixed spatiotemporal derivative as it takes the form
\begin{equation}\label{gen-mod-a-2=0}
u_t = i A_0 u + A_1 u_x + i A_2 u_{xx} + A_3 u_{xxx}
\end{equation}
and the denominator in \eqref{Unk-gen} reduces to $i(\Omega + n\omega_0)$ with equation \eqref{poles-eq} simplifying to
\begin{equation}\label{cubic-a-2=0}
A_3 k^3 + A_2 k^2 -A_1 k - A_0 + n\omega_0 = 0.
\end{equation}

Before discussing the roots of \eqref{cubic-a-2=0}, it is useful to give the precise characterization of the domains $D^-$ and $D^+$. Starting from \eqref{om-def-cubic}, we compute
\begin{equation}
\text{Re}\left(\frac{i\Omega}{ 1+A_{-2} k^2}\right) 
=
\text{Re}\left( i\Omega \right)
=
\text{Im}(k) \big[A_1 - 2 A_2 \text{Re}(k) - 3 A_3 \text{Re}(k)^2 + A_3 \text{Im}(k)^2\big]
\end{equation}
so that, according to the definition \eqref{dpm-gen-def},
\begin{equation}\label{d-ineq}
D^-
=
\left\{
k \in \mathbb C:
\text{Im}(k) < 0
\text{ and }
A_1 - 2 A_2  \text{Re}(k) -  3 A_3 \text{Re}(k)^2  + A_3 \text{Im}(k)^2 > 0
\right\}.
\end{equation}
Thus, for real solutions of \eqref{cubic-a-2=0} that lie along the boundary $\p D^-$  we must have
\begin{equation}\label{real-cond}
3 A_3 k^2 + 2 A_2  k - A_1  \leq  0.
\end{equation}
Similarly, the region $D^+$ is described by 
\begin{equation}\label{d+ineq}
D^+
=
\left\{
k \in \mathbb C:
\text{Im}(k) > 0
\text{ and }
A_1 - 2 A_2  \text{Re}(k) -  3 A_3 \text{Re}(k)^2  + A_3 \text{Im}(k)^2 < 0
\right\}
\end{equation}
 and, therefore, any real solution of \eqref{cubic-a-2=0} that lies along the boundary $\p D^+$ must satisfy
\begin{equation}\label{real-cond+}
3 A_3 k^2 + 2 A_2 k - A_1  \geq  0.
\end{equation}

We now turn our attention to the solution of equation \eqref{cubic-a-2=0}. The roots of this cubic  are given by
\begin{equation}\label{kj-nj}
k_j = \nu_j - \frac{A_2}{3A_3}, \quad j=1,2,3,
\end{equation}
where $\nu_j$ are the solutions of the depressed cubic 
\begin{equation}\label{dc-eq}
\nu^3 + p \nu + q = 0 
\ \ \text{ with } \ \ 
p = - \frac{3A_1A_3+A_2^2}{3A_3^2},
\ 
q = 
\frac{2A_2^3+9A_1A_2A_3-27A_3^2\left(A_0 - n\omega_0\right)}{27A_3^3}.
\end{equation}
The discriminant of this depressed cubic is $\Delta = - 4p^3 - 27q^2$. According to the classical result by del Ferro and Tartaglia, also known as Cardano's formula, there are three cases:
\\[2mm]
\textbf{Case 1:} $\Delta<0$. This is equivalent to $\frac{q^2}{4} + \frac{p^3}{27} > 0$. 
Note that this inequality is satisfied for all $n\in\mathbb Z$ if $p>0$. On the other hand, if $p \leq  0$  then taking square roots and defining
\begin{equation}\label{spm-def}
\omega_{cr}^\pm := \overline{\omega} \mp  2 A_3  \sqrt{\left(-\frac p3\right)^3}
\ \ \text{ where } \ \ 
\overline{\omega} 
:= 
\frac{-2A_2^3 - 9A_1A_2A_3 + 27A_3^2 A_0}{27A_3^2}
\end{equation}
we must have either $n\omega_0 > \omega_{cr}^+$ or $n\omega_0 < \omega_{cr}^-$ (note that $\omega_{cr}^- \leq \overline{\omega} \leq  \omega_{cr}^+$ since $A_3<0$ and $p\leq  0$).
That is, $\Delta<0$ if either  $p>0$, or  $p \leq  0$ and $n\omega_0>\omega_{cr}^+$ or $n\omega_0<\omega_{cr}^-$.
Then, the depressed cubic \eqref{dc-eq} has  the real solution  
$
\nu_1
=
\Big(-\frac q2 + \sqrt{\frac{q^2}{4} + \frac{p^3}{27}} \Big)^{\frac 13} 
+ 
\Big(-\frac q2 - \sqrt{\frac{q^2}{4} + \frac{p^3}{27}} \Big)^{\frac 13} 
$
and the complex conjugate solutions 
$\nu_2
= 
\Big(-\frac q2 + \sqrt{\frac{q^2}{4} + \frac{p^3}{27}} \Big)^{\frac 13}  e^{i\frac{2\pi}{3}}
+ 
\Big(-\frac q2 - \sqrt{\frac{q^2}{4} + \frac{p^3}{27}} \Big)^{\frac 13}  e^{-i\frac{2\pi}{3}}
$
and
$
\nu_3
=
\Big(-\frac q2 + \sqrt{\frac{q^2}{4} + \frac{p^3}{27}} \Big)^{\frac 13}  e^{-i\frac{2\pi}{3}}
+ 
\Big(-\frac q2 - \sqrt{\frac{q^2}{4} + \frac{p^3}{27}} \Big)^{\frac 13}  e^{i\frac{2\pi}{3}}$.
In turn,  the cubic \eqref{cubic-a-2=0} has one real root and two complex conjugate roots, namely
\begin{equation}\label{d-sols}
\begin{split}
&k_1
=
\left(-\frac q2 + \sqrt{\frac{q^2}{4} + \frac{p^3}{27}} \, \right)^{\frac 13} 
+ 
\left(-\frac q2 - \sqrt{\frac{q^2}{4} + \frac{p^3}{27}} \, \right)^{\frac 13} - \frac{A_2}{3A_3},
\\
&k_2
= 
\left(-\frac q2 + \sqrt{\frac{q^2}{4} + \frac{p^3}{27}} \, \right)^{\frac 13} e^{i\frac{2\pi}{3}}
+ 
\left(-\frac q2 - \sqrt{\frac{q^2}{4} + \frac{p^3}{27}} \, \right)^{\frac 13} e^{-i\frac{2\pi}{3}} - \frac{A_2}{3A_3},
\\
&k_3
=
\left(-\frac q2 + \sqrt{\frac{q^2}{4} + \frac{p^3}{27}} \, \right)^{\frac 13} e^{-i\frac{2\pi}{3}}
+ 
\left(-\frac q2 - \sqrt{\frac{q^2}{4} + \frac{p^3}{27}} \, \right)^{\frac 13} e^{i\frac{2\pi}{3}} - \frac{A_2}{3A_3}. 
\end{split}
\end{equation}
Recalling that $\text{Re}(i\Omega)  = 0$ at the roots $k_j$, it follows that precisely one of the complex conjugate roots $k_2$ and $k_3$ lies on $\p D^-$. Moreover, we claim that $k_1 \in \p D^-$. Indeed, we begin by observing that $\nu_1^2 \geq  -\frac p3$. This is clearly true if $p\geq  0$. In addition, if $p<0$ then, noting that for $q \geq 0$ we have $\frac q2 - \sqrt{\frac{q^2}{4} + \frac{p^3}{27}} \geq 0$ while for $q < 0$ we have $-\frac q2 - \sqrt{\frac{q^2}{4} + \frac{p^3}{27}} \geq 0$, we infer that $-\nu_1\geq \left(\frac q2\right)^{1/3}$ for $q\geq 0$ while $\nu_1\geq \left(-\frac q2\right)^{1/3}$ for $q< 0$. Hence, if $p<0$ then $|\nu_1|\geq  \left(\frac{|q|}{2}\right)^{1/3}$ which in view of the fact that $\Delta<0$ implies $\nu_1^2 \geq  -\frac p3$ as desired.
In turn, using \eqref{kj-nj} and substituting for $p$ while recalling that $A_3<0$, we obtain precisely the condition \eqref{real-cond} for $k_1$ and thus prove that $k_1 \in \p D^-$ as claimed.
\\[2mm]
\textbf{Case 2:} $\Delta>0$. This amounts to $\frac{q^2}{4} + \frac{p^3}{27} < 0$, which yields $\omega_{cr}^- < n\omega_0 < \omega_{cr}^+$ with $\sigma_\pm$  given by \eqref{spm-def} (note that $\omega_{cr}^+>\omega_{cr}^-$ because $A_3<0$ and for $\Delta>0$ we must have $p<0$).
Then, the depressed cubic~\eqref{dc-eq} has three distinct real roots, 
$
\nu_1
=
\Big(-\frac q2 + i \sqrt{-\frac{q^2}{4} - \frac{p^3}{27}} \Big)^{\frac 13} 
-
\frac p3 \Big(-\frac q2 + i \sqrt{-\frac{q^2}{4} - \frac{p^3}{27}} \Big)^{-\frac 13}
$,
$\nu_2
=
\Big(-\frac q2 + i \sqrt{-\frac{q^2}{4} - \frac{p^3}{27}} \Big)^{\frac 13} e^{i\frac{2\pi}{3}} 
-
\frac p3 \Big(-\frac q2 + i \sqrt{-\frac{q^2}{4} - \frac{p^3}{27}} \Big)^{-\frac 13} e^{-i\frac{2\pi}{3}}
$
and
$\nu_3
=
\Big(-\frac q2 + i \sqrt{-\frac{q^2}{4} - \frac{p^3}{27}} \Big)^{\frac 13} e^{-i\frac{2\pi}{3}}
-
\frac p3 \Big(-\frac q2 + i \sqrt{-\frac{q^2}{4} - \frac{p^3}{27}} \Big)^{-\frac 13} e^{i\frac{2\pi}{3}}
$,
which correspond to the following distinct real roots for the  cubic \eqref{cubic-a-2=0}: 
\begin{equation}\label{d+sols}
\begin{split}
&k_1
=
\left(-\frac q2 + i \sqrt{-\frac{q^2}{4} - \frac{p^3}{27}} \, \right)^{\frac 13} 
-
\frac p3 \left(-\frac q2 + i \sqrt{-\frac{q^2}{4} - \frac{p^3}{27}} \, \right)^{-\frac 13} - \frac{A_2}{3A_3},
\\
&k_2
=
\left(-\frac q2 + i \sqrt{-\frac{q^2}{4} - \frac{p^3}{27}} \, \right)^{\frac 13} e^{i\frac{2\pi}{3}} 
-
\frac p3 \left(-\frac q2 + i \sqrt{-\frac{q^2}{4} - \frac{p^3}{27}} \, \right)^{-\frac 13} e^{-i\frac{2\pi}{3}} - \frac{A_2}{3A_3},
\\
&k_3
=
\left(-\frac q2 + i \sqrt{-\frac{q^2}{4} - \frac{p^3}{27}} \, \right)^{\frac 13} e^{-i\frac{2\pi}{3}}
-
\frac p3 \left(-\frac q2 + i \sqrt{-\frac{q^2}{4} - \frac{p^3}{27}} \, \right)^{-\frac 13} e^{i\frac{2\pi}{3}} - \frac{A_2}{3A_3}.
\end{split}
\end{equation}
We will now show that $k_1, k_2 \in \p D^-$ while $k_3 \in \p D^+$. Writing
$
-\frac q2 + i \sqrt{-\frac{q^2}{4} - \frac{p^3}{27}} = r e^{i\theta}
$
with $r = \sqrt{-\frac{p^3}{27}}$ and $\tan \theta = \sqrt{-\frac{q^2}{4} - \frac{p^3}{27}}\Big /(-\frac q2)$, we have
\begin{equation}
\begin{split}
&k_1 
= r^{\frac 13} e^{i\frac \theta 3} - \frac p3 \, r^{-\frac 13} e^{-i\frac \theta 3} - \frac{A_2}{3A_3}
=
2 \sqrt{-\frac p3} \cos \left(\frac \theta 3\right) - \frac{A_2}{3A_3},
\\
&k_2 =  r^{\frac 13} e^{i\frac \theta 3} e^{i\frac{2\pi}{3}} - \frac p3 \, r^{-\frac 13} e^{-i\frac \theta 3} e^{-i\frac{2\pi}{3}} - \frac{A_2}{3A_3}
=
2 \sqrt{-\frac p3} \cos \left(\frac{\theta+2\pi} 3\right) - \frac{A_2}{3A_3},
\\
&k_3 =  r^{\frac 13} e^{i\frac \theta 3} e^{-i\frac{2\pi}{3}} - \frac p3 \, r^{-\frac 13} e^{-i\frac \theta 3} e^{i\frac{2\pi}{3}} - \frac{A_2}{3A_3}
=
2 \sqrt{-\frac p3}  \cos \left(\frac{\theta-2\pi} 3\right) - \frac{A_2}{3A_3}.
\end{split}
\end{equation}
Since $\omega_{cr}^- < \overline{\omega} < \omega_{cr}^+$ for $p<0$ and, in addition, $\overline{\omega} = - A_3 q  + n\omega_0$ and $A_3<0$, we have two subcases:
\\[2mm]
(i) If $\overline{\omega} \leq   n\omega_0 < \omega_{cr}^+$, then $q\leq  0$ and so $\theta \in \left(0, \frac \pi 2\right]$. Hence, 
$\frac{\sqrt 3}{2} \leq  \cos \frac \theta 3 < 1$, $-\frac{\sqrt 3}{2} \leq  \cos \frac{\theta+2\pi} 3 < -\frac 12$  and $-\frac 12  < \cos \frac{\theta-2\pi} 3 \leq  0$,
which imply  $\sqrt{-p} - \frac{A_2}{3A_3} \leq  k_1 < 2 \sqrt{-\frac p3} - \frac{A_2}{3A_3}$, $-\sqrt{-p} - \frac{A_2}{3A_3} \leq  k_2 < -\sqrt{-\frac p3} - \frac{A_2}{3A_3}$ and $-\sqrt{-\frac p3} - \frac{A_2}{3A_3} < k_3 \leq  - \frac{A_2}{3A_3}$. Thus, since by \eqref{real-cond} the portion of  $\p D^-$ along $\mathbb R$ is given by $\big(-\infty, -\sqrt{-\frac p3}  - \frac{A_2}{3A_3}\big]
\cup
\big[\sqrt{-\frac p3}-\frac{A_2}{3A_3}, \infty\big)$, 
we deduce that $k_1, k_2 \in \p D^-$ and $k_3 \in \p D^+$. In particular, inequalities \eqref{real-cond} and \eqref{real-cond+} hold strictly at $k_1, k_2$ and at $k_3$, respectively.
\\[2mm]
(ii) If $\omega_{cr}^- <  n\omega_0 \leq  \overline{\omega}$, then $q\geq  0$ and so $\theta \in \left[\frac \pi 2, \pi\right)$.
Hence, $\frac 12 < \cos \frac \theta 3 \leq  \frac{\sqrt 3}{2}$, $-1< \cos \frac{\theta+2\pi} 3 \leq  -\frac{\sqrt 3}{2}$ and $0 \leq  \cos \frac{\theta-2\pi} 3 < \frac 12$, which imply $\sqrt{-\frac p3} - \frac{A_2}{3A_3} < k_1 \leq  \sqrt{-p} - \frac{A_2}{3A_3}$, $-2\sqrt{-\frac p3} - \frac{A_2}{3A_3} < k_2 \leq  -\sqrt{-p} - \frac{A_2}{3A_3}$ and $- \frac{A_2}{3A_3}\leq  k_3 < \sqrt{-\frac p3} - \frac{A_2}{3A_3}$. Thus, as in the previous subcase, we have $k_1, k_2 \in \p D^-$ and $k_3 \in \p D^+$ with inequalities \eqref{real-cond} and \eqref{real-cond+} holding strictly at $k_1, k_2$ and at $k_3$, respectively.
\\[2mm]
\textbf{Case 3:} $\Delta=0$. This is equivalent to $4p^3 + 27q^2 = 0$,
thus either $n \omega_0= \omega_{cr}^+$ or $n\omega_0 =
\omega_{cr}^-$.  In the former case, the cubic \eqref{cubic-a-2=0} admits the real solutions 
\begin{equation}\label{d0sols+}
n\omega_0 = \omega_{cr}^+ \quad \Rightarrow \quad k_1 
= -2\sqrt{-\frac p3} - \frac{A_2}{3A_3}, 
\quad 
k_2 = k_3 
= \sqrt{-\frac p3} - \frac{A_2}{3A_3} . 
\end{equation}
When $n\omega_0 = \omega_{cr}^-$,  the cubic \eqref{cubic-a-2=0} admits the real solutions 
\begin{equation}\label{d0sols-}
n\omega_0 = \omega_{cr}^- \quad \Rightarrow \quad k_1 
= 2\sqrt{-\frac p3} - \frac{A_2}{3A_3}, 
\quad 
k_2 = k_3 
= -\sqrt{-\frac p3} - \frac{A_2}{3A_3} . 
\end{equation}
Hence, in view of the fact that $\p D^-\cap \mathbb R = \big(-\infty, -\sqrt{-\frac p3}  - \frac{A_2}{3A_3}\big]
\cup
\big[\sqrt{-\frac p3}-\frac{A_2}{3A_3}, \infty\big)$, we deduce that  if $p<0$ then $k_1 \in \p D^-$ with strict inequality in \eqref{real-cond} and $k_2=k_3 \in \p D^-\cap \p D^+$ with equality in \eqref{real-cond+}, while if $p=0$ then $k_1 = k_2 = k_3 \in \p D^-\cap \p D^+$ with equality in \eqref{real-cond} and \eqref{real-cond+}.

\vspace*{2mm}
Overall, our analysis  leads to the following conclusions:
\begin{enumerate}[label=$\bullet$, leftmargin=3mm, topsep=0mm, itemsep=1mm]
\item If $3A_1A_3 + A_2^2 < 0$, then the solutions of \eqref{cubic-a-2=0} are given by \eqref{d-sols} for all $n \in \mathbb Z$.
\item If $3A_1A_3 + A_2^2 > 0$, then the solutions of
  \eqref{cubic-a-2=0} are given by \eqref{d-sols} for
  $n\omega_0>\omega_{cr}^+$ or $n\omega_0<\omega_{cr}^-$, by
  \eqref{d+sols} for $\omega_{cr}^- < n\omega_0 < \omega_{cr}^+$, and
  by \eqref{d0sols+} for $n\omega_0 = \omega_{cr}^+$ or \eqref{d0sols-} for $n\omega_0 = \omega_{cr}^-$, where $\omega_{cr}^\pm$ are defined by \eqref{spm-def}. 
\item If $3A_1A_3 + A_2^2 = 0$, then $\omega_{cr}^+=\omega_{cr}^-=\overline{\omega}$ defined by \eqref{spm-def} and the solutions of \eqref{cubic-a-2=0} are given by \eqref{d-sols} for $n\omega_0 \neq \overline{\omega}$  and by \eqref{d0sols+} for $n\omega_0 = \overline{\omega}$. 
\end{enumerate}
Importantly, in all three of the above cases, the contour $\p D^-$
contains precisely two solutions of \eqref{cubic-a-2=0} whereas the
contour $\p D^+$ contains only one, which we label $k_0$.
That is, we have shown that for $A_{-2}=0$ and $A_3<0$ is
satisfied, and hence the D-N map for the general third-order model
\eqref{gen-mod-a-2=0} is indeed given by \eqref{dn-gen}.

This proves Theorem \ref{dn-t} when $A_{-2} = 0$.

\begin{remark}
\textbf{(connection to the radiation condition \eqref{eq:general_radiation_condition})}
The group velocity associated with the third-order equation
\eqref{gen-mod-a-2=0} is  $c_g = \Omega'(k) = 3A_3 k^2 + 2A_2 k - A_1$.
Thus, for a real solution of the cubic \eqref{cubic-a-2=0}, the
condition that this solution lies on $\p D^+$ (equation \eqref{real-cond+}) is equivalent to  the group velocity being non-negative at that solution. 
In particular, recall that the scenario $k_0 \in \p D^+ \cap \mathbb
R$ is possible if either $3A_1A_3+A_2^2 > 0$ and $\omega_{cr}^- <
n\omega_0 < \omega_{cr}^+$, or $3A_1A_3+A_2^2 > 0$ and either
$n\omega_0 = \omega_{cr}^+$ or $n\omega_0 = \omega_{cr}^-$, or $3A_1A_3+A_2^2 = 0$ and $n\omega_0=\overline{\omega}$. Hence, we have a strictly positive group velocity in the first scenario (see ``Case 2'' above) and zero group velocity in the other three scenarios (see ``Case 3'' above), which deals with the boundary case $c_g = 0$.
Furthermore, if $k_0 \in \p D^+ \setminus \mathbb R$ then it follows from the definition \eqref{d+ineq} of $D^+$   that $\text{Im}(k)>0$.

This proves Corollary \ref{rad-cond} when $A_{-2} = 0$.
\end{remark}

\subsubsection*{The linear K\textnormal{d}V equation on the half-line}
\label{lkdv-ss}
The analysis of the roots of the cubic $\Omega + n\omega_0$ (i.e. of equation \eqref{cubic-a-2=0}) for the general third-order model \eqref{gen-mod-a-2=0} provided a complete justification of the assumption used in the derivation of the D-N map \eqref{dn-gen} that precisely two of those three roots lie on the contour $\p D^-$. However, there is one more assumption that was needed for deducing \eqref{dn-gen} and is yet to be justified, namely that the right-hand side of \eqref{Unk-gen} should have removable singularities only at those roots of $\Omega + n\omega_0$ that occur along $\p D^-$. 
Above a certain threshold of values for $|n\omega_0|$, it is easy to see why this latter assumption must be imposed; below that threshold, however, the situation becomes more delicate.  

More precisely, for those values of $n\omega_0$ for which the roots of
$\Omega + n\omega_0$ involve a complex conjugate pair (i.e. in Case~1:
$\Delta<0$ above, $n\omega_0 \in [\omega_{cr}^-,\omega_{cr}^+]^C$),
one of the complex roots lies in the upper half of the $k$-plane where
$\widehat u(k, t)$ (and hence $\widehat U_n(k)$) is not well-defined
(recall discussion below \eqref{ft-hl}). Thus, this root cannot be a
removable singularity of \eqref{Unk-gen}. On the other hand, the
remaining two roots do lie in the domain of definition
$\left\{\text{Im}(k)\leq 0\right\}$ of $\widehat u(k, t)$ and so they
must be removable singularities of \eqref{Unk-gen}, thereby leading to
the system \eqref{dn-gen-sys} and hence to the D-N map
\eqref{dn-gen}. The same is true in the degenerate case of one
distinct and two repeated real roots (i.e. in Case 3: $\Delta=0$
above, $n\omega_0 \in \{\omega_{cr}^-,\omega_{cr}^+\}$) since,
counting without multiplicities, there are exactly two roots that
belong to $\left\{\text{Im}(k)\leq 0\right\}$ and hence must be
removable singularities of \eqref{Unk-gen}.

Nevertheless, in the case of three distinct real roots (i.e. Case 2:
$\Delta>0$ above, $n\omega_0 \in (\omega_{cr}^-,\omega_{cr}^+)$), we
encounter the following issue: while all three roots, namely the two
lying on $\p D^-$ and the one lying on $\p D^+$, belong to
$\left\{\text{Im}(k)\leq 0\right\}$ and hence, according to our
previous reasoning, should be removable singularities of
\eqref{Unk-gen}, requiring that all three of them are removable
results via \eqref{dn-gen-0} in three equations for the two unknowns
$b_n$ and $c_n$, i.e. an overdetermined system. Therefore, the
question is: \textit{which two of the three distinct real roots of
  $\Omega + n\omega_0$ are actually removable singularities of
  \eqref{dn-gen-0}?}
Note that each of the three possible choices leads to a different D-N map since, via Vieta's formulae, the root which is not a removable singularity is the one that eventually appears in the D-N relations \eqref{dn-gen}. From this point of view, the challenge we are facing is directly related to the uniqueness of the D-N map.  

Our conjecture is that, even in the case of three distinct real roots, our earlier assumption that led to \eqref{dn-gen} is correct, namely the roots of $\Omega + n\omega_0$ that must be removable singularities of \eqref{Unk-gen} are precisely the two roots that lie along $\p D^-$ (and not the one along $\p D^+$ --- see remark below \eqref{dpm-gen-def} for a possible explanation of this fact in connection with the domains of definition of $\widehat u(k, t)$ and $\widehat U_n(k)$). In what follows, we prove this conjecture for the special case of the linear KdV equation  by exploiting the explicit solution formula produced by the Fokas method.

\begin{figure}[ht!] 
\centering 
\includegraphics[scale=0.45]{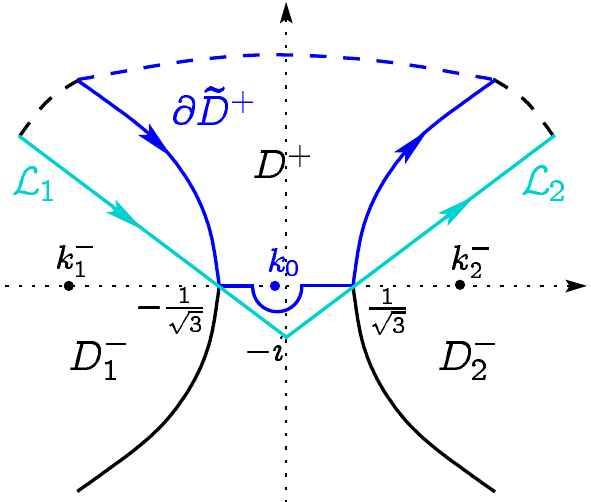}
\caption{The domains $D^+$, $D_1^-$, $D_2^-$ for the linear KdV equation as given by \eqref{dpm-lkdv}-\eqref{dm12-def}. Inside these domains $\text{Re}(i\Omega) < 0$, while outside of them $\text{Re}(i\Omega) > 0$. 
In the case of the wavemaker data \eqref{lkdv-sine-bc}, if $|\omega_0|<2/(3\sqrt 3)$ then the three roots of $\Omega + n\omega_0$ are real and distinct with $k_1^- \in \p D_1^-$, $k_2^- \in \p D_2^-$, $k_0\in \p D^+$ as shown, and the contour of integration in the Fokas method solution \eqref{kdv-sol-wm-x} is deformed from the boundary $\p D^+$ of $D^+$ to the contour $\p \widetilde D^+$ that bypasses the real root $k_0 \in \p D^+$ from below. The resulting formula \eqref{lkdv-sol-def} is used for deriving the D-N map of the wavemaker problem by computing the long-time asymptotics via Cauchy's residue theorem, exponential decay and, when necessary, the deformation of $\p \widetilde D^+$ to the half-lines $\mathcal
L_1$, $\mathcal L_2$ given by \eqref{L1L2-def}, along which
$\text{Re}(i\Omega) > 0$ (except for $k=\pm \frac{1}{\sqrt 3}$ where
$\text{Re}(i\Omega)=0$).}
\label{lkdv-f}
\end{figure}

The linear KdV equation
\begin{equation}\label{lkdv}
u_t + u_x + u_{xxx} = 0
\end{equation}
is included in the four-parameter family of third-order models \eqref{gen-mod-a-2=0} (it corresponds to $A_0=A_2=0$ and $A_1=A_3=-1$), thus the general framework developed earlier can readily be employed. In particular, the dispersion relation is $\Omega(k) = k-k^3$, the regions $D^+$ and $D^-$ are given by
\begin{equation}\label{dpm-lkdv}
D^+ = \left\{k \in \mathbb C: \text{Im}(k) > 0, \ 3\text{Re}(k)^2 - \text{Im}(k)^2 < 1 \right\},
\quad
D^- = D_1^- \cup D_2^-,
\end{equation}
where (see Figure \ref{lkdv-f})
\begin{equation}\label{dm12-def}
\begin{aligned}
&D_1^- = \left\{k \in \mathbb C: \text{Im}(k) < 0, \ 3\text{Re}(k)^2 - \text{Im}(k)^2 > 1, \ \text{Re}(k) < 0 \right\},
\\
&D_2^- = \left\{k \in \mathbb C: \text{Im}(k) < 0, \ 3\text{Re}(k)^2 - \text{Im}(k)^2 > 1, \ \text{Re}(k) > 0 \right\},
\end{aligned}
\end{equation}
and the discriminant of the cubic $\Omega + n\omega_0$ is
$\Delta = 4-27(n\omega_0)^2$, which defines symmetric critical
interval endpoints
$\omega_{cr}^\pm = \pm \omega_{cr} = \pm 2/(3 \sqrt{3})$.
Moreover, in all three scenarios $\Delta<0$, $\Delta>0$ and $\Delta=0$,  each of the contours $\p D^+$, $\p D_1^-$ and $\p D_2^-$ contains precisely one of the three roots  $k_0$, $k_1^-$ and $k_2^-$ of the cubic $\Omega + n\omega_0$ respectively. 

As explained above, the challenge is how to justify that, among the three distinct real roots that arise in the case $\Delta>0$, i.e. $|n\omega_0|<2/(3\sqrt 3)$, the two roots $k_1^- \in \p D_1^- $ and $k_2^- \in \p D_2^-$ that lie on $\p D^- = \p D_1^- \cup \p D_2^-$  should be removable singularities of \eqref{Unk-gen} while the remaining root $k_0 \in \p D^+$ should not, despite the fact that it belongs to the domain of definition $\left\{\text{Im}(k) \leq 0 \right\}$ of $\widehat u(k, t)$ (see remark below \eqref{dpm-gen-def} for a possible explanation of this fact). 
Therefore, we consider equation \eqref{lkdv} in the context of  the wavemaker problem associated with the data
\begin{equation}\label{lkdv-sine-bc}
u(x, 0) = 0, 
\quad 
u(0, t) = \sin(-\omega_0 t)
\end{equation}
and restrict our attention to the scenario $|\omega_0|<2/(3\sqrt 3)$ which, in the case of the data \eqref{lkdv-sine-bc}, corresponds to $\Omega + n\omega_0$ having three distinct real roots since $u(0, t) = - \sum_{n=\pm 1} \frac{\text{sgn}(n)}{2i} e^{in\omega_0 t}$ and, regardless of which two roots end up satisfying \eqref{dn-gen-0}, the resulting D-N relations will imply that the Fourier series coefficients $b_n, c_n$ of the boundary values $u_x(0, t), u_{xx}(0, t)$ vanish whenever the coefficient $a_n$ of $u(0, t)$ does.

In what follows, we will show that choosing $k_1^-$ and $k_2^-$ as the two removable singularities of \eqref{Unk-gen} does indeed produce the correct D-N map, namely the first relation given in \eqref{dn-gen}. We will do so by independently deriving that D-N relation directly from the solution formula obtained for the linear KdV wavemaker problem \eqref{lkdv}-\eqref{lkdv-sine-bc} via the Fokas method, which is given by (see~\cite{fokas2008unified}, Example 1.12)
\begin{equation}\label{kdv-sol-wm-x}
u(x,t)
=
\frac{1}{2i\pi} \sum_{n=\pm 1} a_n \int_{\p D^+} e^{ikx} \Omega' \frac{e^{i n\omega_0 t} - e^{-i\Omega t}}{\Omega + n\omega_0} \, dk, \quad a_{\pm 1} = \mp \frac{1}{2i}.
\end{equation}

Note that the three roots $k_0$, $k_1^-$, $k_2^-$ of $\Omega + n\omega_0$  are removable singularities of the above integrand, which is therefore analytic in $k$. 
Hence, in the particular case that $k_0$, $k_1^-$, $k_2^-$ are real and distinct, we invoke Cauchy's theorem to deform the contour of integration $\p D^+$ in \eqref{kdv-sol-wm-x} to the modified contour $\p \widetilde D^+$ that bypasses the removable singularity $k_0 \in \p D^+$ from below as shown in Figure \ref{lkdv-f} and thus obtain
\begin{equation}\label{lkdv-sol-def} 
u(x,t)
=
\frac{1}{2i\pi} \sum_{n=\pm1} a_n \int_{\p \widetilde D^+} e^{ikx} \Omega' \frac{e^{i n\omega_0 t} - e^{-i\Omega t}}{\Omega + n\omega_0} \, dk.
\end{equation}
Since $k_0$ no longer lies on the contour of integration, we can split the fraction in the integrand to write
\begin{equation}\label{u-break-0}
u(x,t)
=
\frac{1}{2i\pi} \sum_{n=\pm 1} a_n 
\left(
e^{in\omega_0 t} \int_{\p \widetilde D^+} e^{ikx} \frac{\Omega'}{\Omega + n\omega_0} \, dk
-
 \int_{\p \widetilde D^+} e^{ikx} \frac{\Omega' e^{-i\Omega t}}{\Omega + n\omega_0} \, dk
\right).
\end{equation}

Furthermore, for $k\in \left\{\text{Im}(k) \geq 0: k\notin D^+\right\}$ we have $|e^{-i\Omega t}| \leq 1$ so if, in addition, $k$ is large in modulus and stays away from $k_1^-$ and $k_2^-$, then 
$\left|
 \Omega' \frac{e^{-i\Omega t}}{\Omega + n\omega_0}
\right|
\lesssim
\frac{1}{|k|}$.
This decay bound combined with analyticity allows us to employ Cauchy's theorem and Jordan's lemma in the second integral of \eqref{u-break-0} in order to deform $\p \widetilde D^+$ to the union $\mathcal L_1 \cup \mathcal L_2$ of the half-lines (with orientation as shown in Figure \ref{lkdv-f})
\begin{equation}\label{L1L2-def}
\mathcal L_1 = \big\{k \in \mathbb C: k = -i + r e^{i\frac{2\pi}{3}}, \ r \geq 0 \big\},
\quad
\mathcal L_2 = \big\{k \in \mathbb C: k = -i + r e^{i\frac{\pi}{3}}, \ r \geq 0 \big\},
\end{equation}
 and hence express \eqref{u-break-0} in the form
\begin{equation}\label{u-break-fin}
u(x,t)
=
\frac{1}{2i\pi} \sum_{n=\pm 1} a_n 
\left(
e^{in\omega_0 t} \int_{\p \widetilde D^+} e^{ikx} \frac{\Omega'}{\Omega + n\omega_0} \, dk
-
 \int_{\mathcal L_1 \cup \mathcal L_2} e^{ikx}  \frac{\Omega' e^{-i\Omega t}}{\Omega + n\omega_0} \, dk
\right).
\end{equation}
We emphasize the significance of the bound $\left|\frac{\Omega'
    e^{-i\Omega t}}{\Omega + n\omega_0}\right| \lesssim
\frac{1}{|k|}$, which provides uniform decay in $k$ along the circular
arcs connecting $\p \widetilde D^+$ with $\mathcal L_1$ and $\mathcal
L_2$ and thus allows us to argue that the relevant integrals vanish by
invoking Jordan's lemma (or, alternatively, by a direct estimation of
those integrals along the lines of the integral $I_R$ below).

We will now use the modified solution formula \eqref{u-break-fin} in order to compute the leading-order long-time asymptotics of $u_x(x, t)$ for any finite $x > 0$. 
Differentiating with respect to $x$ yields
\begin{equation}\label{ux-break-0}
u_x(x,t)
=
\frac{1}{2\pi} \sum_{n=\pm 1} a_n 
\left(
e^{i n\omega_0 t} \int_{\p \widetilde D^+} e^{ikx} \frac{k \Omega'}{\Omega + n\omega_0} \, dk
-
 \int_{\mathcal L_1 \cup \mathcal L_2} e^{ikx} \frac{k \Omega' e^{-i\Omega t}}{\Omega + n\omega_0} \, dk
\right).
\end{equation}
Let $\p \widetilde D_R^+$ denote the finite portion of the contour $\p \widetilde D^+$ that lies below the circular arc $C_R = \left\{Re^{i\theta}: \frac{\pi}{3} \leq \theta \leq \frac{2\pi}{3}, R>1 \right\}$. 
With this notation, the first integral on the right-hand side of \eqref{ux-break-0} can be regarded as the limit
$$
\int_{\p \widetilde D^+} e^{ikx} \frac{k \Omega'}{\Omega + n\omega_0} \, dk
=
\lim_{R\to\infty} \int_{\p \widetilde D_R^+} e^{ikx} \frac{k \Omega'}{\Omega + n\omega_0} \, dk.
$$
Therefore, by Cauchy's residue theorem, 
\begin{equation}\label{crt-lkdv}
\int_{\p \widetilde D^+} e^{ikx} \frac{k \Omega'}{\Omega + n\omega_0} \, dk
=
2i\pi \, \text{Res}\left[e^{ikx} \frac{k \Omega'}{\Omega + n\omega_0}\right]_{k=k_0}
-
\lim_{R\to\infty} I_R(x, t)
\end{equation}
where  
$I_R(x, t) := \int_{C_R} e^{ikx} \frac{k \Omega'}{\Omega + n\omega_0} \, dk$.
We claim that $\lim_{R\to\infty} I_R(x, t) = 0$ for any $x>0$, $t>0$ and hence, calculating the residue in \eqref{crt-lkdv}, 
\begin{equation}\label{crt-lkdv-2}
\int_{\p \widetilde D^+} e^{ikx} \frac{k \Omega'}{\Omega + n\omega_0} \, dk
=
2i\pi k_0 e^{ik_0 x}.
\end{equation}
Indeed, for $k \in C_R$ we have $\left|e^{ikx}\right| = e^{-xR\sin\theta}$, $|\Omega'| \leq 1+3R^2$ and, for $R^2\geq 2 [1+2/(3\sqrt 3)]$ (recall that $n=\pm 1$ and $|n\omega_0| = |\omega_0|<2/(3\sqrt 3)$), 
$|\Omega + n\omega_0| 
=
\left|Re^{i\theta}-R^3e^{3i\theta} + n\omega_0\right|
\geq
\big|Re^{i\theta}-R^3e^{3i\theta}\big| -  |n\omega_0| 
\geq
R^3 - R -  |n\omega_0|
\geq
\frac 12 R^3$.
Thus,  
\begin{equation*}
\left|I_R(x, t)\right|
\leq 
\frac{2R^2 (1+3R^2)}{R^3}
\int_{\frac{\pi}{3}}^{\frac{2\pi}{3}} e^{-xR\sin\theta}  d\theta
\end{equation*}
and, using the symmetry of $\sin\theta$ about $\frac \pi 2$ as well as the well-known inequality $\sin\theta \geq \frac{2\theta}{\pi}$, $0 \leq \theta \leq \frac \pi 2$, we find
\begin{equation}
\left|I_R(x, t)\right|
\leq
\frac{2R^2 (1+3R^2)}{R^3} \cdot 
\frac{\pi}{2xR} \left(e^{-\frac{2xR}{3}} - e^{-xR}\right),
\nn
\end{equation}
which implies the claim for any $x>0$, $t > 0$.
In turn, combining \eqref{ux-break-0} with \eqref{crt-lkdv-2} we obtain
\begin{equation}\label{ux-break}
u_x(x,t) - \sum_{n=\pm 1} i k_0 a_n e^{ik_0x+in\omega_0 t}
=
-\frac{1}{2\pi} \sum_{n=\pm 1} a_n 
 \int_{\mathcal L_1 \cup \mathcal L_2} e^{ikx} \frac{k \Omega' e^{-i\Omega t}}{\Omega + n\omega_0} \, dk.
\end{equation}
\begin{remark}
We could have chosen the contour $\p \widetilde D^+$ to bypass $k_0$ from above. This would have made the first integral of \eqref{ux-break-0} equal to zero and would have moved the residue contribution over to the second integral.
\end{remark}

Next, we will prove that, for any $x>0$, the integral on the right-hand side of \eqref{ux-break} vanishes as $t\to\infty$. 
We emphasize that the deformation from $\p \widetilde D^+$ to $\mathcal L_1 \cup \mathcal L_2$ in \eqref{u-break-fin} was possible because $\text{Re}(i\Omega) \geq 0$ in the upper half-plane and outside $D^+$. In fact, $\text{Re}(i\Omega) \geq 0$ also in the lower half-plane and outside $D_1^-$ and $D_2^-$ and, more precisely, $\text{Re}(i\Omega) > 0$  for $k\in\mathcal L_1 \cup \mathcal L_2 \setminus \big\{\pm 1/\sqrt 3\big\}$ so that the exponential $e^{-i\Omega t}$ is strictly decaying there. This fact will allow us to  show that the integral on the right-hand side of \eqref{ux-break} decays like $O(1/\sqrt t)$ as $t\to\infty$.
We remark that our argument will not work without the presence of the exponential $e^{-i\Omega t}$, which decays as $|k| \to \infty$ along $\mathcal L_1 \cup \mathcal L_2$. Thus, the deformation from $\p \widetilde D^+$ to $\mathcal L_1 \cup \mathcal L_2$ is only possible and helpful for the second integral in \eqref{u-break-0} and not for the first one, which does not involve $e^{-i\Omega t}$. 

Since the estimation along the half-lines $\mathcal L_1$ and $\mathcal L_2$ is entirely analogous, we provide the details only for   $\mathcal L_2$. Then, $k=-i + r e^{i\frac{\pi}{3}} = \frac r2 + i\big(\frac{\sqrt 3 r}{2} - 1\big)$ and so
$\text{Re}(i\Omega) 
=
\frac 32\big(r-\frac{2}{\sqrt 3}\big)^2$.
Importantly, observe that due to the definition of $\mathcal L_2$ the expression for  $\text{Re}(i\Omega)$ is quadratic (as opposed to cubic) in $r$.
Also, similar to the proof of \eqref{ux-break}, for $r$ sufficiently large, say $r\geq 10$, we have $\left|\Omega+n\omega_0\right| \gtrsim r^3$. 
Furthermore, $\left|\Omega+n\omega_0\right|$ is strictly positive by the definition of  $\mathcal L_2$ and the fact that $|n\omega_0| = |\omega_0| < 2/(3 \sqrt 3)$, i.e. $\left|\Omega+n\omega_0\right| \geq C(n\omega_0)>0$ for some positive constant depending on $n\omega_0$. Hence,
\begin{equation}
\begin{aligned}
&\quad
\left|
\int_{\mathcal L_2} e^{ikx} \frac{k \Omega' e^{-i\Omega t}}{\Omega + n\omega_0} \, dk
\right|
\leq
\int_{0}^{10} e^{-\big(\frac{\sqrt 3}{2} r - 1\big)x}  \frac{\left|k \Omega'\right| e^{-\text{Re}(i\Omega) t}}{\left|\Omega + n\omega_0\right|} \, dr
+
\int_{10}^\infty e^{-\big(\frac{\sqrt 3}{2} r - 1\big)x} \frac{\left|k \Omega'\right| e^{-\text{Re}(i\Omega) t}}{\left|\Omega + n\omega_0\right|} \, dr
\\
&\lesssim
\int_{0}^{10} e^x  \frac{10 \cdot \big(1+3\cdot 10^2\big) \cdot e^{-\frac 32\big(r-\frac{2}{\sqrt 3}\big)^2 t}}{C(n\omega_0)} \, dr
+
\int_{10}^\infty e^{-\big(\frac{\sqrt 3}{2} r - 1\big)x} \frac{r \cdot r^2 \cdot e^{-\frac 32\big(r-\frac{2}{\sqrt 3}\big)^2 t}}{r^3} \, dr
\\
&\lesssim
e^x \int_{0}^\infty \left(\tfrac{1}{C(n\omega_0)} + 1\right) e^{-\frac 32\big(r-\frac{2}{\sqrt 3}\big)^2 t} dr
\leq
e^x \int_{-\infty}^\infty \left(\tfrac{1}{C(n\omega_0)} + 1\right) e^{-\frac 32 t r^2} dr
\simeq
e^x \, \frac{\frac{1}{C(n\omega_0)} + 1}{\sqrt t}.
\label{I2-est}
\end{aligned}
\end{equation}
A similar result holds for the integral along $\mathcal L_1$ and hence for \eqref{ux-break}. 

Overall, we have shown that, in the scenario $|\omega_0| < 2/(3\sqrt 3)$ of three distinct real roots of $\Omega+n\omega_0$, for any $x, t > 0$ the spatial derivative of the solution  to the linear KdV wavemaker problem \eqref{lkdv}-\eqref{lkdv-sine-bc} satisfies
\begin{equation}\label{lkdv-wm-ux-est}
\Big|
u_x(x,t) - \sum_{n=\pm 1} i k_0 a_n e^{ik_0x+in\omega_0 t}
\Big|
\lesssim
\frac{1}{2\pi} \sum_{n=\pm 1} |a_n| \left(\tfrac{1}{C(n\omega_0)} + 1\right) \frac{e^x}{\sqrt t}.
\end{equation}
In other words, for each $x \in [0,L]$, some $L > 0$,
\begin{equation}
  \label{eq:kdv_subcritical_Neumann}
u_x(x,t) \sim \sum_{n=\pm 1} i k_0 a_n e^{ik_0x+in\omega_0 t},
\quad t \to \infty.
\end{equation}
As we will see in Sec.~\ref{kdv-asymp-s}, the bound
\eqref{lkdv-wm-ux-est} is not tight and can be refined to
$O(t^{-3/2})$. Taking the limit, we find
$\lim_{x\to 0^+} u_x(x, t) \sim \sum_{n=\pm 1} i k_0 a_n e^{in\omega_0
  t}$ as $t \to \infty$,
in agreement with the first of the D-N relations \eqref{dn-gen}. A
similar analysis of the solution representation \eqref{u-break-fin}
can be used to prove that the second spatial derivative $u_{xx}$
satisfies the second of the D-N relations \eqref{dn-gen}.

We have shown that the D-N relations are obtained from the solution
formula \eqref{u-break-fin}.  In fact, essentially the same analysis
implies the asymptotic form of the solution
\begin{equation}
  \label{eq:6}
  u(x,t) \sim \sum_{n=\pm 1} a_n e^{ik_0x + in\omega_0 t}, \quad t \to
  \infty, \quad x = \mathcal{O}(1) .
\end{equation}
By linearity, this result extends to the more general boundary data
$g_0(t) \sim \sum_{n \in \mathbb{Z}} a_n e^{in\omega_0 t}$.  In
particular, we have also proven that the solution is asymptotically
time periodic for each fixed $x \ge 0$.  Thus, for the linear KdV
wavemaker problem, we have proven Theorem \ref{dn-t} and the radiation
condition \eqref{eq:general_radiation_condition} (Corollary
\ref{rad-cond}) for any frequency $\omega_0 \in \mathbb{R}$ without
the assumption of asymptotic time periodicity of the solution.

The asymptotic solution \eqref{eq:6} is only valid for fixed $x$ as
$t \to \infty$.  In Sec.~\ref{kdv-asymp-s}, we extend the analysis of
the linear KdV wavemaker problem to obtain an asymptotic expansion of the
solution for large $t$ that is uniform in $x \ge 0$.

\begin{remark}
  The above arguments and, in particular, estimate \eqref{I2-est},
  crucially rely on the fact that $\text{Re}(i\Omega) \geq 0$ along
  $\mathcal L_1$ and $\mathcal L_2$. Attempting to deform the contour
  $\p D^+$ differently in order to also enclose one or both of the
  roots $k_1, k_2$ would force us to include contours in regions where
  $\text{Re}(i\Omega) < 0$ (see Figure \ref{lkdv-f}) and hence would
  result in exponentially growing terms, thus an estimate like
  \eqref{I2-est} would no longer hold. This is why the only root
  contributing to the D-N map (through the residue calculation
  \eqref{crt-lkdv}) is the pole along $\p D^+$, namely $k_0$.
\end{remark}

\subsection{A general nonlocal second-order model: $A_{-2} > 0$ and $A_3 = 0$}
\label{2-nonlocal}

For $A_{-2} > 0$ and $A_3 = 0$, the third-order equation \eqref{gen-mod} becomes of second order but still involves the mixed spatiotemporal derivative term, namely
\begin{equation}\label{gen-mod-a3=0}
\big(1-A_{-2} \p_x^2\big) u_t = i A_0 u + A_1 u_x + i A_2 u_{xx}.\end{equation}
Moreover, the equation $\Omega + n \omega_0 \big(1+A_{-2}k^2\big) = 0$, which in general corresponds to the cubic \eqref{poles-eq}, now reduces to the quadratic
\begin{equation}\label{cubic-a3=0}
\left(A_2 + n\omega_0 A_{-2}\right) k^2 - A_1 k - A_0 + n\omega_0 = 0.
\end{equation}
In addition, according to definition \eqref{dpm-gen-def}, 
\begin{equation}\label{dm-bbm}
D^- =
\left\{ k\in\mathbb C:
\text{Im}(k) < 0
\text{ and } 
A_1 - 2 \left(A_2 + A_0 A_{-2}\right) \text{Re}(k) - A_1 A_{-2}  \text{Re}(k)^2 - A_1 A_{-2}  \text{Im}(k)^2 > 0
\right\}.
\end{equation}
Thus, for real solutions of \eqref{cubic-a3=0} that lie along the boundary $\p D^-$, we must have
\begin{equation}\label{real-cond-a3=0}
A_1 A_{-2}  k^2 + 2 \left(A_2 + A_0 A_{-2}\right) k -A_1 \leq  0.
\end{equation}
Likewise, the region $D^+$ is given by
\begin{equation}\label{dp-bbm}
D^+ =
\left\{ k\in\mathbb C:
\text{Im}(k) > 0
\text{ and } 
A_1 - 2 \left(A_2 + A_0 A_{-2}\right) \text{Re}(k) - A_1 A_{-2}  \text{Re}(k)^2 - A_1 A_{-2}  \text{Im}(k)^2 < 0
\right\}
\end{equation}
and, therefore, any real solution of \eqref{cubic-a3=0} that lies along the boundary $\p D^+$ must satisfy
\begin{equation}\label{real-cond-a3=0+}
A_1 A_{-2}  k^2 + 2 \left(A_2 + A_0 A_{-2}\right) k -A_1 \geq  0.
\end{equation}

Next, we focus on the solution of equation \eqref{cubic-a3=0}. 
We begin with the scenario in which the parameters $A_{-2}$, $A_2$,  $\omega_0$ are such that $A_2 + n\omega_0 A_{-2} \neq 0$ for all $n\in\mathbb Z$.
Then, \eqref{cubic-a3=0} is a quadratic equation with solutions  
\begin{equation}\label{cubic-a3=0-roots}
k_\pm = k_\pm(n) := \frac{A_1 \pm \sqrt{\Delta}}{2\left(A_2 + n\omega_0 A_{-2}\right)},
\end{equation}
where the discriminant $\Delta$ is given by
$$
\Delta = A_1^2+4\left(A_0-n\omega_0\right)\left(A_2 + n\omega_0 A_{-2}\right).
$$
Observe that $\Delta$ is itself a quadratic in the quantity $n\omega_0$ with corresponding discriminant
$$
\widetilde \Delta 
=
16\left[ \left(A_0A_{-2} + A_2\right)^2 + A_{-2} A_1^2 \right].
$$
Note that $\widetilde \Delta\geq  0$ since $A_{-2}>0$. 
\\[2mm]
\textbf{Case 1: $\widetilde \Delta > 0$.} This case  amounts to either $A_1 \neq 0$ or $A_2 + A_0A_{-2} \neq 0$. Then,  $\Delta$ has two distinct real roots
\begin{equation}
\omega_{cr}^\pm = \frac{A_0A_{-2} - A_2 \pm \sqrt{\left(A_2 + A_0A_{-2}\right)^2 + A_1^2 A_{-2}}}{2A_{-2}},
\end{equation}
that determine a bifurcation in solution behavior.
If $n\omega_0>\omega_{cr}^+$ or $n\omega_0<\omega_{cr}^-$ then $\Delta < 0$ and the solutions $k_+$,  $k_-$ of equation \eqref{cubic-a3=0} given by  \eqref{cubic-a3=0-roots} form a complex conjugate pair; 
if $\omega_{cr}^- < n\omega_0 < \omega_{cr}^+$ then $\Delta > 0$ and $k_+$,  $k_-$ are real and distinct;
and if $n\omega_0 = \omega_{cr}^+$ or $n\omega_0 = \omega_{cr}^-$ then
$k_+$, $k_-$ reduce to a double real solution.

Importantly, for each $n\in\mathbb Z$ there is precisely one solution of \eqref{cubic-a3=0} that lies on $\p D^-$. 
This is straightforward to see if $n\omega_0>\omega_{cr}^+$ or $n\omega_0<\omega_{cr}^-$ since, by the definition of the regions $D^+$ and $D^-$, any non-real solution of \eqref{cubic-a3=0} must lie either on $\p D^+$ or on $\p D^-$. 
On the other hand, if $k_\pm \in \mathbb R$  then more details are needed: 
\\[2mm]
(i) If $A_1=0$ and $A_2 + A_0 A_{-2} > 0$, then \eqref{real-cond-a3=0} requires $k \leq  0$. However, by Vieta's formulae, $k_+ + k_- = 0$. Thus precisely one (possibly double, if equal to zero) solution of \eqref{cubic-a3=0} lies on $\p D^-$ and, since $\mathbb R \subset \left(\p D^+ \cap \p D^-\right)$, the remaining solution lies on $\p D^+$.
\\[2mm]
(ii) If $A_1=0$ and $A_2 + A_0 A_{-2} < 0$, then \eqref{real-cond-a3=0} requires $k \geq  0$ and, as in (i), precisely one (possibly double, if equal to zero) solution of \eqref{cubic-a3=0} lies on $\p D^-$ and the remaining solution lies on $\p D^+$. 
\\[2mm]
(iii) If $A_1 > 0$ then, since $A_{-2}>0$, \eqref{real-cond-a3=0} requires that $\beta_- \leq  k \leq  \beta_+$, where
\begin{equation}\label{bpm-def}
\beta_\pm = \frac{-\left(A_2 + A_0 A_{-2}\right) \pm \sqrt{\left(A_2 + A_0 A_{-2}\right)^2+A_1^2A_{-2}}}{A_1A_{-2}}.
\end{equation}

First, we claim that $\beta_- \leq  k_- \leq  \beta_+$ (and hence, via \eqref{real-cond-a3=0}, that $k_- \in \p D^-$).
Indeed, note that
$\omega_{cr}^\pm = \frac{A_1 \beta_\pm}{2} + A_0$
and, since $k_\pm \in \mathbb R$ amounts to having $\omega_{cr}^- \leq  n\omega_0 \leq  \omega_{cr}^+$, it follows that
$
A_1 \beta_- \leq  2 \left(n\omega_0 - A_0\right) \leq  A_1 \beta_+$.
Hence, observing that
$k_- = \frac{2\left(n\omega_0-A_0\right)}{A_1+\sqrt \Delta}$
and $A_1+\sqrt \Delta > 0$, we infer that
$
\frac{A_1 \beta_-}{A_1+\sqrt \Delta} \leq  k_- \leq  \frac{A_1 \beta_+}{A_1+\sqrt \Delta}
$
which, in view of the inequalities $0<\frac{A_1}{A_1+\sqrt \Delta}\leq  1$, $\beta_+ > 0$, $\beta_-  < 0$, allows us to conclude that $\beta_- \leq  k_- \leq  \beta_+$ as claimed. Note that the equalities are attained only when $\Delta=0$, i.e. only when $n\omega_0 = \omega_{cr}^+$ or $n\omega_0 = \omega_{cr}^-$, in which case $k_-(\omega_{cr}^+) = k_+(\omega_{cr}^+) = \beta_+$ or $k_-(\omega_{cr}^-) = k_+(\omega_{cr}^-) = \beta_-$, respectively.

In addition, we claim that $k_+ \geq  \beta_+$ or $k_+ \leq  \beta_-$, with equalities attained only when $\Delta =0$  (and hence that $k_+ \notin \p D^-$ except when $k_+ = k_-$). 
Indeed, in view of \eqref{cubic-a3=0-roots},  \eqref{bpm-def} and the fact that $A_1 A_{-2}>0$, if $A_2 + n\omega_0 A_{-2} >0$ then the inequality 
$A_1^2 A_{-2}  
\geq 
2\left(A_2 + n\omega_0 A_{-2}\right) \Big[ \sqrt{\left(A_2 + A_0 A_{-2}\right)^2+A_1^2A_{-2}} -\left(A_2 + A_0 A_{-2}\right) \Big]$,
which is valid since $n\omega_0 \leq  \omega_{cr}^+$ implies
$2\left(A_2 + n\omega_0 A_{-2}\right)
\leq 
A_0A_{-2} + A_2 + \sqrt{\left(A_2 + A_0A_{-2}\right)^2 + A_1^2 A_{-2}}$,
yields $k_+ \geq  \beta_+$. 
Similarly, if $A_2 + n\omega_0 A_{-2} < 0$ then, noting that $n\omega_0 \geq \omega_{cr}^-$ implies 
$-2\left(A_2 + n\omega_0 A_{-2}\right)
\leq 
\sqrt{\left(A_2 + A_0A_{-2}\right)^2 + A_1^2 A_{-2}} -\left(A_0A_{-2} + A_2\right)$,
it follows that
$A_1^2 A_{-2}  
\geq 
-2\left(A_2 + n\omega_0 A_{-2}\right) \Big[ \sqrt{\left(A_2 + A_0 A_{-2}\right)^2+A_1^2A_{-2}} + \left(A_2 + A_0 A_{-2}\right) \Big]$,
and, in turn, $k_+ \leq \beta_-$, establishing the claim.
\\[3mm]
(iv) Finally, if $A_1 < 0$ then, since $A_{-2}>0$, we have $\beta_+<0$ and $\beta_->0$, and \eqref{real-cond-a3=0} requires either $k\geq  \beta_-$ or $k \leq  \beta_+$.
Via entirely analogous steps as in (iii), using the fact that $\omega_{cr}^-\leq  n\omega_0 \leq  \omega_{cr}^+$ we obtain $\beta_+ \leq  k_+ \leq  \beta_-$ and $k_-\geq  \beta_-$  (for $A_2 + n\omega_0 A_{-2}<0$) or $k_-\leq  \beta_+$ (for $A_2 + n\omega_0 A_{-2}>0$), with equalities possible only when $\Delta=0$, i.e. only when $n\omega_0 = \omega_{cr}^+$ or $n\omega_0 = \omega_{cr}^-$, in which case $k_-(\omega_{cr}^+) = k_+(\omega_{cr}^+) = \beta_+$ or $k_-(\omega_{cr}^-) = k_+(\omega_{cr}^-) = \beta_-$, respectively.
Thus, once again we see that there exists precisely one (possibly double) root along $\p D^-$, namely $k_-$, while the other root, i.e. $k_+$, lies on $\p D^+$.

Therefore, in the case $\widetilde \Delta > 0$ we have shown that for each $n\in\mathbb Z$ there is precisely one (possibly double) solution of \eqref{cubic-a3=0} that lies on $\p D^-$, with the remaining solution lying on $\p D^+$.
\\[2mm]
\textbf{Case 2: $\widetilde \Delta = 0$.} This case is much simpler because it requires $A_1= A_2 + A_0A_{-2} = 0$ and so it reduces our second-order model \eqref{gen-mod-a3=0} to the equation
\begin{equation}\label{red-mod}
\big(1-A_{-2} \p_x^2\big) \left(u_t- i A_0 u\right) = 0,
\end{equation}
which can be solved by elementary techniques. Indeed, the  ODE in $x$ along with the boundary conditions $u(0, t) = g_0(t)$ and $\lim_{x\to\infty} u(x, t) = 0$ can be integrated to yield an expression for the quantity $u_t- i A_0 u$, which can be integrated in $t$ to imply, upon imposing the initial condition $u(x, 0) = u_0(x)$,
\begin{equation}\label{exp-sol}
u(x, t) = e^{iA_0t}  u_0(x) +  \big[g_0(t) - e^{iA_0t}  g_0(0)\big] e^{-\frac{x}{\sqrt{A_{-2}}}}.
\end{equation}

The analysis of the scenario $A_2 + n\omega_0 A_{-2} \neq 0$ for all $n\in\mathbb Z$ is complete and leads to the following conclusions:
\\[2mm]
$\bullet$ If either $A_1 \neq 0$ or $A_2 + A_0A_{-2} \neq 0$, then one of the solutions $k_+$, $k_-$ of the quadratic \eqref{cubic-a3=0} given by \eqref{cubic-a3=0-roots} lies on the contour $\p D^-$ while the other solution does not lie on $\p D^-$ (with the exception of the special case $k_+ = k_-$).  
Denote the solution along $\p D^-$ by $k_1^-$ and the solution off $\p D^-$ (which, as explained above, lies on $\p D^+$) by $k_0$. 
Then, setting $k=k_1^-$ into the D-N map condition \eqref{dn-gen-0}, which for the second-order model \eqref{gen-mod-a3=0} reads
\begin{equation}\label{dn-a3=0}
i \left(A_2 + n \omega_0 A_{-2}\right) b_n 
= \left[ \left(A_2 + n \omega_0 A_{-2}\right) k - A_1 \right]  a_n,
\end{equation}
and using the Vieta formula $k_0 + k_1^- = \frac{A_1}{A_2 + n \omega_0 A_{-2}}$ for equation \eqref{cubic-a3=0}, we obtain the D-N relation $b_n = i k_0  a_n$ which was conjectured in \eqref{dn-gen}.   
\\[2mm]
$\bullet$ If $A_1=A_2 + A_0A_{-2} = 0$, then the solution to the Dirichlet problem of equation \eqref{gen-mod-a3=0} on the half-line  is given by the simple formula \eqref{exp-sol} and the D-N map can be computed by directly differentiating this formula.

\vskip 3mm

Finally, we consider the scenario of $A_2 + n_*\omega_0 A_{-2} = 0$ for some $n_*\in\mathbb Z$, i.e. the possibility that the parameters $A_{-2}$, $A_2$, $\omega_0$ are such that $n_* := -A_2/(\omega_0A_{-2}) \in \mathbb Z$. Then, equation \eqref{dn-gen-0} does not yield information about $b_{n_*}$, as either it becomes an identity (if $A_1=0$) or it forces $a_{n_*} =0$ (if $A_1 \neq 0$). 

If $A_1=0$, then equation \eqref{cubic-a3=0} evaluated at $n_*$ implies $A_0A_{-2} + A_2=0$ so we recover the simple model \eqref{red-mod} with the explicit solution formula \eqref{exp-sol}.

If $A_1 \neq 0$, then we note that under the transformation $u \mapsto u + a_{n_*}  e^{in_* \omega_0 t}$, which in view of \eqref{g0-per} eliminates the term corresponding to $n_*$ from the series expansion of $u(0, t)$, equation \eqref{gen-mod-a3=0} becomes
\begin{equation}\label{gen-mod-a3=0-v}
\big(1-A_{-2} \p_x^2\big) \big(u_t  + in_* \omega_0 a_{n_*}  e^{in_* \omega_0 t}\big)  = i A_0 \big(u + a_{n_*} e^{in_* \omega_0 t}\big) + A_1 u_x + i A_2  u_{xx}.
\end{equation}
If $A_0 A_{-2} + A_2 = 0$ then $n_* \omega_0 = A_0$ and equation \eqref{gen-mod-a3=0-v} is identical to equation \eqref{gen-mod-a3=0}, so without loss of generality we can take $a_{n_*} = 0$. In particular, this is the case for the linear BBM equation considered below.
On the other hand, if $A_0 A_{-2} + A_2 \neq 0$ then $n_* \omega_0
\neq A_0$ and \eqref{gen-mod-a3=0-v} is different from
\eqref{gen-mod-a3=0}, so the case $n_* = -A_2/(\omega_0A_{-2}) \in
\mathbb Z$ and $A_1 \left(A_0 A_{-2} + A_2\right) \neq 0$ is not
covered by our analysis unless the associated Dirichlet datum actually
satisfies $a_{n_*} = 0$.

Concerning the value of $b_{n_*}$, if the Dirichlet datum is actually periodic (and not just asymptotically periodic) then using the compatibility between the initial and boundary values at $x=t=0$ and assuming that $u(x, t)$ is continuously differentiable up to the boundary, we find 
$b_{n_*} = u_0'(0) - \sum_{n\in\mathbb Z, n\neq n_*} b_n$,
with the coefficients $b_n$ for $n\neq n_*$ given by the first of the D-N relations \eqref{dn-gen}.
Furthermore, if $A_2=0$ and hence $n_*=0$, then $b_{n_*} = b_0$ can be determined even if the solution is only asymptotically periodic. In particular, evaluating equation \eqref{gen-mod-a3=0} at $x=0$ and then integrating over $\big[t, t + \frac{2\pi}{\omega_0}\big]$ with $t$ arbitrary, we exploit the periodicity of the relevant exponentials to deduce, by taking the limit $t\to\infty$, 
\begin{equation}
b_0 = -i \frac{A_0}{A_1} \, a_0.
\end{equation}
Note that, in the particular case of the linear BBM equation, the above expression yields $b_0 = 0$. 

This proves Theorem \ref{dn-t} when $A_3 = 0$.

\begin{remark}
  \textbf{(connection to the radiation condition
\eqref{eq:general_radiation_condition})}  
The group velocity associated with the second-order equation
\eqref{gen-mod-a3=0} is  $c_g = \omega'(k) = \frac{A_1A_{-2} k^2 + 2(A_2+A_0A_{-2})k - A_1}{(1+ A_{-2} k^2)^2}$.
Thus, for a real solution of the quadratic \eqref{cubic-a3=0}, the
condition that this solution lies on $\p D^+$ (equation \eqref{real-cond-a3=0+}) amounts precisely to the group velocity being non-negative at that solution. 
In particular, in the non-trivial scenario of either $A_1\neq 0$ or $A_2 + A_0A_{-2}\neq 0$, recall that  $k_0 \in \p D^+ \cap \mathbb
R$ is possible if either $\omega_{cr}^- <  n\omega_0 < \omega_{cr}^+$  or $n\omega_0 = \omega_{cr}^\pm$. In the first case we have a strictly positive group velocity, while in the second case the group velocity is zero.
Furthermore, if $k_0 \in \p D^+ \setminus \mathbb R$ then 
$\textnormal{Im} \, k_0  >
0$.

  This proves Corollary \ref{rad-cond} when $A_{3} = 0$.
\end{remark}

\subsubsection*{The linear BBM on the half-line}
\label{lbbm-ss}
Recall that the main idea behind deriving the D-N map \eqref{dn-gen} is that each root of $\Omega + n \omega_0 \big(1+A_{-2}k^2\big)$ that occurs along $\p D^-$ should be a removable  singularity of the right-hand side of \eqref{Unk-gen} and hence a solution of \eqref{dn-a3=0}.
In this connection, the analysis provided above for the quadratic \eqref{cubic-a3=0} associated with the second-order model \eqref{gen-mod-a3=0} justifies the assumption made in the derivation of  \eqref{dn-gen} that precisely one of the two roots of $\Omega + n \omega_0 \big(1+A_{-2}k^2\big)$ lies on the contour $\p D^-$. 
Nevertheless, the additional assumption that \textit{only} this root
must satisfy \eqref{dn-a3=0} has not yet been justified and, although
this is an easy task for $n \omega_0$ outside the interval
$[\omega_{cr}^-,\omega_{cr}^+]$ --- since then the root off $\p D^-$
is in the upper half-plane where \eqref{ft-hl} is generally not
defined --- it becomes more subtle for values within the interval.
In what follows, we address this outstanding uniqueness issue in the case of the wavemaker problem \eqref{lkdv-sine-bc} for the linear BBM equation, which reads
\begin{equation}\label{lbbm}
u_t + u_x - u_{xxt} = 0
\end{equation}
and corresponds to the choice of parameters $A_0=A_2=0$, $A_{-2}=1$, $A_1=-1$ in the general model \eqref{gen-mod-a3=0}.
We remark that, despite the absence of the second spatial derivative term, the mixed derivative term --- an essential element of the model \eqref{gen-mod-a3=0} when compared to the model \eqref{gen-mod-a-2=0} --- is still present in equation~\eqref{lbbm} and makes the analysis quite different from that of the linear KdV equation.

\begin{figure}[ht!] 
\centering
\includegraphics[scale=0.45]{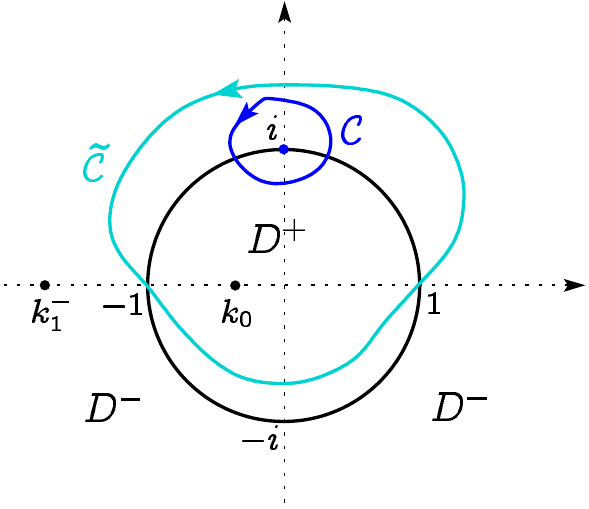}
\caption{The regions $D^+$ and $D^-$ for the linear BBM equation as
  given by \eqref{dpm-lbbm}.   In the case of the wavemaker data \eqref{lkdv-sine-bc}, for
  $0<|2\omega_0|<1$ the two roots of $\omega+n\omega_0$ are real and
  distinct with $k_1^- \in \p D^-$, $k_0 \in \p D^+$ as shown. Then, in
  order to compute the long-time asymptotics of the solution, we
  deform the contour of integration in the Fokas formula
  \eqref{lbbm-sol} from $\mathcal C$ to the contour
  $\widetilde{\mathcal C}$ along which $\text{Re}(i\omega) > 0$
  (except for $k=\pm 1$ where $\text{Re}(i\omega)=0$). This
  deformation yields a residue contribution from $k_0 \in \p D^+$ (as
  opposed to $k_1^- \in \p D^-$), justifying the assumption that only
  $k_1^-$ (and not $k_0$) should satisfy \eqref{dn-a3=0}.}
\label{lbbm-f}
\end{figure}

In the case of the linear BBM equation \eqref{lbbm}, the regions $D^-$ and $D^+$ of \eqref{dm-bbm} and \eqref{dp-bbm} are given by
\begin{equation}\label{dpm-lbbm}
D^-= \left\{ k\in\mathbb C: \text{Im}(k) < 0, \ \text{Re}(k)^2 + \text{Im}(k)^2 > 1 \right\},
\quad
D^+ = \left\{ k\in\mathbb C: \text{Im}(k) > 0, \ \text{Re}(k)^2 + \text{Im}(k)^2 < 1 \right\},
\end{equation}
i.e. $D^+$ is the interior of the upper half of the unit disk centered
at the origin and $D^-$ is the exterior of that disk in the lower half
of the complex plane, as shown in Figure \ref{lbbm-f}.  The critical
interval endpoints are symmetric
$\omega_{cr}^\pm = \pm \omega_{cr} = \pm 1/2$.
Furthermore,
\begin{equation}
  \label{eq:7}
  \Omega + n \omega_0 \big(1+A_{-2}k^2\big) = k + n\omega_0(1+k^2) = 0
\end{equation}
when (cf.~eqs.~\eqref{cubic-a3=0} and
\eqref{cubic-a3=0-roots})
\begin{equation}\label{lbbm-zeros}
k = k_\pm(n) =  
\left\{
\def\arraystretch{0.8}
\begin{array}{ll}
\dfrac{-1\pm \sqrt{1-(2n\omega_0)^2}}{2n\omega_0}, &n\neq 0, 
\\
0, &n=0.
\end{array}
\right.
\end{equation}

If $|2n\omega_0| > 1$ then \eqref{lbbm-zeros} form a pair of complex
conjugate roots that lie on the unit circle.
Thus, due to the domain of \eqref{ft-hl}, only the root lying in the lower half of the complex plane --- and hence on $\p D^-$ --- must satisfy \eqref{dn-a3=0}. The same is true if $|2n\omega_0| = 1$, since  then the roots \eqref{lbbm-zeros} reduce to a single double root equal to $\mp 1 \in \p D^-$.
On the other hand, if $0<|2n\omega_0| < 1$ then the roots \eqref{lbbm-zeros} are real and distinct, with one of them, denoted by $k_0$, lying inside $(-1, 1) \subset \p D^+$ and the other one, denoted by $k_1^-$, lying in $[-1, 1]^c \subset \p D^-$ (see Figure \ref{lbbm-f}). 
In this case, as noted above, since both roots are in the domain of \eqref{ft-hl}, it is not immediately clear why only $k_1^-$, and not $k_0$, should satisfy \eqref{dn-a3=0} as a removable singularity of \eqref{Unk-gen}. 
Next, we will show that this is indeed the case by taking advantage of the powerful solution representation provided by the Fokas method.

In particular, let us consider the linear BBM equation \eqref{lbbm} in
the context of the wavemaker problem associated with the data
\eqref{lkdv-sine-bc}. In this case, observe that the scenario of two
distinct real roots of \eqref{eq:7} corresponds to the range
$0<|2\omega_0|<1$, since
$u(0, t) = - \sum_{n=\pm 1} \frac{\text{sgn}(n)}{2i} e^{in\omega_0 t}$
and, regardless of which root ends up satisfying \eqref{dn-a3=0}, the
resulting D-N relation will imply that the Fourier series coefficient
$b_n$ of the Neumann value $u_x(0, t)$ vanishes whenever the
coefficient $a_n$ of $u(0, t)$ does. Hence, we work with
$0<|2\omega_0|<1$ and two distinct real roots as shown in Figure
\ref{lbbm-f}.
The solution formula is given by \cite{vd2013}
\begin{equation}\label{lbbm-sol}
u(x, t) 
= 
\frac{1}{2i\pi} \sum_{n=\pm 1} a_n \int_{\mathcal C}  e^{ikx}  \omega' \frac{e^{i n\omega_0 t}-e^{-i\omega t}}{\omega+n\omega_0} \, dk 
- e^{-x} \sin(\omega_0 t),
\quad
a_{\pm 1} = \mp \frac{1}{2i},
\end{equation}
where the contour $\mathcal{C}$ is a positively oriented, closed
contour in the upper half of the complex plane and enclosing $k=i$, as
shown in Figure \ref{lbbm-f} and
\begin{equation}
\omega = \omega(k) := \frac{k}{1+k^2}
\end{equation}
with $\omega'$ denoting the derivative of $\omega$. 
Note that the integrand in \eqref{lbbm-sol} has essential singularities at $k=\pm i$ due to the exponential $e^{-i\omega t}$. In addition, for $k\neq \pm i$ we have $\text{Re}\left(\frac{i\Omega}{ 1+ k^2}\right)
=
\text{Re}(i\omega)$ and so  inside the regions $D^-$ and $D^+$ of \eqref{dpm-lbbm} we have $\text{Re}(i\omega)<0$, while outside those regions we have $\text{Re}(i\omega)>0$. 

Thanks to analyticity and Cauchy's theorem, the contour of integration $\mathcal C$ in \eqref{lbbm-sol} can be deformed to the contour $\widetilde{\mathcal C}$ depicted in Figure \ref{lbbm-f}. 
Thus, differentiating in $x$, we have  
\begin{equation}\label{ux-lbbm}
u_x(x, t) 
= 
\frac{1}{2\pi} \sum_{n=\pm 1} a_n \int_{\widetilde{\mathcal C}}  e^{ikx}  k \omega' \frac{e^{i n\omega_0 t}-e^{-i\omega t}}{\omega+n\omega_0} \, dk + e^{-x} \sin(\omega_0 t).
\end{equation}
Since $\widetilde{\mathcal C}$ does not cross any singularities, we can split the integral with respect to $k$ in two parts:
\begin{equation}\label{ux-lbm-break}
u_x(x, t) 
= 
\frac{1}{2\pi} \sum_{n=\pm 1} a_n 
\left(
e^{i n\omega_0 t}
\int_{\widetilde{\mathcal C}}  e^{ikx}  \frac{k \omega'}{\omega+n\omega_0} \, dk 
-
 \int_{\widetilde{\mathcal C}}  e^{ikx} \frac{k \omega' e^{-i\omega t}}{\omega+n\omega_0}  \, dk 
\right)
+ e^{-x} \sin(\omega_0 t).
\end{equation}
The contour $\widetilde{\mathcal C}$ encloses the points $k=i$ and $k = k_0$, which are simple poles of the integrand of the first of the above integrals. Hence, by Cauchy's residue theorem and the fact that $k_0$ satisfies \eqref{lbbm-zeros}, 
\begin{equation}\label{ux-lbm-break-1}
\int_{\widetilde{\mathcal C}}  e^{ikx}  \frac{k\omega'}{\omega+n\omega_0} \, dk 
=
2\pi \left( e^{-x} + ik_0 e^{ik_0x} \right).
\end{equation}

Moreover, note that along $\widetilde{\mathcal C}$ we have $\text{Re}(i\omega)> 0$ except for  $k=\pm 1$ where $\text{Re}(i\omega)=0$. Thus, we expect that the second integral in \eqref{ux-lbm-break} admits a similar bound to \eqref{I2-est} for the integral on the right-hand side of \eqref{ux-break}. Below, we establish this rigorously.
Thanks to analyticity, without loss of generality we work with $\widetilde C = \left\{k \in \mathbb C: k = i + \sqrt 2 e^{i\theta}, -\frac \pi 2 \leq \theta \leq \frac{3\pi}{2} \right\}$. Then, $\text{Re}(k) = \sqrt 2 \cos \theta$, $\text{Im}(k) = 1 + \sqrt 2 \sin \theta$ and 
$
|1+k^2|^2 
=
4 \big(3 + 2\sqrt 2 \sin\theta\big)$,
implying 
$$
\phi(\theta) := \text{Re}(i\omega)
=
\frac{\big(1 + \sqrt 2 \sin \theta\big)^2}{2 \big(3 + 2\sqrt 2 \sin\theta \big)}
.
$$
Furthermore,
$
|k\omega'| \leq 1 + 1/\sqrt{2}$,
and $|\omega + n\omega_0| \geq |1-2n\omega_0|/2 > 0$ since
$|n\omega_0| < 1/2$.
Hence,
$$
\left|
\int_{\widetilde{\mathcal C}}  e^{ikx} \frac{k \omega' e^{-i\omega t}}{\omega+n\omega_0}  \, dk
\right|
\leqslant
\int_{-\frac \pi 2}^{\frac{3\pi}{2}} 
e^{-(1+\sqrt 2 \sin\theta) x} 
\frac{|k \omega'| e^{-\text{Re}(i\omega) t}}{|\omega+n\omega_0|} 
\sqrt 2 d\theta
\leqslant
\frac{4(\sqrt{2}+1)e^{(\sqrt{2}-1)x}}{|1-2n\omega_0|}
\int_{-\pi/2}^{\pi/2} e^{-t \phi(\theta)} d\theta.
$$
For $\theta \in [-\frac \pi 2, \frac \pi 2]$, $\phi(\theta) \geq 0$ with the minimum attained only at $-\frac{\pi}{4}$. This integral can be estimated using Laplace's method (e.g. see Lemma 6.2.3 in \cite{af2003}), we obtain
$
\int_{-\pi/2}^{\pi/2}  e^{-t\phi(\theta)} d\theta
\sim
\sqrt{\frac{2\pi}{t}} + O\left(\frac{1}{t\sqrt t}\right)
$
as $t\to\infty$.
Therefore, for any $x\geq 0$ and $|n\omega_0| < 1/2$, we conclude that
\begin{equation}\label{ux-lbm-break-2}
\left|
\int_{\widetilde{\mathcal C}}  e^{ikx} \frac{k \omega' e^{-i\omega t}}{\omega+n\omega_0}  \, dk
\right|
=
O\left(\frac{e^{(\sqrt 2-1)x}}{\sqrt t}\right), \quad t\to\infty.
\end{equation}

In view of \eqref{ux-lbm-break-1} and \eqref{ux-lbm-break-2}, the expression \eqref{ux-lbm-break} yields
\begin{equation}\label{ux-fin}
\left|
u_x(x, t) - \sum_{n=\pm 1} i k_0 a_n \, e^{ik_0x + i n\omega_0 t}
\right|
\lesssim
O\left(\frac{e^{(\sqrt 2-1)x}}{\sqrt t}\right), \quad t\to\infty.
\end{equation}
In particular, first taking the limit $x\to 0^+$ and then performing the above asymptotic analysis as $t\to\infty$, we obtain 
$\lim_{x\to 0^+} u_x(x, t) \sim \sum_{n=\pm 1} i k_0 a_n \, e^{ik_0x + i n\omega_0 t}$ as $t\to\infty$,
in agreement with the D-N map \eqref{dn-gen}.

A similar analysis of the solution formula \eqref{lbbm-sol} implies
the asymptotic solution for fixed $x \ge 0$
\begin{equation}
  \label{eq:8}
  u(x,t) \sim \sum_{n=\pm1} a_n e^{ik_0x + in\omega_0t}, \quad t \to
  \infty , \quad x = \mathcal{O}(1) .
\end{equation}
By linearity, this result extends to the more general boundary data
$g_0(t) \sim \sum_{n \in \mathbb{Z}} a_n e^{in\omega_0t}$.  We have
proven the asymptotic time periodicity of the solution to the linear
BBM wavemaker problem so that Theorem 1 and the radiation condition
\eqref{eq:general_radiation_condition} (Corollary \ref{rad-cond}) hold
for any frequency $\omega_0 \in \mathbb{R}$ without the need to assume
asymptotic time periodicity of the solution.

The asymptotic solution \eqref{eq:8} is only valid for fixed $x$ as
$t \to \infty$.  In Sec.~\ref{bbm-asymp-s}, we extend the analysis of
the linear BBM wavemaker problem to obtain an asymptotic expansion of
the solution for large $t$ that is uniform in $x \ge 0$.  We then
observe that this analysis quantitatively agrees with viscous
core-annular experiments in Sec.~\ref{exp-s}.

\begin{remark}
In order to establish estimate \eqref{ux-lbm-break-2}, it is crucial for $\widetilde{\mathcal C}$ to lie in regions where $\textnormal{Re}(i\omega)\geq 0$. This requirement implies that only $k_0$ (and not $k_1^-$) can be enclosed by $\widetilde{\mathcal C}$, since attempting to also enclose $k_1^-$ would push $\widetilde C$ into regions where $\textnormal{Re}(i\omega)<0$. Consequently, only $k_0$ (and not $k_1^-$) contributes in the asymptotics \eqref{ux-fin} (via the residue calculation \eqref{ux-lbm-break-1}). 
\end{remark}

\section{Uniform in $x$ long-time asymptotics for the linear KdV equation with a sinusoidal boundary condition}
\label{kdv-asymp-s}

Section \ref{third-ss} derives the D-N map and the asymptotic solution to
the linear KdV wavemaker problem \eqref{lkdv}--\eqref{lkdv-sine-bc}
for $t\to\infty$ and a fixed $x > 0$. Herein, we study the KdV
wavemaker problem in a different asymptotic limit where $t\to\infty$
while $\xi=x/t=\mathcal{O}(1)$ is fixed, i.e., we consider the solution
along rays in the $x$-$t$ plane. This limit provides information about
how waves propagate away from the origin in the wavemaker problem.

\begin{figure}[tb!]
    \centering
        \sidesubfloat[]{\includegraphics[width=0.35\linewidth]{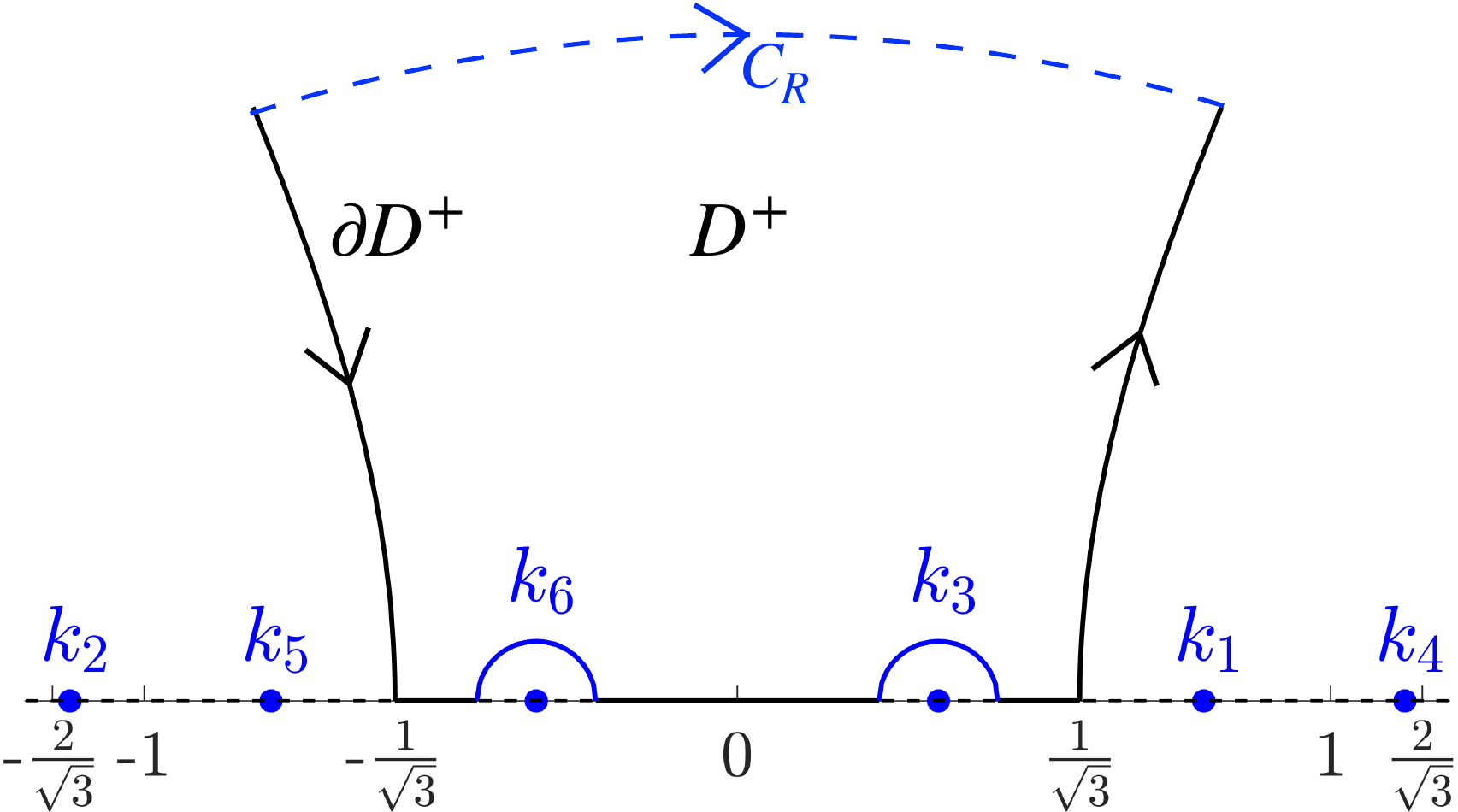} \label{fig:KdV_pole_a} }
        \sidesubfloat[]{\includegraphics[width=0.35\linewidth]{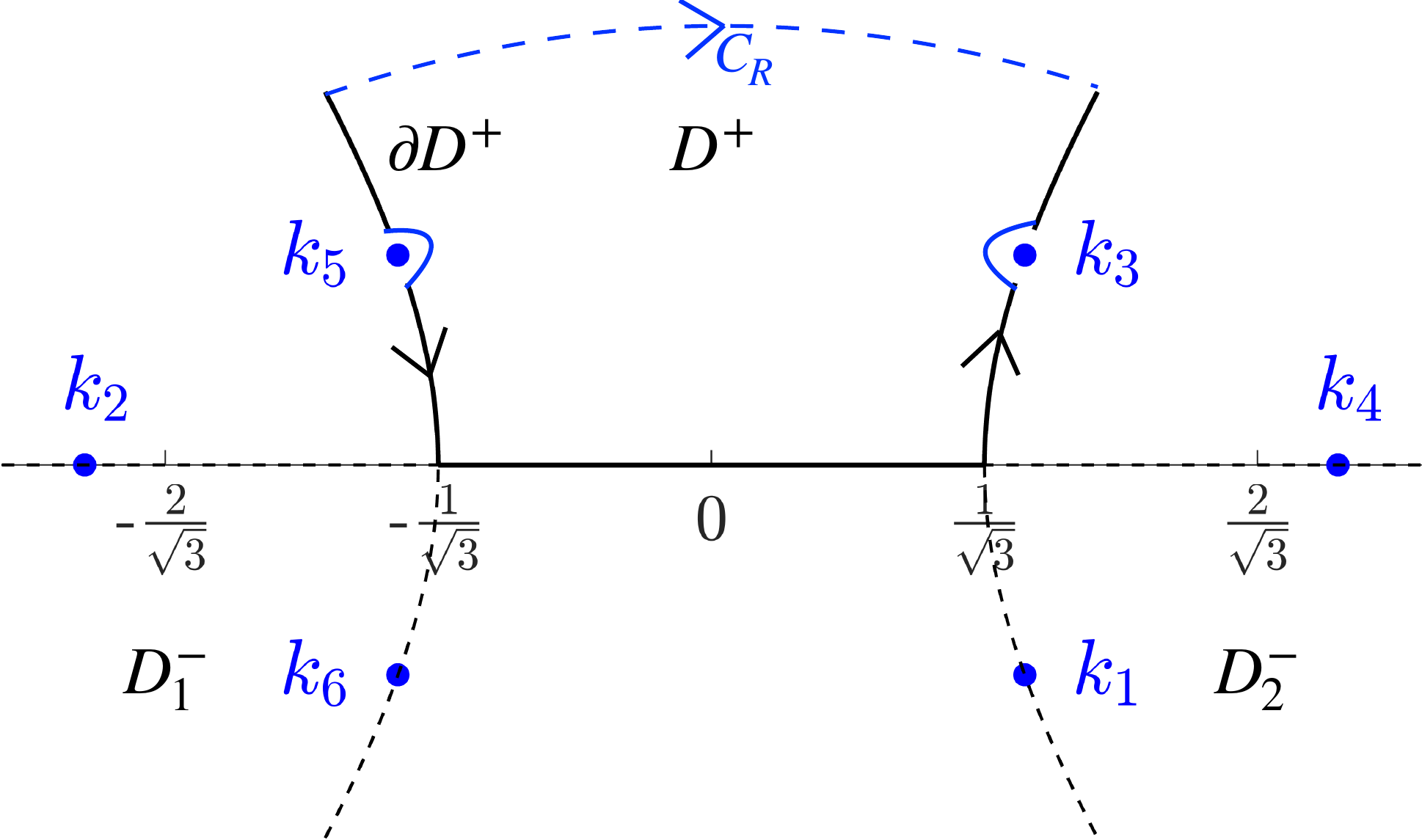} \label{fig:KdV_pole_b} }
        \caption{The positions of the simple poles $k_{1,\dots,6}$ in
          \eqref{eq:KdV_IBVP_soln} given by $\omega=\pm\omega_0$ for
          (a) the subcritical regime $0<\omega_0<2/(3\sqrt{3})$ and
          (b) the supercritical regime $\omega_0>2/(3\sqrt{3})$,
          relative to the contour $\p D^+$ for the linear KdV
          equation. The contour is indented by small semi-circles to
          avoid the poles and closed from the top with the
          semicircular arc $C_R$. }
    \label{fig:KdV_poles}
\end{figure}

Recall the wavemaker initial-boundary value problem
\eqref{lkdv-sine-bc} in which $u(x,0) = 0$, $x \ge 0$ and
$u(0,t) = -\sin(\omega_0 t)$, $t \ge 0$ where $\omega_0 > 0$ is the
input frequency.  To obtain the long-time asymptotics, we use the
solution formula for the linear KdV wavemaker problem
\eqref{lkdv}--\eqref{lkdv-sine-bc} given by \eqref{kdv-sol-wm-x},
where the contour of integration $\p D^+$ is shown in Figure
\ref{lkdv-f}.  Note that the integrand in \eqref{kdv-sol-wm-x} has
removable singularities at the six simple zeros $k_1, \ldots, k_6$ of
the two cubic equations $\Omega+n\omega_0=0$, $n=\pm1$ where $\Omega$
is defined in \eqref{om-def-cubic}.
For $0 < \omega_0 < \frac{2}{3\sqrt{3}}$, we have  
\begin{equation}\label{eq:KdV-poles-sub}
\begin{aligned}
    k_1 &= \frac{2}{\sqrt{3}}\cos \left(\frac{\alpha}{3} \right) \in \left(\frac{1}{\sqrt{3}},1\right),  \quad
    k_2 = -\frac{1}{\sqrt{3}}\cos \left(\frac{\alpha}{3} \right) - \sin \left( \frac{\alpha}{3} \right) \in \left(-\frac{2}{\sqrt{3}}, -1\right), \\
    k_3 &= -\frac{1}{\sqrt{3}}\cos \left(\frac{\alpha}{3} \right) + \sin \left( \frac{\alpha}{3} \right) \in \left(0, \frac{1}{\sqrt{3}}\right), \quad
    k_4 =  \frac{2}{\sqrt{3}}\cos \left(\frac{\theta}{3} \right) \in \left(1, \frac{2}{\sqrt{3}}\right), \\
    k_5 &= -\frac{1}{\sqrt{3}}\cos \left(\frac{\theta}{3} \right) - \sin \left( \frac{\theta}{3} \right) \in \left(-1, -\frac{1}{\sqrt{3}}\right), \quad
    k_6 = -\frac{1}{\sqrt{3}}\cos \left(\frac{\theta}{3} \right) + \sin \left( \frac{\theta}{3} \right) \in \left( -\frac{1}{\sqrt{3}}, 0\right),
\end{aligned}
\end{equation}
where
$\alpha = \pi - \tan^{-1}\left( \frac{\sqrt{12-81\omega_0^2}}{9\omega_0}
\right) \in \left(\frac{\pi}{2}, \pi \right)$ and
$\theta = \tan^{-1}\left( \frac{\sqrt{12-81\omega_0^2}}{9\omega_0}
\right)=\pi-\alpha \in \left(0, \frac{\pi}{2} \right)$.  Because all
the zeros are real as shown in Figure \ref{fig:KdV_poles}(a), we call this the subcritical regime. 
On the other hand, in the supercritical regime $\omega_0 > \frac{2}{3\sqrt{3}}$  we have
\begin{equation}\label{eq:KdV-poles-sup}
\begin{aligned}
    k_1 &= \frac{1}{2^{2/3} 3^{1/3} r_1^{1/3}} + \frac{r_1^{1/3}}{2^{4/3}3^{2/3}} + i \left(-\frac{3^{1/6}}{2^{2/3}r_1^{1/3}} + \frac{r_1^{1/3}}{2^{4/3}3^{1/6}} \right), \quad
    k_2 = -\frac{2^{1/3}}{3^{1/3}r_1^{1/3}} - \frac{r_1^{1/3}}{2^{1/3}3^{2/3}}, \\
    k_3 &= \frac{1}{2^{2/3} 3^{1/3} r_1^{1/3}} + \frac{r_1^{1/3}}{2^{4/3}3^{2/3}} + i \left(\frac{3^{1/6}}{2^{2/3}r_1^{1/3}} - \frac{r_1^{1/3}}{2^{4/3}3^{1/6}} \right), \quad
    k_4 = \frac{2^{1/3}}{3^{1/3}r_2^{1/3}} + \frac{r_2^{1/3}}{2^{1/3}3^{2/3}}, \\
    k_5 &= -\frac{1}{2^{2/3} 3^{1/3} r_2^{1/3}} - \frac{r_2^{1/3}}{2^{4/3}3^{2/3}} + i \left(-\frac{3^{1/6}}{2^{2/3}r_2^{1/3}} + \frac{r_2^{1/3}}{2^{4/3}3^{1/6}} \right),
    \\
        k_6 &= -\frac{1}{2^{2/3} 3^{1/3} r_2^{1/3}} - \frac{r_2^{1/3}}{2^{4/3}3^{2/3}} + i \left(\frac{3^{1/6}}{2^{2/3}r_2^{1/3}} - \frac{r_2^{1/3}}{2^{4/3}3^{1/6}} \right),
\end{aligned}
\end{equation}
where $r_1= 9\omega_0 - \sqrt{81\omega_0^2- 12})$ and
$r_2=9\omega_0 + \sqrt{81\omega_0^2-12})$ are real numbers.  In this case, the
location of the zeros is shown in Figure \ref{fig:KdV_poles}(b) and
described as follows: $k_1$ is complex and on the boundary of
$D_2^-$, $k_2 \in (-\infty,-2/\sqrt{3})$ is real,
$k_3$ is complex and on the right boundary of $D^+$,
$k_4 \in (2/\sqrt{3},\infty)$ is real, $k_5$ is
complex and on the left boundary of $D^+$, and $k_6$ is also
complex and on the boundary of $D_1^-$. The poles are symmetric with
respect to the imaginary and real axes and lie on the ellipse
$\frac{\mathrm{Re}(k)^2}{a^2} + \frac{3}{4}
\frac{\mathrm{Im}(k)^2}{b^2}=1$, where $a=|k_2|=|k_4|$ and
$b=|\text{Im}(k_1)|= |\text{Im}(k_3)|=
|\text{Im}(k_5)|= |\text{Im}(k_6)|$. In addition,
$k_{1,3,5,6}$ are also on the hyperbola
$3\text{Re}(k)^2 - \text{Im}(k)^2 = 1$, $j=1,3,5,6$,
given as $\p D^+$.  We define the critical frequency to be
\begin{equation}
  \label{eq:kdv_wcrit}
  \omega_0 = \omega_{cr} = \frac{2}{3\sqrt{3}} .
\end{equation}

Importantly, once the two $t$-exponentials in the solution formula \eqref{kdv-sol-wm-x} are separated, the removable singularities $k_1, \ldots, k_6$ turn into simple poles for the respective integrands. Therefore, in order to avoid crossing these poles, prior to splitting \eqref{kdv-sol-wm-x} we first indent $\p D^+$ so that the poles are excluded from the interior of $\p D^+$, as shown in Figure \ref{fig:KdV_poles}. 
In what follows, $\p D^+$ denotes the indented contour shown in Figure \ref{fig:KdV_poles}.
Then, by substituting the linear dispersion relation
$\Omega=k-k^3$ and expressing the complex exponentials as
trigonometric functions, we are able to rewrite the solution
\eqref{kdv-sol-wm-x} as
\begin{subequations}\label{eq:KdV_IBVP_soln}
\begin{align}
    u(x,t) & = u_1(x,t) + u_2(x,t),     \label{eq:KdV_IBVP_soln_2} \\
    u_1(x,t)  &= \frac{1}{2\pi}    \int_{k\in \partial D^+} \,
                f(k) e^{ikx}  \, dk , \quad f(k):= \frac{1-3k^2}{(k-k^3)^2-\omega_0^2}\left[-\omega_0 \cos(\omega_0 t)  + i(k-k^3)\sin(\omega_0 t) \right] \label{eq:KdV_IBVP_soln_2_u12}\\
    u_2(x,t) &= \frac{\omega_0}{2\pi} \int_{k\in \partial D^+}
               \, g(k) e^{ikx-i(k-k^3)t} \, dk , \quad g(k) := \frac{1-3k^2}{(k-k^3)^2-\omega_0^2} . \label{eq:KdV_IBVP_soln_2_u3}
\end{align}
\end{subequations}
Note that the dependence of $f$ on $t$ has been suppressed.   

We begin by observing that, since all of the poles $k_1, \ldots, k_6$
lie outside the contour $\p D^+$, by Cauchy's theorem inside the
region enclosed by $\p D^+$ and $C_R$ (the semicircular arcs of radius
$R$ in Fig.~\ref{fig:KdV_poles}), we have
$u_1(x, t) = \frac{1}{2\pi} \int_{k\in C_0} \, f(k) e^{ikx} \, dk$. To
evaluate along the contour $C_R$, we first notice that the integrand
$f(k)$ in \eqref{eq:KdV_IBVP_soln_2_u12} is uniformly bounded in $|k|$
($f(k)= \mathcal{O}(k^{-1})$, $k \to \infty$). Since $x\geq0$, we
apply Jordan's lemma in the upper half plane by sending $R \to \infty$
and conclude that
\begin{equation} 
u_1(x, t) \equiv 0 \quad \Rightarrow \quad u(x, t) = u_2(x, t).
\end{equation} 
The remaining term $u_2(x,t)$ defined by \eqref{eq:KdV_IBVP_soln_2_u3}
will be asymptotically evaluated in the large-time limit by the method
of steepest descent. The difference is that $e^{ikx}$ in $u_1(x,t)$
decays in $D^+$ so that one can apply Jordan's lemma along $C_R$,
whereas $u_2(x,t)$ involves $e^{ikx-i(k-k^3)t}$ which does not
necessarily decay in $D^+$.

We now asymptotically approximate $u_2(x,t)$ in
\eqref{eq:KdV_IBVP_soln_2_u3} using the method of steepest
descent. Define $\xi=x/t$ and the phase $\phi(k)=ik\xi-i(k-k^3)$, then
$u_2 = \frac{\omega_0}{2\pi} \int_{\p D^+}g(k) e^{\phi(k)t}dk$. The
saddle points are found by requiring $\phi'(\rho_j)=0$ so that
\begin{equation}
  \rho_1 = \sqrt{\dfrac{1-\xi}{3}}, \quad \rho_2 = -\sqrt{\dfrac{1-\xi}{3}}. 
  \label{eq:LKdV_saddles}
\end{equation}
The saddle points $k_{1,2}$ are real for $0\leq \xi\leq 1$ and
imaginary for $\xi>1$. Let $k=p+iq$, $p,q \in \mathbb{R}$, so that
$\phi(k)= q-3p^2q+q^3-q\xi + i(-p+p^3-3pq^2+p\xi)$. We now deform the
original contour $\partial D^+$ onto the steepest descent paths that
are defined by $\text{Im}\{\phi(k)\}=\text{Im}\{\phi(\rho_j)\}$ and
$\text{Re}\{\phi(k)-\phi(\rho_j)\}<0$, $j=1,2$.

\begin{figure}[tb!]
    \centering
    \begin{minipage}{.5\textwidth}
        \centering
        \sidesubfloat[]{\includegraphics[width=0.8\linewidth]{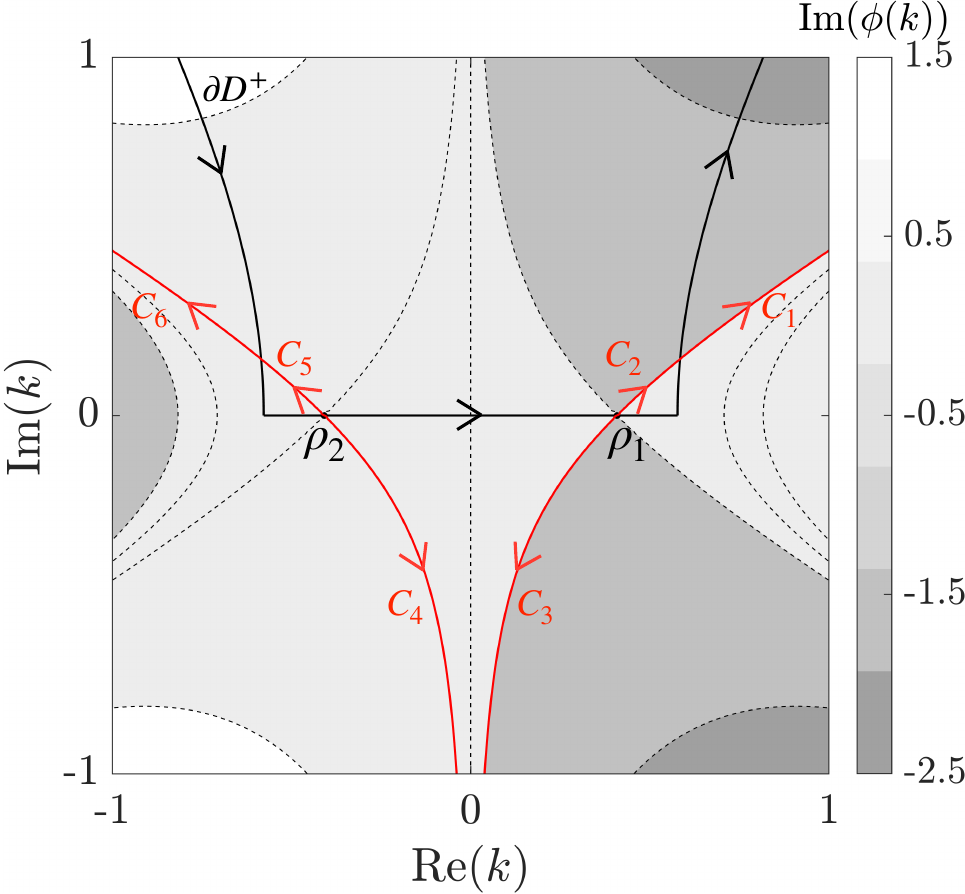} \label{fig:LKdV_contour_sub}}
    \end{minipage}%
    \begin{minipage}{0.5\textwidth}
        \centering
        \sidesubfloat[]{\includegraphics[width=0.8\linewidth]{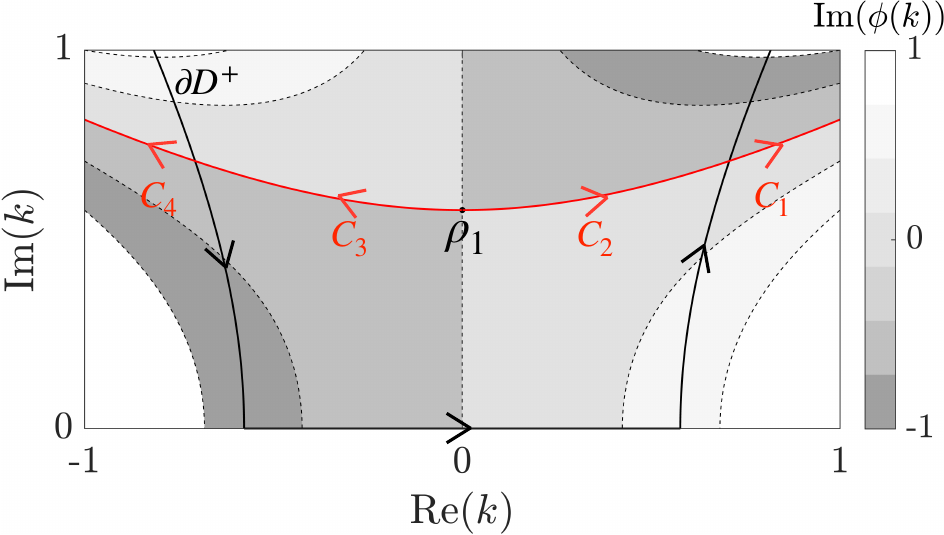} \label{fig:LKdV_contour_super} } 
    \end{minipage}
    \caption{Illustration of the saddles $\rho_{j}$ from
      \eqref{eq:LKdV_saddles} and the steepest descent contour (red)
      for the solution \eqref{eq:KdV_IBVP_soln_2_u3} of the linear KdV
      wavemaker problem in the case of (a) $0<\xi<1$ and (b) $\xi>1$,
      against the original contour $\p D^+$ (black). The contours
      $C_1 \cup C_2 \cup C_3$ and $C_4 \cup C_5 \cup C_6$ in (a) are
      given by Equations \eqref{eq:LKdV_contour_xi<1_1} and
      \eqref{eq:LKdV_contour_xi<1_2} respectively, and the contours
      $C_1 \cup C_2 \cup C_3 \cup C_4$ in (b) are described by
      \eqref{eq:LKdV_contour_xi>1}. The color map indicates levels of
      $\text{Im}(\phi(k))$ with $\text{Im}(\phi(k))$ constant along
      the steepest descent path. Recall that the contour $\p D^+$
      avoids the poles $k_1, \ldots, k_6$ as shown in Figure
      \ref{fig:KdV_poles}.}
\end{figure}

\begin{enumerate}[label=(\alph*), wide, labelwidth=!, labelindent=0pt]

\item Regime $0 < \xi < 1$.

  For the saddles $0 < \rho_1 < \frac{1}{\sqrt{3}}$ and
  $-\frac{1}{\sqrt{3}} < \rho_2 < 0$, we have the steepest descent
  contour shown in Figure \ref{fig:LKdV_contour_sub} in red with local
  steepest descent directions $\frac{\pi}{4},\frac{5\pi}{4}$ from the
  saddle $\rho_1$ and $\frac{3\pi}{4},\frac{7\pi}{4}$ from the saddle
  $\rho_2$. With $\phi(\rho_1) = i\frac{2}{9}\sqrt{3-3\xi} (-1+\xi)$
  and $\phi(\rho_2) = -i\frac{2}{9}\sqrt{3-3\xi} (-1+\xi)$, the two
  curves are
  \begin{subequations}
    \begin{align}
      \text{Im}(\phi(k)) = \text{Im}(\phi(\rho_1))
      &\implies -p+p^3-3pq^2+p\xi = \frac{2}{9}\sqrt{3-3\xi}
        (-1+\xi), \label{eq:LKdV_contour_xi<1_1}
      \\ 
      \text{Im}(\phi(k)) = \text{Im}(\phi(\rho_2))
      &\implies -p+p^3-3pq^2+p\xi = -\frac{2}{9}\sqrt{3-3\xi}
        (-1+\xi), \label{eq:LKdV_contour_xi<1_2}
    \end{align}
  \end{subequations}
  respectively. Note that the exponential term satisfies
  $|e^{\phi(k)t}| = 1$ on the steepest descent contour for
  $0 < \xi < 1$ so it is purely oscillatory.  On the right steepest
  descent path $C_1 \cup C_2 \cup C_3$ in Figure
  \ref{fig:LKdV_contour_sub}, the dominant contribution comes from the
  neighborhood of $\rho_1$. Parameterizing this neighborhood by
  $k=\rho_1 + (1+i)s $ with small $s \in (-\epsilon,\epsilon)$, we
  obtain
  \begin{align}\label{eq:KdV-saddle-cont1}
    \int_{C_1+C_2-C_3} g(k) e^{\phi(k)t} \, dk 
    &\sim \int_{-\epsilon}^{\epsilon}
      g(\rho_1 + (1+i)s) e^{t\left( \phi(\rho_1) +
      \frac{1}{2}(1+i)^2s^2 \phi''(\rho_1) \right)} \, ds \nonumber\\ 
    & \sim e^{t\phi(\rho_1)} g(\rho_1)
      \int_{-\infty}^{\infty} e^{is^2 \phi''(\rho_1)} \, ds, \quad
      \phi''(\rho_1)=2i\sqrt{3(1-\xi)},    \nonumber\\
    & = \frac{e^{ - i\frac{2}{9}\sqrt{3(1-\xi)}(1-\xi) t}}{\sqrt{t}} \frac{27 \xi
      }{4-(\xi +3) \xi ^2-27 \omega_0^2} 
      \sqrt{\frac{\pi}{2\sqrt{3(1-\xi)}}}, \quad t \to \infty .
  \end{align}
  Similarly, on the left steepest descent path $C_4 \cup C_5 \cup C_6$
  in Figure \ref{fig:LKdV_contour_sub}, we parameterize
  $k=\rho_2 + (1-i)s$, $s \in (-\epsilon,\epsilon)$. The dominant
  contribution is the complex conjugate of \eqref{eq:KdV-saddle-cont1}
  \begin{align}\label{eq:KdV-saddle-cont2}
    \int_{C_4-C_5-C_6} g(k) e^{\phi(k) t} \, dk 
    & \sim \frac{1}{\sqrt{t}} e^{\frac{2}{9}i\sqrt{3(1-\xi)}(1-\xi)t} \frac{27 \xi
      }{-(\xi +3) \xi ^2-27 \omega_0^2+4}
      \sqrt{\frac{\pi}{2\sqrt{3(1-\xi)}}},  \quad t \to \infty. 
  \end{align}
  We now bring in the contribution from the poles $k_{1,\dots,6}$ by
  considering different cases for the input frequency $\omega_0$.
    
  \begin{enumerate}[label=(a.\arabic*), wide, labelwidth=!, labelindent=0pt]
    
    \bigskip
  \item Subcritical regime:  $0<\omega_0<\frac{2}{3\sqrt{3}}$.
    
    The two possible scenarios in this case correspond to Regions I
    and II.a in Figure \ref{fig:LKdV_soln_descent}. In particular,
    with both $\rho_1,k_3 \in [0,\frac{1}{\sqrt{3}}]$ and
    $\rho_2,k_6 \in [-\frac{1}{\sqrt{3}},0]$, their relative ordering
    is important.  The boundary between Regions I and II.a is defined
    by the coalescence $k_3=\rho_1$ and $k_6=\rho_2$ so that
    \begin{equation}
      \label{eq:kdv_grp_velocity_curve}
      \xi = 1-3k_3^2 =\omega'(k_3), \quad
      0<\omega_0<2/(3\sqrt{3}) ~\text{and}~ 0<\xi<1.
    \end{equation}
    We are now in a position to deform the original contour
    $\partial D^+$ onto the steepest descent paths.
    There are two cases to consider. First, the region
    $0 < \xi<\omega'(k_3)$, which corresponds to the poles lying
    between the saddle points as shown in Region I of Figure
    \ref{fig:LKdV_soln_descent}. Notice that the integrand
    $g(k) e^{\phi(k)t}$ of \eqref{eq:KdV_IBVP_soln_2_u3} decays for
    $k\to\infty$ in the region between $D_1$ and $C_1$ (also $D_5$ and
    $C_6$) with
    $ \left| e^{\phi(k)t} \right| = e^{-\text{Im}(k)\xi t +
      \text{Im}(k-k^3)t}$, $\xi>0,t>0$, where $\text{Im}(k)>0$ and
    $\text{Im}(k-k^3) < 0$ in this region. Thus, one can simply deform
    $D_1$ onto $C_1$ and $D_5$ onto $C_6$ by Cauchy's theorem and
    Jordan's lemma so that
    $\int_{D_1}g(k) e^{\phi(k)t}\,dk=\int_{C_1} g(k) e^{\phi(k)t}\,dk$
    and
    $\int_{D_5}g(k) e^{\phi(k)t}\,dk=-\int_{C_6}g(k)
    e^{\phi(k)t}\,dk$.  The absence of poles between $D_2$ and $C_2$
    ($D_4$ and $C_5$) and Cauchy's theorem imply
    $\int_{D_2} g(k) e^{\phi(k)t}\,dk=\int_{C_2}g(k) e^{\phi(k)t}\,dk$
    and
    $\int_{D_4} g(k) e^{\phi(k)t}\,dk=-\int_{C_5} g(k)
    e^{\phi(k)t}\,dk$. Finally, since the poles $k_3$ and $k_6$ lie
    inside the region enclosed by $D_3$, $C_3$ and $C_4$, using
    Cauchy's residue theorem we compute
    \begin{align*}
      \int_{D_3+C_3-C_4} g(k) e^{\phi(k)t}\,dk
      &= -2\pi i \left( \text{Res}(k_3) + \text{Res}(k_6) \right), \\
      \text{Res}(k_3) = \frac{1}{2\omega_0} e^{i(k_3 \xi t - \omega_0
      t)},
      &\quad \text{Res}(k_6) = -\text{Res}(k_3)^*, 
    \end{align*}
    where $k_3=-k_6>0$ satisfies the linear dispersion relation
    $k_3-k_3^3=\omega_0$. Combining these results, we obtain
    $u_2(\xi;t) = \sin \left[(k_3 \xi -\omega_0) t \right] +
    \frac{\omega_0}{2\pi} \int_{C_1+C_2-C_3+C_4-C_5-C_6} g(k)
    e^{\phi(k)t}\,dk$.  The remaining integral on the right-hand side
    corresponds to the saddle point contributions
    \eqref{eq:KdV-saddle-cont1} and \eqref{eq:KdV-saddle-cont2}, which
    are dominated by the pole contributions. Therefore, the large $t$
    solution to the wavemaker problem is
    \begin{equation}
      u(\xi t,t) = \sin \left[ (k_0 \xi -\omega_0 ) t \right] +
      \mathcal{O}(t^{-1/2}), \quad t \to \infty,
      \label{eq:LKdV_soln_plane} 
    \end{equation}
    where $k_0 \equiv k_3$ in this case, which satisfies the radiation
    condition \eqref{eq:general_radiation_condition}.
        
    In contrast, the case $\omega'(k_0) < \xi < 1$ occurs when the
    saddles are in between the poles depicted in Region II.a of Figure
    \ref{fig:LKdV_soln_descent}. Using Cauchy's theorem, we can deform
    $\p D^+$ to the steepest descent contour $\cup_{j=1}^6 C_j$ with
    no pole contributions.     Hence, we obtain
    $u_2(x,t)= \frac{\omega_0}{2\pi} \int_{C_1+C_2-C_3+C_4-C_5-C_6}
    g(k) e^{\phi(k)t} dk$ so the
    saddle point contributions \eqref{eq:KdV-saddle-cont1} and
    \eqref{eq:KdV-saddle-cont2} are asymptotically dominant
    \begin{equation}
      u(\xi t,t) \sim \frac{\cos\left( \frac{2}{9} t 
          \sqrt{3(1-\xi)} (1-\xi) \right) }{\sqrt{t}} \frac{27 \xi
        \omega_0 }{4-(\xi +3) \xi ^2-27 
        \omega_0^2}  \sqrt{\frac{1}{2\pi\sqrt{3(1-\xi)}}}, \quad t
      \to \infty. \label{eq:LKdV_xi<1_decay_soln}
    \end{equation}
        
    In summary, when $0 < \omega_0 < 2/(3\sqrt{3})$, the leading order
    solution \eqref{eq:LKdV_soln_plane} is a traveling wave with real
    wavenumber $k_3$ and frequency $\omega_0$ satisfying the
    dispersion relation $\omega_0 = k_3 - k_3^3$ as $t \to \infty$
    when $0 < \xi = x/t < \omega'(k_3)$. This regime is dominated by
    the pole contributions from contour deformation.  When
    $x/t > \omega'(k_3)$, i.e., ahead of the boundary-generated,
    invading wavefront, the leading order solution
    \eqref{eq:LKdV_xi<1_decay_soln} exhibits algebraic decay in time
    with $u(x,t)=\mathcal{O}(t^{-1/2})$, $\xi = x/t = \mathcal{O}(1)$,
    $t \to \infty$.  Here, the solution is dominated by saddle point
    contributions.  In contrast to the previous case
    \eqref{eq:LKdV_soln_plane} when $0 < \xi < \omega'(k_0)$, the
    asymptotic solution \eqref{eq:LKdV_xi<1_decay_soln} is not time
    periodic.  
    
    Note that the obtained asymptotic approximations are not uniform in $\xi$.  For example, eq.~\eqref{eq:LKdV_xi<1_decay_soln} is singular when $\xi = \omega'(k_3)$.  The cause of nonuniformity is the collisions of the poles $k_{3,6}$ with the saddle points $\rho_{1,2}$.  Uniform asymptotic approximations involving the error function with complex argument can be obtained \cite{clemmow_extensions_1950,van_der_waerden_method_1952,felsen_asymptotic_1994}.
        
    \bigskip
  \item Supercritical regime: $\omega_0>2/(3\sqrt{3})$.
    
    In this case, some poles are complex, including $k_{3,5}$, which
    are avoided by semicircle indentations of the contour $\p D^+$
    (c.f.~Fig.~\ref{fig:KdV_pole_a}).  In order to determine the pole
    locations relative to the steepest descent contours $C_j$,
    $j = 1,\ldots,6$ depicted in red in Figure
    \ref{fig:LKdV_contour_sub} (case $0<\xi<1$), notice that $C_{1,2}$
    and $C_{5,6}$ intersect with $\p D^+$ at
    \begin{align*}
      (p_1,q_1) &= \left( \sqrt{\frac{2+\xi}{6}} \cos
                  \frac{\varphi}{3} , \sqrt{-1 + \frac{2+\xi}{2}
                  \cos^2 \frac{\varphi}{3} } \right), \quad 
                  (p_2,q_2) = \left( -\sqrt{\frac{2+\xi}{6}} \cos
                  \frac{\varphi}{3}, \sqrt{-1 + \frac{2+\xi}{2} \cos^2
                  \frac{\varphi}{3} }  \right),
    \end{align*}
    where
    $\varphi= \arctan \left[ \frac{3 \sqrt{\xi(4-2\xi+\xi^2)}}
      {2\sqrt{2} (\xi-1)\sqrt{1-\xi}} \right] $. Here, $(p_1,q_1)$ is
    associated with $C_{1,2} \cap \p D^+$ and $(p_2,q_2)$ is
    associated with $C_{5,6} \cap \p D^+$. We now determine for what
    value of $\xi$ the poles $k_{3,5}$ coincide with
    $(p_{1,2},q_{1,2})$.  Letting
    $\text{Re}(k_3) = p_1, \text{Re}(k_5) = p_2$ and simplifying,
    results in
    \begin{equation}
      l_1(\xi,\omega_0) := \frac{1}{2^{2/3} 3^{1/3} r_1^{1/3}} +
      \frac{r_1^{1/3}}{2^{4/3}3^{2/3}} -  \sqrt{\frac{2+\xi}{6}} \cos
      \frac{\varphi}{3}=0, \quad \omega_0>2/(3\sqrt{3}) \text{ and }
      0<\xi<1,  \label{eq:LKdV_l1}
    \end{equation}
    where $r_1$ is defined below \eqref{eq:KdV-poles-sup}. Equation
    \eqref{eq:LKdV_l1} determines the curve in the
    $\omega_0$-$\xi$ plane that splits Regions II.b and III in
    Figure \ref{fig:LKdV_soln_descent}.
    
    When $l_1(\xi,\omega_0)<0$, the poles $k_{3,5}$ are avoided by the
    corner contours $D_{2,4}$ with $\text{Im}(k_{3,5})>0$ (Region II.b
    in Figure \ref{fig:LKdV_soln_descent}). Using Cauchy's theorem and
    Jordan's lemma where appropriate, we deform each sector of
    $\p D^+$ onto the steepest descent path with no pole
    contributions.  This gives
    $u_2(\xi;t)=\frac{\omega_0}{2\pi} \int_{C_1+C_2-C_3+C_4-C_5-C_6}
    g(k) e^{\phi(k)t}\,dk$ so that the leading-order behavior of
    $u(\xi t, t)$ is determined by the saddle contributions as in
    \eqref{eq:LKdV_xi<1_decay_soln}, with algebraic decay
    $\mathcal{O}(t^{-1/2})$.
        
    When $l_1(\xi,\omega_0)>0$ (see Region III in Figure
    \ref{fig:LKdV_soln_descent}), the poles $k_3$ and $k_5$ lie
    between the contours $D_1$, $C_1$ and $D_5$, $C_6$
    respectively. Thus, when invoking Cauchy's theorem, the residue
    contributions from these poles must be taken into
    account. Computing
    \begin{equation}\label{eq:KdV-res-sup}
      \text{Res}(k_3) = \frac{1}{2\omega_0} e^{i\left(k_3  \xi t -
          \omega_0 t \right)}, \quad \text{Res}(k_5) =
      -\text{Res}(k_3)^* 
    \end{equation}
    because $k_3 = -k_5^*$.  The pole $k_3$ corresponds to the
    singularity on $\p D^+$ that determines the D-N map $k_0 := k_3$
    and satisfies the dispersion relation $k_3 - k_3^3 = \omega_0$.
    We find
    $u_2(\xi t,t)= e^{-\mathrm{Im}(k_0) \xi t} \sin \left[
      (\mathrm{Re}(k_0) \xi - \omega_0 t) \right] +
    \frac{\omega_0}{2\pi} \int_{C_1+C_2-C_3+C_4-C_5-C_6} g(k)
    e^{\phi(k)t}\,dk$. In light of \eqref{eq:KdV-saddle-cont1} and
    \eqref{eq:KdV-saddle-cont2}, we deduce
    \begin{equation}
      u(\xi t,t) \sim  e^{-\mathrm{Im}(k_0) \xi t} \sin \left[
        (\mathrm{Re}(k_0) \xi - \omega_0) t \right] +
      \frac{1}{\sqrt{t}} \frac{27 \xi
        \omega_0 \cos\left( \frac{2}{9} t \sqrt{3(1-\xi)} (1-\xi)
        \right) }{-(\xi +3) \xi ^2-27 \omega_0^2+4}
      \sqrt{\frac{1}{2\pi\sqrt{3(1-\xi)}}}, \quad t \to
      \infty. \label{eq:LKdV_xi<1_super_soln}  
    \end{equation}
    The second term in \eqref{eq:LKdV_xi<1_super_soln} is the dominant
    contribution, with $\mathcal{O}(t^{-1/2})$ for
    $\xi = \mathcal{O}(1)$.  The first term is beyond all orders. In
    other words, the asymptotic behavior of the solution in Region III
    is the same as the asymptotic behavior
    \eqref{eq:LKdV_xi<1_decay_soln} in Regions II.a and II.b.
    However, as we now analyze, both terms in
    \eqref{eq:LKdV_xi<1_super_soln} are required for obtaining a
    uniform asymptotic expansion that is also valid as $\xi\to0^+$,
    i.e., $x$ fixed and $t \to \infty$.
        
    Rewriting \eqref{eq:LKdV_xi<1_super_soln} in terms of $x=\xi t$
    and $t$, we obtain
    \begin{equation*}
      \begin{split}
        u(x,t) \sim &e^{-\text{Im}(k_0) x} \sin(\mathrm{Re}(k_0)x -
                      \omega_0t) \\
                    &+ \frac{1}{t^{3/2}} \frac{27 \omega_0 x }{-(x/t +3)
                      (x/t) ^2-27 \omega_0^2+4} \cos\left[ \frac{2}{9} t
                      \sqrt{3(1-x/t)} (1-x/t) \right] \sqrt{\frac{1}{2\pi
                      \sqrt{3(1-x/t)}}},
      \end{split}
    \end{equation*}
    for $x \in [0,L]$, some $L > 0$ and $t \to \infty$.  This
    expression is asymptotically equivalent to
    $$u(x,t)= e^{-\text{Im}(k_0) x}
    \sin(\mathrm{Re}(k_0)x - \omega_0t) + \mathcal{O}(t^{-3/2}), \quad
    t \to \infty .$$ The asymptotics of $u(x,t)$ match our earlier
    analysis using the D-N map for fixed $x$ and large $t$ given in
    Equation \eqref{eq:LKdV_series_soln}. Therefore, keeping both
    terms in \eqref{eq:LKdV_xi<1_super_soln} gives a uniform in $x$
    asymptotic expansion as $t \to \infty$.
  \end{enumerate}
  
  \begin{figure}[tb!]
    \centering
    \includegraphics[width=0.8\linewidth]{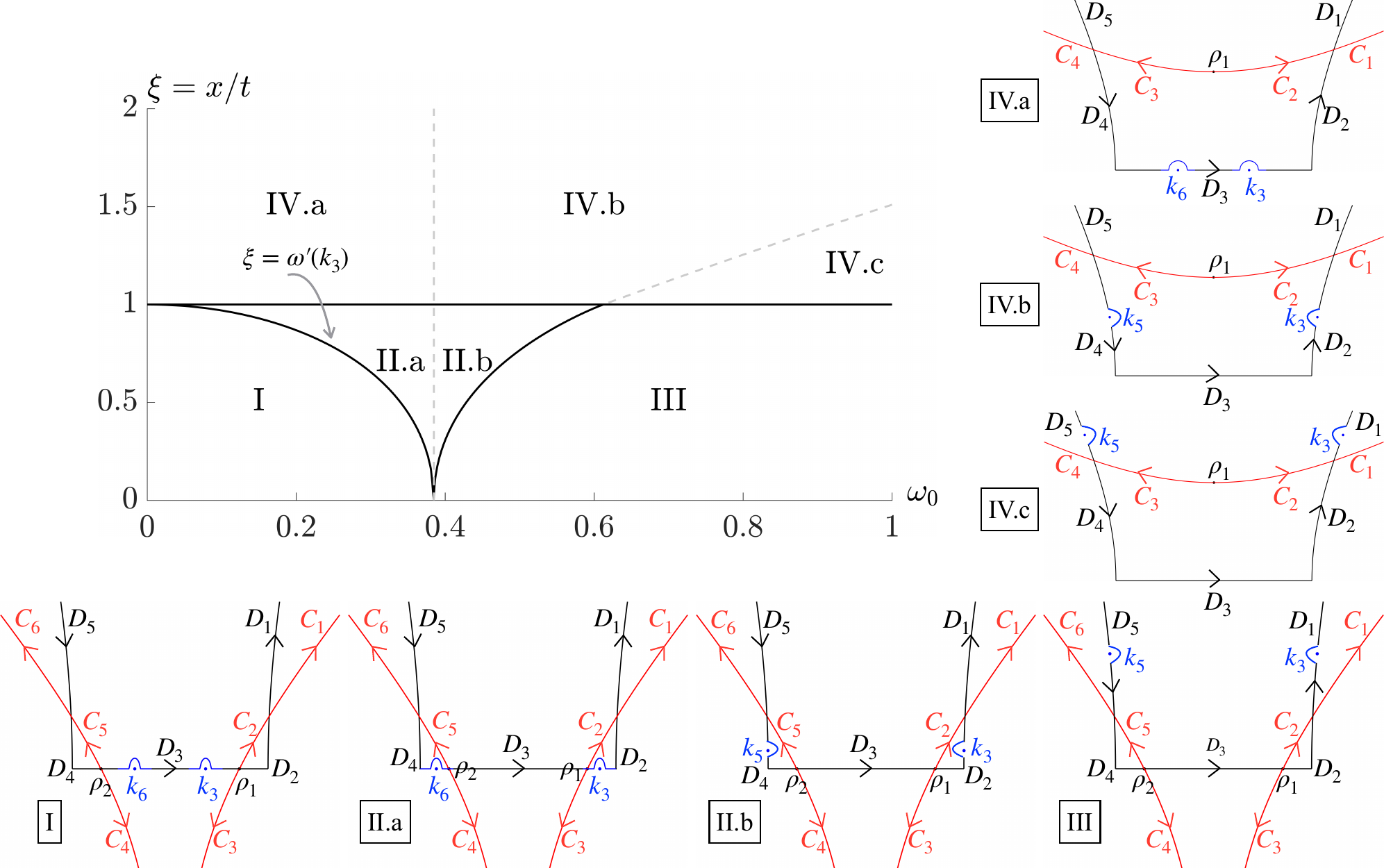}
    \caption{Phase diagram of the long-time asymptotics for the linear
      KdV equation \eqref{lkdv} subject to the sinusoidal boundary
      condition \eqref{lkdv-sine-bc}. The solution in each region is
      described in Section \ref{kdv-asymp-s}. The insets show the
      locations of the saddles and the poles relative to the original
      contour $\partial D^+ = \cup_{j=1}^5 D_j$ and the steepest
      descent contours $C_1 \cup C_2 \cup C_3$ and
      $C_4 \cup C_5 \cup C_6$ that are determined by the values of
      $\xi = x/t$ and the boundary frequency $\omega_0$.  The solid curves
      separating various regions in the $\omega_0$-$\xi$ plane are
      $\xi = 1$, Equations \eqref{eq:kdv_grp_velocity_curve} and
      \eqref{eq:LKdV_l1}.  The dashed curves are $\omega_0 =
      2/(3\sqrt{3})$ and \eqref{eq:KdV_l2}.  }
    \label{fig:LKdV_soln_descent}
  \end{figure}

  \bigskip
\item Regime $\xi > 1$.  
  
  Now the saddles are purely imaginary
  $\rho_1=i\sqrt{ \frac{\xi-1}{3}}$,
  $\rho_2=-i \sqrt{\frac{\xi-1}{3}}$ and the local steepest descent
  directions are $\theta_1 \in \{0,\pi\}$ at $\rho_1$ and
  $\theta_2\in\{\frac{\pi}{2},\frac{3\pi}{2}\}$ at $\rho_2$, because
  $\phi''(\rho_1)=- 2\sqrt{3(\xi-1)} < 0$,
  $\phi''(\rho_2)=2\sqrt{3(\xi-1)} > 0$. However, $\rho_2$ is
  disregarded since
  $\phi(\rho_2) = \frac{2}{9} \sqrt{3\xi-3} (\xi-1) > 0$ so
  $e^{t \phi(\rho_2)}$ is unbounded as $t\to \infty$. On the other
  hand, since $\phi(\rho_1) = -\frac{2}{9} \sqrt{3\xi-3} (\xi-1) < 0$,
  $e^{t \phi(\rho_1)}$ decays as $t\to\infty$. The deformed contour,
  which follows the steepest descent directions emanating from the
  saddle $\rho_1$, is shown in red in Figure
  \ref{fig:LKdV_contour_super} and is described by setting
  $k = p + i q$, $p,q \in \mathbb{R}$ and
  \begin{equation}
    \text{Im}(\phi(k)) = \text{Im}(\phi(\rho_1)) \implies
    p(p^2-1-3q^2+\xi) = 0. \label{eq:LKdV_contour_xi>1}
  \end{equation}
  Knowing that the main contribution along the contours
  $C_1, \ldots, C_4$ in Figure \ref{fig:LKdV_contour_super} comes from
  the neighborhood of the saddle $\rho_1$, we parameterize
  $k=\rho_1 + s$ with small $s \in (-\epsilon, \epsilon)$ so that,
  using the Laplace method, we find
  \begin{align}\label{eq:KdV-saddle-xi>1}
    \int_{C_1+C_2-C_3-C_4} g(k) e^{\phi(k)t} \, dk
    &\sim \int_{-\epsilon}^{\epsilon}
      \frac{(1-3\left(\rho_1+s\right)^2)\omega_0}
      {\left(\left(\rho_1+s\right) -\left(\rho_1+s\right)^3\right)^2
      -\omega_0^2} e^{t \phi(\rho_1+s)  } \, ds  \nonumber\\
    &\sim \frac{(1-3\rho_1^2)
      \omega_0}{(\rho_1-\rho_1^3)^2-\omega_0^2} e^{t\phi(\rho_1)}
      \int_{-\infty}^\infty e^{t\frac{s^2}{2}\phi''(\rho_1)} \, ds
      \nonumber\\ 
    &\sim \frac{e^{-\frac{2}{9}
      \sqrt{3\xi-3} (\xi-1) t}}{\sqrt{t}} \frac{27 \xi  \omega_0
      }{4-(\xi +3) \xi ^2-27 \omega_0^2}
      \sqrt{\frac{\pi}{\sqrt{3(\xi-1)}}}, \quad t \to \infty. 
  \end{align}
  Now we deform $\p D_+$ onto $C_1, \ldots, C_4$ taking the poles into
  consideration.
    
  \begin{enumerate}[label=(b.\arabic*), wide, labelwidth=!, labelindent=0pt]
    
    \bigskip
  \item Subcritical regime: $0<\omega_0<2/(3\sqrt{3})$.
    
    This case corresponds to Region IV.a in Figure
    \ref{fig:LKdV_soln_descent}. The poles $k_{3,6}$ are real and
    lie below $\p D^+$. Consequently, they lie outside the region
    enclosed by $C_2, C_3, D_2, D_3, D_4$, so that Cauchy's theorem
    implies
    $u_2(\xi t,t) = \frac{\omega_0}{2\pi} \int_{C_1+C_2-C_3-C_4}
    g(k) e^{\phi(k)t}\,dk$.  Using \eqref{eq:KdV-saddle-xi>1}, the
    leading-order solution in this region decays exponentially
    \begin{align}
      u(\xi t,t) 
      \sim \frac{e^{-\frac{2}{9} \sqrt{3\xi-3}
      (\xi-1) t}}{\sqrt{t}} \frac{27 \xi  \omega_0 }{8-2(\xi +3) \xi
      ^2-54 \omega_0^2} 
      \sqrt{\frac{1}{\pi \sqrt{3(\xi-1)}}}, \quad t \to
      \infty.  \label{eq:KdV_xi>1_soln_2} 
    \end{align}
    
    \bigskip   
  \item Supercritical regime: $\omega_0>2/(3\sqrt{3})$.
    
    This case corresponds to Regions IV.b and IV.c in Figure
    \ref{fig:LKdV_soln_descent}.  In particular, the poles
    $\kappa_{3,5}$ are now near the right and left edge of the contour
    $\p D^+$, respectively, with $\mathrm{Re}(k_3)=-\mathrm{Re}(k_5)$
    and $\mathrm{Im}(k_3)=\mathrm{Im}(k_5)>0$. Thus, we need to
    distinguish how they are located relative to the deformed contours
    $C_1, \ldots, C_4$. Note from Figure \ref{fig:LKdV_contour_super}
    that these contours intersect with $\p D^+$ when $k = p+iq$,
    $-p + p^3 - 3p q^2 + p\xi \Big|_{q=\sqrt{3p^2-1}} = 0$, i.e., at
    the points
    \begin{align*}
      (p_3,q_3) = \left( \frac{\sqrt{2+\xi}}{2\sqrt{2}},
      \sqrt{-1+\frac{3}{8}(2+\xi)} \right), \quad 
      (p_4,q_4) = \left( - \frac{\sqrt{2+\xi}}{2\sqrt{2}},
      \sqrt{-1+\frac{3}{8}(2+\xi)} \right). 
    \end{align*}
    The poles $k_{3,5}$ coincide with these intersection points when
    $p_3 = \text{Re}(k_3)$ or, equivalently, $p_4 = \text{Re}(k_5)$,
    i.e. when
    \begin{equation}
      l_2:=  \frac{1}{2^{2/3} 3^{1/3} r_1^{1/3}} +
      \frac{r_1^{1/3}}{2^{4/3}3^{2/3}} -
      \frac{\sqrt{2+\xi}}{2\sqrt{2}} =0,  \quad \omega_0>2/(3\sqrt{3})
      \text{ and } \xi>1,  \label{eq:KdV_l2}
    \end{equation}
    where $r_1$ is defined below \eqref{eq:KdV-poles-sup}. Equation
    \eqref{eq:KdV_l2} determines the curve in the $\omega_0$-$\xi$
    plane that splits Regions IV.b and IV.c in Figure
    \ref{fig:LKdV_soln_descent}. When $l_2(\xi,\omega_0)<0$, then the
    poles $\kappa_{3,5}$ lie outside the region enclosed by
    $C_2, C_3, D_2, D_3, D_4$ (see Region IV.b in Figure
    \ref{fig:LKdV_soln_descent}), thus
    $u_2(\xi t,t) = \frac{\omega_0}{2\pi} \int_{C_1+C_2-C_3-C_4} g(k)
    e^{\phi(k)t}\,dk $ so that the leading-order behavior of
    $u(\xi t,t)$ is given by \eqref{eq:KdV_xi>1_soln_2}. On the other
    hand, when $l_2(\xi,\omega_0)>0$ then $k_3$ and $k_5$ are inside
    the regions enclosed by $C_1, D_1$ and $C_4, D_5$ respectively
    (Region IV.c in Figure \ref{fig:LKdV_soln_descent}) so in view of
    \eqref{eq:KdV-res-sup} we have the additional residue contribution
    $u_2(\xi t,t)= e^{-\mathrm{Im}(k_3)\xi t} \sin(\mathrm{Re}(k_3) \xi
    t - \omega_0 t) + \frac{\omega_0}{2\pi} \int_{C_1+C_2-C_3-C_4}
    g(k) e^{\phi(k)t}\, dk$. Hence, using once again
    \eqref{eq:KdV-saddle-xi>1}, the leading-order asymptotics of the
    solution in this region is described by
    \begin{equation}
      \begin{split}
        u(\xi t,t)  \sim
        &e^{-\mathrm{Im}(k_3)\xi t}
          \sin(\mathrm{Re}(k_3) \xi t - \omega_0 t) \\
        &+ \frac{e^{-\frac{2}{9} \sqrt{3\xi-3}
          (\xi-1) t}}{\sqrt{t}} \frac{27 \xi\omega_0  }{8-2(\xi +3)
          \xi ^2-54 \omega_0^2} 
          \sqrt{\frac{1}{\pi \sqrt{3(\xi-1)}}}, \quad t \to
          \infty.
      \end{split}
      \label{eq:KdV_solnV}
    \end{equation}
    Observe that both terms on the right-hand side decay
    exponentially. In order to figure out which one of them is
    dominant, we notice that when $l_2>0$,
    i.e. $\mathrm{Re}(k_3)>\sqrt{(2+\xi)/8}$, we have
    $\mathrm{Im}(k_3)^2=3\mathrm{Re}(k_3)^2-1 > (3\xi-2)/8$ so that
    $\mathrm{Im}(k_3)\xi>\frac{2}{9} \sqrt{3\xi-3} (\xi-1)$.  Thus,
    \eqref{eq:KdV_solnV} is dominated by the second term.  Therefore,
    in the regime $\xi>1$, the leading-order behavior of $u(\xi t,t)$ is
    always given by \eqref{eq:KdV_xi>1_soln_2} in both the subcritical
    and the supercritical regime, i.e. regardless of the range of
    values of $\omega_0$.

    When $\xi$ passes through 1, the asymptotic expansion abruptly changes from oscillatory to exponential, e.g., \eqref{eq:LKdV_xi<1_super_soln} to \eqref{eq:KdV_solnV}.  This nonuniformity is due to the collision of the saddle points $\rho_{1,2}$, which can be cast as a uniform asymptotic expansion by use of Airy-type integrals \cite{felsen_asymptotic_1994}.
  \end{enumerate}
\end{enumerate}

The large $t$ solution given in eqs.~\eqref{eq:LKdV_soln_plane},
\eqref{eq:LKdV_xi<1_decay_soln}, \eqref{eq:LKdV_xi<1_super_soln},
\eqref{eq:KdV_xi>1_soln_2}, and \eqref{eq:KdV_solnV} proves Theorem
\ref{lkdv-t}.  We now carry out a much simpler calculation of a subset
of these results using modulation theory.

\begin{remark}[Linear modulation theory] 
  Wave modulation theory is one of the primary techniques to describe
  slow modulations of periodic waves. Herein, we use linear modulation
  theory to obtain the wave's amplitude, wavenumber and frequency for
  the linear KdV boundary value problem
  \eqref{lkdv}-\eqref{lkdv-sine-bc} and compare with the long-time
  asymptotic approximation just obtained.  The modulation equations of
  wave conservation are \cite{whitham2011linear}
  \begin{subequations}
    \label{eq:Whitham_modulation}
    \begin{align}
      k(\omega)_t + \omega_x &=0, \label{eq:Whitham_1}\\
      (a^2)_t+(\omega'(k)a^2)_x &=0, \label{eq:Whitham_2}
    \end{align}
  \end{subequations}
  where $\omega$ is the angular frequency and $k(\omega)$ is one of
  the roots of the linear dispersion relation
  $k-k^3=\omega$. Translating the initial-boundary value problem
  \eqref{lkdv}-\eqref{lkdv-sine-bc} for $u(x,t)$ to the
  initial-boundary value problem for $k(x,t)$ and $a(x,t)$, we have
  \begin{equation} \label{eq:Whitham_BC}
    \begin{split}
      \omega(x,0) &= 0, \quad a(x,0) = 0, \quad x > 0, \\
      \omega(0,t) &= \omega_0, \quad 
                    a(0,t) = 1, \quad t>0.
    \end{split}
  \end{equation}
  Equations \eqref{eq:Whitham_modulation} are a system of quasi-linear
  equations with the double characteristic speed (eigenvalue)
  $c_g(\omega_0)$ that has a one-dimensional eigenspace.  Equation
  \eqref{eq:Whitham_1} and the initial-boundary conditions
  \eqref{eq:Whitham_BC} are invariant under the hydrodynamic scaling
  symmetry $x \to b x$, $t \to b t$ for any $b > 0$.  Consequently, we
  seek a self-similar solution that preserves this symmetry:
  $\omega=\omega(\xi)$, where $\xi=x/t$.  Since
  $k(\omega)_t=\omega_t/c_g(\omega)$, either $\omega'(\xi)=0$ or
  $c_g(\omega)=\xi$, where
  $c_g(\omega) = 1/k'(\omega) = 1 - 3 k(\omega)^2$ is the group
  velocity. Combining these in a piecewise fashion, we obtain the
  solution
  \begin{equation}
    \label{eq:modulation_solution_kdv}
    \omega(\xi)= 
    \begin{cases}
      \omega_0, & 0\leq\xi\leq c_g(\omega_0)\\
      \frac{1}{9} \sqrt{3-3\xi}(2+\xi), & c_g(\omega_0)<\xi\leq1
    \end{cases} , \quad 
    k(\xi)= 
    \begin{cases}
      k_0, & 0\leq\xi\leq c_g(\omega_0)\\
      \sqrt{\frac{1-\xi}{3}}, & c_g(\omega_0)<\xi\leq1
    \end{cases},
  \end{equation}
  where $k_0 = k(\omega_0)$ satisfies the linear dispersion relation
  $k_0-k_0^3=\omega_0$.

  Here we encounter the issue of selecting $k(\omega_0)$ from the
  three possible wavenumbers that satisfy the cubic polynomial linear
  dispersion relation.  In the context of modulation theory, there are
  two admissibility conditions for the solution
  \eqref{eq:modulation_solution_kdv} to be well-defined: i)
  $c_g(\omega_0) > 0$ so that $\xi = x/t$ varies across non-negative
  values and ii) $\mathrm{Im}\, k(\omega_0) = 0$ so that the solution
  is real-valued.  These are equivalent to the radiation condition
  \eqref{eq:general_radiation_condition} when
  $0 < \omega_0 < 2/(3\sqrt{3})$.  We proved in Sec.~\ref{third-ss}
  the adherence of the linear KdV equation wavemaker problem to the
  radiation condition \eqref{eq:general_radiation_condition} so a
  unique $k(\omega_0)$ branch exists with these properties
  ($k_0 = k_3$ from eq.~\eqref{eq:KdV-poles-sub}).  With an input
  frequency $\omega_0 \in (0,2/(3\sqrt{3}))$, the solution consists of
  a plane wave with constant frequency $\omega_0$ and wavenumber $k_0$
  for $0<\xi<c_g(\omega_0)$.  The group velocity demarcates the
  invading front of the plane wave into the quiescent region.

  By integrating the generalized wavenumber $\theta_x=k(x/t)$ and
  generalized frequency $\theta_t=-\omega(x/t)$ relations, as well as
  invoking the boundary data $u(0,t) = -\sin(\omega_0 t)$
  \eqref{lkdv-sine-bc}, we can obtain the phase $\theta(x,t)$ up to an
  integration constant $\theta_0 \in \mathbb{R}$ (more generally,
  $\theta_0$ can be slowly varying \cite{haberman_modulated_1988})
  \begin{equation}
    \label{eq:10}
    \theta(x,t) =
    \begin{cases}
      k_0x - \omega_0t, & 0 \le x \le c_g(\omega_0)t \\
      -\frac{2\sqrt{3}}{9} t(1-x)^{3/2} + \theta_0, & c_g(\omega_0)t <
                                                      x \le t .
    \end{cases}
  \end{equation}
  In addition, solving \eqref{eq:Whitham_2} yields the amplitude
  \begin{equation} \label{eq:Whitham_a}
    a(\xi,t) = 
    \begin{cases}
      1, & 0\leq\xi\leq c_g(\omega_0)\\
      \frac{F(\xi)}{\sqrt{t}}, & c_g(\omega_0)<\xi\leq1
    \end{cases}.
  \end{equation}
  The self-similar solution \eqref{eq:Whitham_a} for $a$ does not
  preserve the hydrodynamic scaling symmetry when
  $c_g(\omega_0) < \xi \le 1$ and $F(\xi)$ can be any smooth function.
  This is due to the degeneracy of the double characteristic speed
  $c_g(\omega_0)$.

  These results agree with the large $t$ asymptotic solutions in
  eqs.~\eqref{eq:LKdV_soln_plane} and \eqref{eq:LKdV_xi<1_decay_soln}.
  Note that neither the phase $\theta$ nor the amplitude $a$ are
  continuous at $\xi = c_g(\omega_0)$.  Here, $\theta_0$ and $F(\xi)$
  cannot be determined through leading order modulation
  theory. However, we can appeal to our long time asymptotic result
  \eqref{eq:LKdV_xi<1_decay_soln} in order to identify
  \begin{equation}
    \label{eq:9}
    \theta_0 = 0, \quad F(\xi) = \frac{27 \omega_0 \xi}{(4-27\omega_0^2 -
      \xi^2(\xi+3))\sqrt{2\pi \sqrt{3(1-\xi)}}} .
  \end{equation}
  It is interesting to note that, for $c_g(\omega_0)t < x \le t$, the
  phase $\theta(x,t)$ is independent of the boundary data.  In this
  region, the phase is determined by a structural property of the
  linear KdV equation, namely the saddle points $\rho_{1,2}$, which
  are independent of $\omega_0$.  See the discussion leading up to
  \eqref{eq:LKdV_xi<1_decay_soln}.

  In summary, modulation theory is a simple method for determining
  some features of the long-time asymptotic solution.  However, its
  leading order manifestation does not provide information on the
  integrated phase $\theta_0$, the amplitude function $F(\xi)$, nor
  any information on the supercritical regime $\omega_0>2/(3\sqrt{3})$
  or when $\xi>1$.  Yet, the properties that can be ascertained from
  modulation theory such as the frequency $\omega$, wavenumber $k$,
  and amplitude $a$ agree with the leading order, long-time asymptotic
  solutions \eqref{eq:LKdV_soln_plane} and
  \eqref{eq:LKdV_xi<1_decay_soln} of Regions I and II.a, respectively
  in Fig.~\ref{fig:LKdV_soln_descent}.
\end{remark}

\section{Uniform in $x$ long-time asymptotics for the linear BBM equation with a sinusoidal boundary condition}
\label{bbm-asymp-s}

This section studies the linear BBM wavemaker problem in the
asymptotic limit $t\to\infty$ and $\xi=x/t=\mathcal{O}(1)$.  The
analysis is similar to that for the linear KdV wavemaker problem
studied in the previous section.  Consider the linear BBM
initial-boundary value problem \eqref{lbbm} posed on the semi-infinite
line with initial data $u(x,0)=0$, $x\geq0$ and the time-periodic
sinusoidal boundary condition $u(0,t)=-\sin(\omega_0 t)$, $t\geq 0$,
where $\omega_0>0$. The solution of the problem can be derived via the
Fokas method given by Equation \eqref{lbbm-sol}, where the contour of
integration $\mathcal{C}$ is a small, closed contour around $k = i$ as
shown in Figure \ref{lbbm-f}. We rewrite the solution \eqref{lbbm-sol}
in terms of trigonometric functions
\begin{subequations}     \label{eq:BBM_soln_exact}
    \begin{align}
        u(x,t) &= u_1(x,t) + u_2(x,t) + e^{-x} \sin(-\omega_0  t ), \\
        u_1(x,t) &=  \frac{1}{2\pi} \int_\mathcal{C} \frac{(1-k^2)e^{ikx}}{k^2-(1+k^2)^2\omega_0^2}\left[ -\omega_0 \cos(\omega_0 t )+ i\frac{k}{1+k^2}\sin(\omega_0 t )  \right] \, dk, \label{eq:BBM_soln_exact_1}\\
        u_2(x, t) &= \frac{\omega_0}{2\pi}\int_\mathcal{C} h(k)  e^{ikx-i\frac{k}{1+k^2} t } \, dk, \quad h(k):=\frac{1-k^2}{k^2-(1+k^2)^2\omega_0^2}. \label{eq:BBM_soln_exact_2}
\end{align} 
\end{subequations}
The integrand in equation \eqref{eq:BBM_soln_exact_1} has the simple
pole $k=i$ inside the contour $\mathcal{C}$ so computing the
corresponding residue gives $u_1(x,t)=\sin(\omega_0 t
)e^{-x}$. Therefore, the solution of the linear BBM wavemaker problem
is reduced to evaluating \eqref{eq:BBM_soln_exact_2} whose integrand
possesses an essential singularity at $k=i$. Then, the solution is
simply $u(x,t)=u_2(x,t)$.

We appeal to the method of steepest descent for studying the long-time
asymptotic behavior of the integral \eqref{eq:BBM_soln_exact_2}. The
poles of $h(k)$ occur when
$k^2-(1+k^2)^2\omega_0^2=0$ and are given by $k =
k_j$ where
\begin{equation}
k_1 = \frac{1+\sqrt{1-4\omega_0^2}}{2\omega_0}, \quad  k_2=\frac{1-\sqrt{1-4\omega_0^2}}{2\omega_0}, \quad 
k_3 = \frac{-1+\sqrt{1-4\omega_0^2}}{2\omega_0}, \quad k_4=\frac{-1-\sqrt{1-4\omega_0^2}}{2\omega_0}.
\label{eq:BBM_singularity}
\end{equation}
When $0 < \omega_0 < 1/2$, all the poles are real and this corresponds
to the subcritical regime.  When $\omega_0 > 1/2$, the poles are
complex and corresponds to the supercritical regime.  We therefore
define
\begin{equation}
  \label{eq:bbm_wcrit}
  \omega_0 =\omega_{cr} = 1/2
\end{equation}
as the critical frequency.  Defining the phase
$\phi(k)=ik\xi -i\frac{k}{1+k^2}$, where $\xi=x/ t$, then the saddle
points satisfy $\phi'(k)=i\xi -i\frac{1-k^2}{(1+k^2)^2} =0$ and hence
are located at
\begin{equation}
\begin{aligned}
  &\rho_1 = \sqrt{-1-\frac{1-\sqrt{1+8\xi }}{2\xi}}, \quad \rho_2 =
  -\sqrt{-1-\frac{1-\sqrt{1+8\xi }}{2\xi}},
  \\
  &\rho_3=
  \sqrt{-1-\frac{1+\sqrt{1+8\xi }}{2\xi }}, \quad \rho_4 = -
  \sqrt{-1-\frac{1+\sqrt{1+8\xi }}{2\xi }}. \label{eq:LBBM_saddles} 
  \end{aligned}
\end{equation}
If $0 < \xi < 1$, we have $\rho_1,\rho_2\in\mathbb{R}$ with
$\rho_1\in(0,1)$, $\rho_2\in(-1,0)$, and $\rho_3,\rho_4$ are purely
imaginary with $\text{Im}(\rho_3)>\sqrt{3}$ and
$\text{Im}(\rho_4)<-\sqrt{3}$. If $\xi>1$, all four saddle points are
purely imaginary and symmetrically distributed along the imaginary
axis with $\text{Im}(\rho_1)\in(0,1)$, $\text{Im}(\rho_2)\in(-1,0)$,
$\text{Im}(\rho_3)\in(1,\sqrt{3})$, and
$\text{Im}(\rho_4)\in(-\sqrt{3},-1)$. Let $k=p+iq$ so that
$ \phi(k)=\left(-\xi q + \frac{-p^2 q-q^3+q}{4 p^2
    q^2+\left(p^2-q^2+1\right)^2}\right) + i\left(\xi p + \frac{-p^3-p
    q^2-p}{4 p^2 q^2+\left(p^2-q^2+1\right)^2}\right).  $ The steepest
descent paths occur when
$\text{Im}\{\phi(k)\} = \text{Im}\{\phi(\rho_j)\}$ is a constant and
$\text{Re}\{\phi(k)-\phi(\rho_j)\}<0$, $j=1,2,3,4$.

\begin{figure}[tb!]
    \centering
    \begin{minipage}{.5\textwidth}
        \centering
        \sidesubfloat[]{\includegraphics[width=0.8\linewidth]{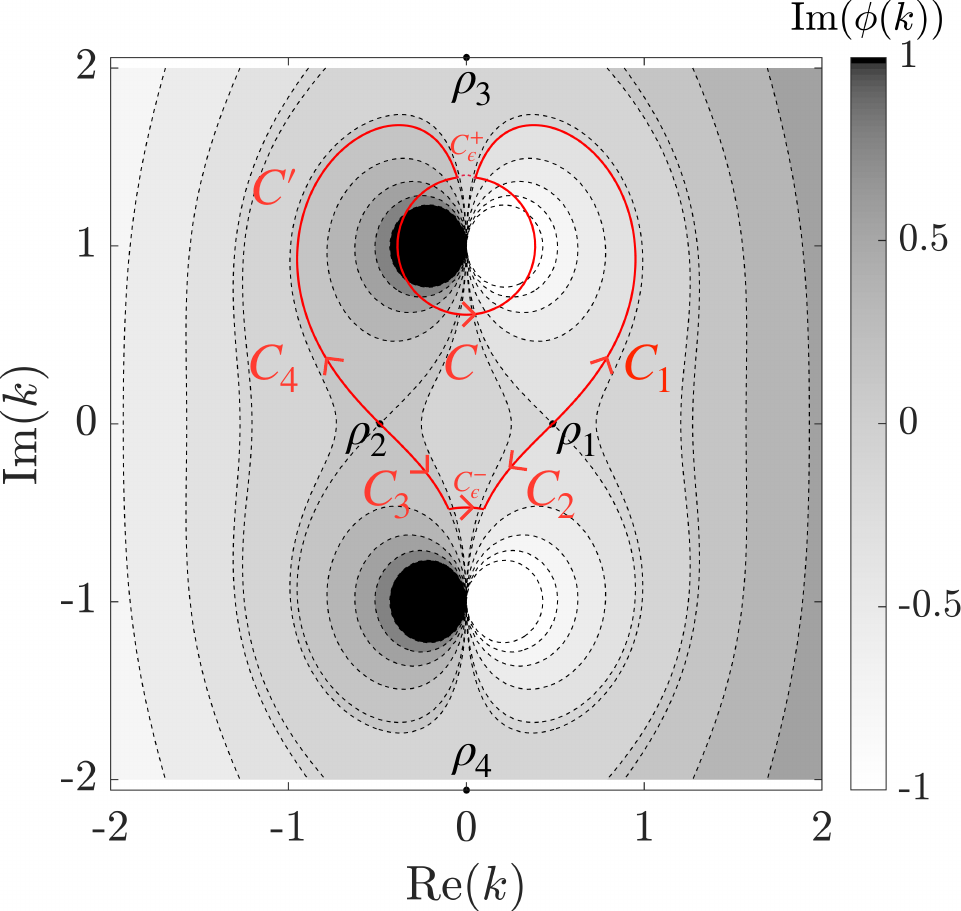} \label{fig:LBBM_contour_sub} } 
    \end{minipage}%
    \begin{minipage}{0.5\textwidth}
        \centering
        \sidesubfloat[]{\includegraphics[width=0.8\linewidth]{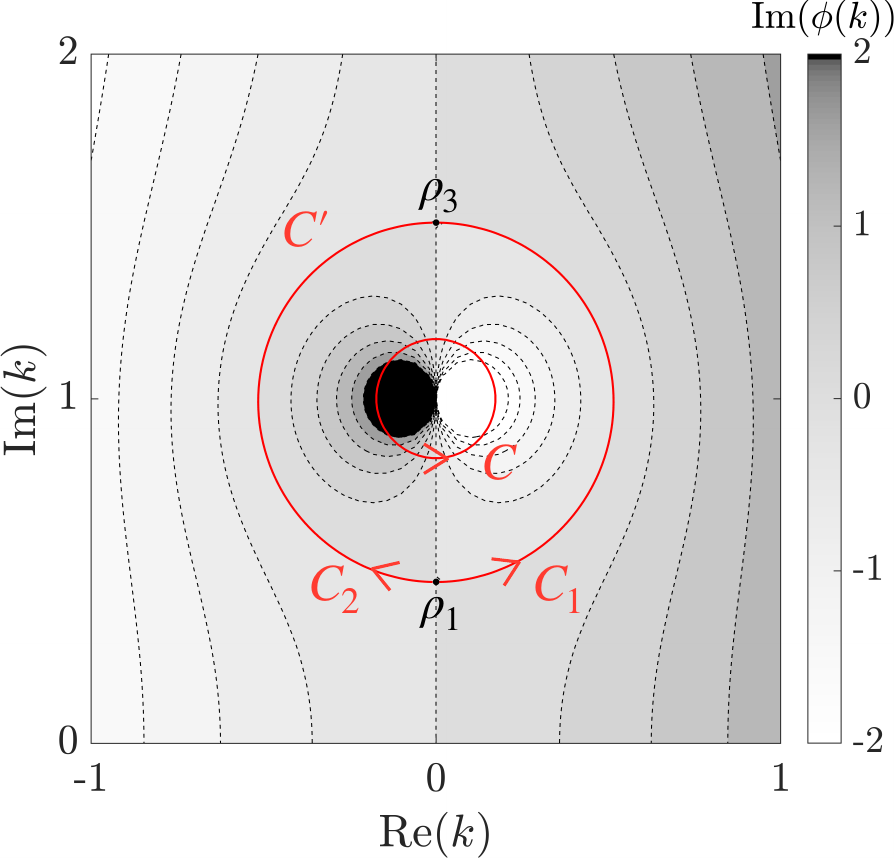} \label{fig:LBBM_contour_super}} 
      \end{minipage}
      \caption{The deformed, steepest descent contour $\mathcal{C}'$
        (red) for the linear BBM wavemaker problem solution
        \eqref{eq:BBM_soln_exact_2} in the case of (a) $0<\xi<1$ and
        (b) $\xi>1$. The positions of the saddles $\rho_j$ are given by
        \eqref{eq:LBBM_saddles}.  The closed curve $\mathcal{C}$ is the original small
        contour around $k=i$. The steepest descent contours
        $\mathcal{C}_1 \cup \mathcal{C}_2$,
        $\mathcal{C}_3 \cup \mathcal{C}_4$ in (a) are given by
        Equations \eqref{eq:BBM_contour_k1},\eqref{eq:BBM_contour_k2},
        and $\mathcal{C}_1 \cup \mathcal{C}_2$ in (b) is
        \eqref{eq:LBBM_contour_xi>1}. The color map indicates levels
        of $\text{Im}(\phi(k))$.}
      \label{fig:LBBM_contours}
\end{figure}

Now we select the steepest descent contours $\mathcal{C}'$ passing
through the saddles that avoid the essential singularities at
$k = \pm i$.  Then, we deform $\mathcal{C}$ onto $\mathcal{C}'$ and
consider the large $t$ behavior.  Our general strategy will be to
first consider the limit $t \to \infty$ and then show that the
remaining contributions from avoiding the essential singularities are
negligible.  Depictions of the contours in
Figure \ref{fig:LBBM_contours} will be discussed after splitting the
range of $\xi$ into different intervals.

\begin{enumerate}[label=(\alph*), wide, labelwidth=!, labelindent=0pt]
\item Regime $0 < \xi < 1$.

  The deformed heart-shaped contour (reminiscent of the {\it Catalpa} tree leaf) is shown in Figure
  \ref{fig:LBBM_contour_sub}, which passes through the saddles
  $\rho_1$ and $\rho_2$.  It is not necessary to deform onto the
  saddles $\rho_3$ and $\rho_4$.  The steepest descent curves passing
  through the saddles $\rho_{1,2}$ are respectively given by
  \begin{subequations}
    \label{eq:BBM_contour}
    \begin{align}
      \text{Im}(\phi(k)) = \text{Im}(\phi(\rho_1)) \implies
      \xi p +  \frac{-p^3-p q^2-p}{4 p^2 q^2+\left(p^2-q^2+1\right)^2}
      &= -\frac{2 \sqrt{2\xi} \left(-2 \xi +\sqrt{8 \xi
        +1}-1\right)^{3/2}}{\left(\sqrt{8 \xi
        +1}-1\right)^2},  \label{eq:BBM_contour_k1}\\ 
      \label{eq:BBM_contour_k2}  
      \text{Im}(\phi(k)) = \text{Im}(\phi(\rho_2)) \implies
      -\xi p +  \frac{-p^3-p q^2-p}{4 p^2
      q^2+\left(p^2-q^2+1\right)^2} &= \frac{2 \sqrt{2\xi} \left(-2
                                      \xi +\sqrt{8 \xi +1}-1
                                      \right)^{3/2}}{\left(\sqrt{8 \xi
                                      +1}-1\right)^2}.
    \end{align}
  \end{subequations}
  The heart-shaped steepest descent paths \eqref{eq:BBM_contour}
  intersect at $k=\pm i$ for all $\xi$; however, the contours we use
  avoid $k=\pm i$ by introducing a small opening
  $\mathcal{C}_\epsilon^+$ of the circle near $k = i$ and a small arc
  $\mathcal{C}_\epsilon^-$ near $k = -i$ as shown in
  Fig.~\ref{fig:LBBM_contour_sub}.  By parameterizing $k$ on
  $\mathcal{C}_\epsilon^+$ as $k=i+\epsilon e^{i\theta}$ for small
  $\epsilon$, we have the restriction
  $\theta \in [\pi/2-\delta(\epsilon),\pi/2+\delta(\epsilon)]$ and
  $0 < \delta(\epsilon) < \pi/2$ is defined by the intersection of $C$
  and $C'$ (see Fig.~\ref{fig:LBBM_contour_sub}).  Then, we obtain
  \begin{align*}
    \left| \int_{\mathcal{C}_\epsilon^+} h(k) e^{t\phi(k)} \, dk
    \right| \le \epsilon
    \int_{\pi/2-\delta(\epsilon)}^{\pi/2+\delta(\epsilon)}  \left|
    h(i+\epsilon e^{i\theta}) \right| \,  e^{t \rm Re(\phi(i+\epsilon
    e^{i\theta}))} \, d\theta.
  \end{align*}
  Because
  $ \rm Re(\phi(i+\epsilon e^{i\theta}))= -
  \frac{\sin\theta}{2\epsilon} - (\frac{1}{4} +\xi) +
  \mathcal{O}(\epsilon) < 0$ when $0< \epsilon \ll 1$ and
  $\left| h(i+\epsilon e^{i\theta}) \right| < M$ is bounded, we obtain
  $\left| \int_{\mathcal{C}_\epsilon^+} h(k) e^{t\phi(k)} \, dk
  \right| \le \epsilon 2 M \delta(\epsilon)$.  Since this bound is
  independent of $t \ge 0$, taking $t \to \infty$ and then
  $\epsilon \to 0$ results in a vanishing contribution.  A similar
  argument implies that the contribution to the integral from the
  small arc $\mathcal{C}_\epsilon^-$ above $k = -i$ vanishes
  $\epsilon \to 0$. Therefore, we can refine the deformed contour
  $\mathcal{C}'$ as
$$
\begin{aligned}
	\int_{\mathcal{C}'}  h(k) e^{t\phi(k)} \, dk  &= \int_{
      \left(\mathcal{C}-\mathcal{C}_\epsilon^+ \right) - \mathcal{C}_1
      +  \mathcal{C}_2 -  \mathcal{C}_3 +  \mathcal{C}_4 -
      \mathcal{C}_\epsilon^- } h(k) e^{t\phi(k)} \, dk  
      \\
      &=
    \int_{\mathcal{C} - \mathcal{C}_1 +  \mathcal{C}_2 -
      \mathcal{C}_3 +  \mathcal{C}_4 } h(k) e^{t\phi(k)} \, dk  +
    \int_{- \mathcal{C}_\epsilon^+ - \mathcal{C}_\epsilon^-} h(k) e^{t\phi(k)} \, dk,
    \end{aligned}
    $$
    so that the solution $u$ in Equation \eqref{eq:BBM_soln_exact} can
    be represented as
    \begin{equation}
      \label{eq:bbm_cauchy_deformation}
        u(\xi t, t ) 
        = \int_{\mathcal{C}} h(k) e^{t \phi(k)}\,dk 
        = \int_{\mathcal{C}'} h(k) e^{t \phi(k)}\,dk +
        \int_{\mathcal{C}_1 - \mathcal{C}_2 +
          \mathcal{C}_3 - \mathcal{C}_4}h(k) e^{t \phi(k)}\,dk + 
          \int_{\mathcal{C}_\epsilon^+ + \mathcal{C}_\epsilon^-}
          h(k) e^{t \phi(k)}\,dk .
    \end{equation}
    We have just proven that the third integral on the right hand side
    vanishes in the limit $\epsilon \to 0$ for any $t > 0$.  We now
    consider the second integral on the right hand side and prove that
    its dominant contributions come from the saddles on
    $\mathcal{C}_1 \cup \mathcal{C}_2$ and
    $\mathcal{C}_3 \cup \mathcal{C}_4$.  The first integral may give
    pole contributions, which we consider afterward.

    Now, we require $0 < \epsilon < 1$.  On the contour
    $\mathcal{C}_1\cup\mathcal{C}_2$ of Figure
    \ref{fig:LBBM_contour_sub}, the neighborhood of $\rho_1$
    dominates. In a neighborhood of $\rho_1$, the path
    $\mathcal{C}_1 \cup \mathcal{C_2}$ is approximately parameterized by
    $k=\rho_1\pm(1+i)s$, $s\in(-\epsilon,\epsilon)$, so that by the
    Laplace method, we obtain
    \begin{equation} \label{eq:BBM_xi<1_saddle1}
      \begin{aligned}
        \int_{\mathcal{C}_1-\mathcal{C}_2} h(k) e^{t\phi(k)} \, dk
        & \sim \frac{e^{i  t \text{Im}(\phi(\rho_1))}}{\sqrt{t}} 
          \frac{1-\rho_1^2}{\rho_1^2(1-2\omega_0^2)-\omega_0^2(1+\rho_1^4)} 
          \sqrt{\frac{\pi}{ |\phi''(\rho_1)|}}, \quad t\to \infty.
      \end{aligned}
    \end{equation}
    Similarly, on the contour $\mathcal{C}_3\cup\mathcal{C}_4$ of
    Figure \ref{fig:LBBM_contour_sub}, we obtain the dominant
    contribution by integrating over a neighborhood of $\rho_2$.  Since
    $\rho_1^2=\rho_2^2$,
    $\rm Im(\phi(\rho_1))= - \rm Im(\phi(\rho_2))$ and
    $|\phi''(\rho_1)|=|\phi''(\rho_2)|$, we sum up the saddle point
    contributions as
    \begin{equation}  \label{eq:BBM_xi<1_saddle2}
      \begin{aligned}
        \int_{\mathcal{C}_1 - \mathcal{C}_2 + \mathcal{C}_3 -
        \mathcal{C}_4} h(k) e^{t\phi(k)}   \,dk
        &\sim
          \frac{\cos\left[\text{Im}(\phi(\rho_1)) t \right]}{\sqrt{t}}
          \frac{2(1-\rho_1^2)}{\rho_1^2(1-2\omega_0^2)-\omega_0^2(1+\rho_1^4)}
          \sqrt{\frac{\pi}{|\phi''(\rho_1)|}}
          , \quad t \to
          \infty,
      \end{aligned}
    \end{equation} 
    which is $\mathcal{O}( t ^{-1/2})$.  To determine how the poles
    $k_j$, $j=1,2,3,4$, in Equation \eqref{eq:BBM_singularity} enter
    the contour integral, we need to consider different cases for
    $\omega_0$.

    \begin{enumerate}[label=(a.\arabic*), wide, labelwidth=!, labelindent=0pt]

      \bigskip
    \item Subcritical regime:  $0<\omega_0<1/2$.

      In this case, all the poles lie on the real line with
      $k_1\in(1,\infty)$, $k_2\in(0,1)$, $k_3\in(-1,0)$ and
      $k_4\in(-\infty,-1)$.  While it is easy to see that
      $k_{1,4}$ are outside of the heart-shaped contour $\mathcal{C}'$,
      $k_{2,3}$ are inside $\mathcal{C}'$ when $k_2<\rho_1$ and
      $k_3>\rho_2$. Observe that
      \begin{equation} 
        \omega'(k_2) = \frac{1-4\omega_0^2+\sqrt{1-4\omega_0^2}}{2}
        , \quad 0<\omega_0<1/2 ~ \text{ and } ~
        0<\xi< \omega'(k_2).  \label{eq:BBM_kappa_pos_1} 
      \end{equation}
      When $\xi < \omega'(k_2)$, $k_{2,3}$ are in $\mathcal{C}'$
      (see Region I in Figure \ref{fig:LBBM_soln_descent}). Applying
      Cauchy's residue theorem and using
      $$
      \text{Res}_h(k_2) 
      = - \frac{1}{2\omega_0} \exp\left[i\left(-\frac{x}{2\omega_0} +
          \frac{x\sqrt{1-4\omega_0^2}}{2\omega_0} + \omega_0 t \right)
      \right] , \quad \text{Res}_h(k_4) = \text{Res}_h(k_2)^*, 
      $$
      we obtain the asymptotic solution to the linear BBM sinusoidal
      initial-boundary value problem
      \begin{equation}
        \begin{split}
          u(\xi t, t ) 
          &= \sin \left[ \left(
            \frac{1-\sqrt{1-4\omega_0^2}}{2\omega_0} \xi - \omega_0
            \right) t   \right] + \mathcal{O}(t^{-1/2}), \quad t \to
            \infty ,
        \end{split}
        \label{eq:LBBM_soln_I}
      \end{equation}
      where the saddle contribution in Equation
      \eqref{eq:BBM_xi<1_saddle2} is the $\mathcal{O}(t^{-1/2})$
      correction.  We notice that $k_0:=k_2=\frac{1-\sqrt{1-4\omega_0^2}}{2\omega_0}$ corresponds to the unique root that is given by the D-N map and satisfies the radiation condition \eqref{eq:general_radiation_condition}.

      When $\xi > \omega'(k_2)$, $k_{2,3}$ are outside
      $\mathcal{C}'$ (see Region II in Figure
      \ref{fig:LBBM_soln_descent}) so that there are no pole
      contributions and the saddle contributions dominate. Using
      eq.~\eqref{eq:BBM_xi<1_saddle2}, we obtain the leading-order
      asymptotic solution as $t \to \infty$
      \begin{small}
        \begin{equation}
          \begin{split}
            u( \xi t, t)
            &\sim \frac{\cos\left[\frac{2
              \sqrt{2\xi}  }{\left(\Xi-1\right)^2}
              \left(-2 \xi +\Xi-1 \right)^{3/2}  t
              \right] }{\sqrt{t}}
              \frac{ \omega _0 \left(-4 \xi +\Xi-1\right)
              \xi^{1/4} \left(\Xi-1\right)^{3/2}
              }{2^{5/4} \sqrt{\pi } \left(\left(4 \xi -\Xi+1\right)
              \omega _0^2+\xi  \left(2 \xi - \Xi + 1\right)\right)
              \left(\Xi^2 -\Xi\right)^{1/2} 
              \left ( -2 \xi -1 +\Xi \right )^{1/2}}, \\
            &:= u_{s1}(\xi t,t) ,
          \end{split}
          \label{eq:LBBM_soln_II}
        \end{equation}
        \end{small}
        where $\Xi := \sqrt{8 \xi + 1}$.  The solution
        \eqref{eq:LBBM_soln_II} is $\mathcal{O}(t^{-1/2})$, indicating
        the typical algebraic decay in time of dispersive waves
        propagating faster than the group velocity.

        \bigskip
      \item Supercritical regime:  $\omega_0>1/2$.

        The poles $k_j$, $j=1,2,3,4$, are all complex and lie on the unit
        circle $\text{Re}(k_j)^2 + \text{Im}(k_j)^2=1$.  Furthermore,
        $k_{2,4}$ are in the lower half plane and are not enclosed by
        $\mathcal{C}'$, while $k_{1,3}$ are inside $\mathcal{C}'$ when
        $\text{Im} \phi(k_1)<\text{Im}\phi(\rho_1)$. Setting
        $\text{Im} \phi(k_1) = \text{Im} \phi(\rho_1)$ defines the
        curve
        \begin{equation}
          l_3(\xi,\omega_0) := \frac{\xi}{2\omega_0}-\omega_0 + \frac{2
            \sqrt{2\xi}  \left( -2 \xi + \Xi(\xi) -1
            \right)^{3/2}}{\left(\Xi(\xi)-1\right)^2} = 0, \quad
          \omega_0>1/2 ~\text{ and }~ 0<\xi<1, \quad \Xi := \sqrt{8 \xi +1},
          \label{eq:BBM_kappa_pos_2}
        \end{equation}
        in the $\omega_0$-$\xi$ plane.  When $l_3(\xi,\omega_0)>0$
        (Region III in Figure \ref{fig:LBBM_soln_descent}), $k_{1,3}$
        are inside the heart-shaped contour $\mathcal{C}'$ and give
        the pole contributions
        $$
        \text{Res}_h(k_1) = \frac{1}{2 \omega_0}
        \exp\left[i\left(\frac{\xi}{2\omega_0} - \omega_0 \right) -
          \frac{\xi t\sqrt{4\omega_0^2-1}}{2\omega_0} \right] , \quad
        \text{Res}_h(k_3) = - \text{Res}_h(k_1)^*.
        $$
        The pole $k_1$ is consistent with the singularity $k_0$ that lies on $\p D^+$ in the D-N map and satisfies the radiation condition \eqref{eq:general_radiation_condition}.
        Combining the pole and saddle contributions, we obtain the
        asymptotic solution
        \begin{equation}
          u( \xi t, t) \sim \exp \left[-\frac{t \xi
              \sqrt{4\omega_0^2-1}}{2\omega_0}\right] \sin\left[\left(
              \frac{\xi}{2\omega_0}-\omega_0 \right) t \right] + u_{s1}(\xi
          t,t), \quad t \to \infty,
          \label{eq:BBM_soln_III}
        \end{equation}
        where $u_{s1}$ is the contribution from the saddles given in
        \eqref{eq:LBBM_soln_II}. Similar to the long-time asymptotic
        solution of the linear KdV wavemaker problem
        \eqref{eq:LKdV_xi<1_super_soln}, both terms in
        \eqref{eq:BBM_soln_III} are required for a uniform asymptotic
        expansion in $x = \xi t$.

        On the other hand, when
        $l_3(\xi,\omega_0)<0$ (Region II in Figure
        \ref{fig:LBBM_soln_descent}), i.e. no poles are inside the
        deformed contour $\mathcal{C}'$, the neighborhoods of the
        saddles $\rho_{1,2}$ dominate and the leading-order solution is
        given by \eqref{eq:LBBM_soln_II}.  The transition from regions I or III to region II in $\xi$ corresponds to the collision of the poles $k_{2,3}$ with saddle points $\rho_{1,2}$ and can be made uniform in $\xi$ \cite{felsen_asymptotic_1994}.
      \end{enumerate}

      \bigskip
      \begin{figure}
        \centering
        \includegraphics[width=0.9\linewidth]{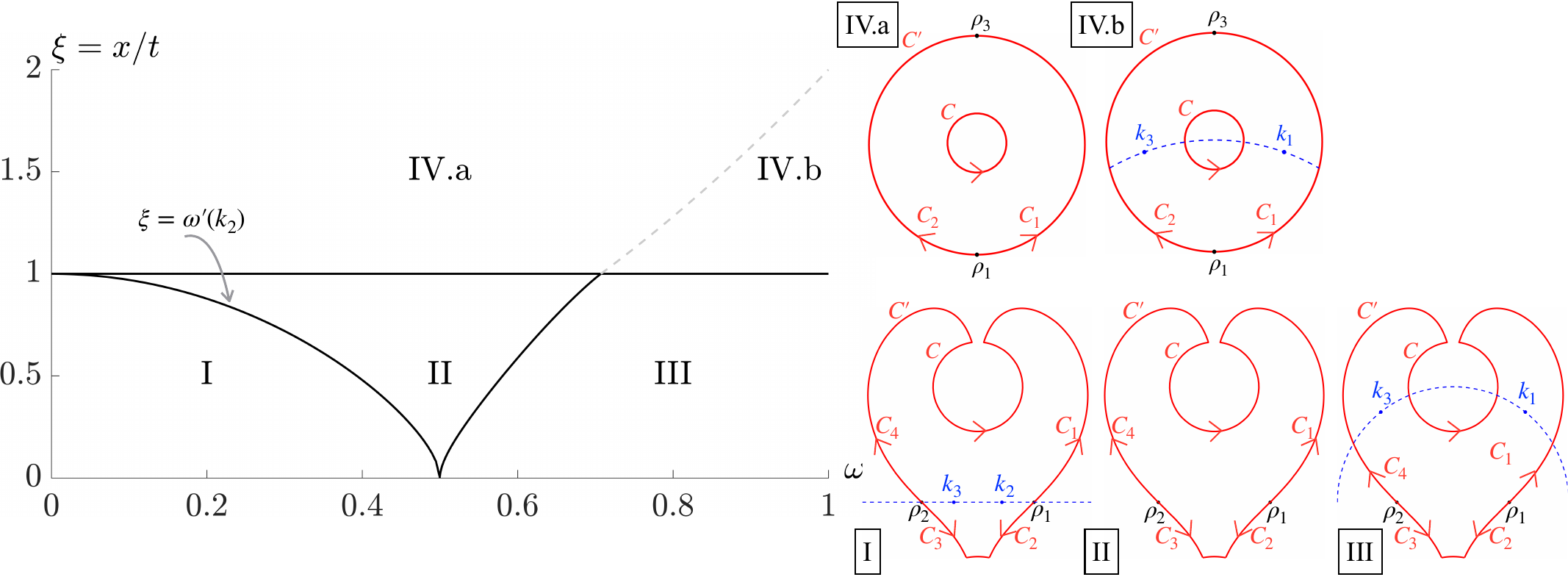}
        \caption{Phase diagram of the long-time asymptotic solutions
          to the linear BBM equation \eqref{lbbm} time-periodic
          sinusoidal boundary value problem. The leading-order
          solutions of each region are given in Section
          \ref{bbm-asymp-s}. The insets present the contour
          deformation and the positions of the poles for different
          values of $\xi$ and $\omega_0$. }
        \label{fig:LBBM_soln_descent}
      \end{figure}
    \item Regime $\xi > 1$.

      Under the condition $\xi>1$, all the saddle points $\rho_j$ are
      imaginary and the second derivative $\phi''(\rho_j)$ at each
      saddle $\rho_j$, $j=1,2,3,4$, is real.  Shown in Figure
      \ref{fig:LBBM_contour_super}, the curve
      $\mathcal{C}' = \mathcal{C}_1 \cup \mathcal{C}_2$ passing
      through $\rho_1$ and $\rho_3$ has constant $\text{Im}(\phi(k))$
      defined by
      \begin{equation} 
        \text{Im}(\phi(k))=\text{Im}(\phi(\rho_1))=0 \implies
        \xi p +  \frac{-p^3-p q^2-p}{4 p^2 q^2+\left(p^2-q^2+1\right)^2} = 0,
        \label{eq:LBBM_contour_xi>1}
      \end{equation}
      where $k = p + i q$.  The steepest descent direction follows the
      circular contour $\mathcal{C}_1 \cup \mathcal{C_2}$, and
      therefore, the neighborhood of the saddle point $\rho_1$ gives
      the dominant contribution for the integral. With
      $\phi(\rho_1),\phi''(\rho_1)\in \mathbb{R}$ and
      $\phi(\rho_1)<0,\phi''(\rho_1)<0$, we apply the Laplace method and
      obtain
    \begin{align}
      \int_{\mathcal{C}_1-\mathcal{C}_2}  h(k)  e^{t\phi(k)} \,dk
      &\sim
        \frac{\exp\left[ t \phi(\rho_1)\right]}{\sqrt{t}} 
        \frac{(1-\rho_1^2)}{\rho_1^2(1-2\omega_0^2)-\omega_0^2(1+\rho_1^4)} 
        \sqrt{\frac{2\pi}{-
        \phi''(\rho_1)}}, \quad t \to \infty.
        \label{eq:BBM_xi>1_saddle}
    \end{align}
    Now we deform the original small closed contour $\mathcal{C}$ onto
    the steepest descent curve $\mathcal{C}'$; in doing so, we
    incorporate the pole contributions from $k_j$ in
    \eqref{eq:BBM_singularity} when necessary.

    \begin{enumerate}[label=(b.\arabic*), wide, labelwidth=!, labelindent=0pt]
      \bigskip
    \item Subcritical regime: $0<\omega_0<1/2$.

      All the poles \eqref{eq:BBM_singularity} are on the real line
      and so are not enclosed in $\mathcal{C}'$ (see Region IV.a in
      Figure \ref{fig:LBBM_soln_descent}). This indicates that
      $\int_{\mathcal{C}} h(k) e^{t\phi(k)} \, dk =
      \int_{\mathcal{C}'} h(k) e^{t\phi(k)} \, dk$.  We use the saddle
      contribution from \eqref{eq:BBM_xi>1_saddle} to obtain the
      leading-order solution
      \begin{equation}
        \begin{aligned}
          u( \xi t , t)
          &\sim  -\frac{e^{-\alpha(\xi) t}}{\sqrt{t}} \frac{ \omega _0
            \left(\frac{ \xi^{1/2} 
            \left(\Xi-1\right)^3  }{\sqrt{2 \xi -\Xi+1} \left(\Xi^2
            -\Xi \right)   }
            \right)^{1/2} \left(4 \xi +\Xi+8 \omega
            _0^2-1\right) }{4\ 2^{3/4}
            \sqrt{\pi } \left((4 \xi -1) \omega _0^2+(\xi -1) \xi +4
            \omega_0^4\right)}, \quad t \to \infty,  \\
          & := u_{s2}(\xi t, t)
        \end{aligned}
        \label{eq:BBM_xi>1_soln_1}
      \end{equation}
      where
      \begin{equation}
        \label{eq:alpha}
        \alpha(\xi) = \frac{ \sqrt{2 \xi^2 - \xi \Xi + \xi}
          \left(\Xi-3\right)  }{ \sqrt{2} 
          \left(\Xi-1\right)}, \quad \Xi := \sqrt{8\xi + 1} .
      \end{equation}
      The solution in this region has rapid decay in time of order
      $\mathcal{O}( t ^{-1/2}e^{- \alpha(\xi) t })$.

      \bigskip
    \item Supercritical regime:  $\omega_0>1/2$.

      The poles $k_{1,3}$ are in the first and second quadrants on a
      unit circle centered at the origin, so are enclosed in
      $\mathcal{C'}$ when
      $\text{Im}(\phi(k_1))<\text{Im}(\phi(\rho_1))=0$ and,
      equivalently,
      $\text{Im}(\phi(k_3))>\text{Im}(\phi(\rho_1))=0$. This implies,
      $k_{1,3}$ are inside $\mathcal{C}'$ when
      $ \omega_0>\frac{1}{\sqrt{2}}$, $1<\xi<2\omega_0^2$ (see Region
      IV.b in Figure \ref{fig:LBBM_soln_descent}) and outside
      $\mathcal{C}'$ when
      $ \frac{1}{2}<\omega_0\leq \frac{1}{\sqrt{2}}$, $\xi>1$ or
      $\omega_0>\frac{1}{\sqrt{2}}$, $\xi>2\omega_0^2$ (Region IV.a in
      Figure \ref{fig:LBBM_soln_descent}).  The poles $k_{2,4}$ are in
      the third and fourth quadrants so are outside of the deformed
      contour $\mathcal{C}'$.
      
      For $k_{1,3}$ outside the contour in Region IV.a, the asymptotic
      solution is given by \eqref{eq:BBM_xi>1_soln_1}.  In Region
      IV.b, the contour integral over $\mathcal{C}$ can be expressed
      as
      \begin{equation}
        \label{eq:2}
        \int_{\mathcal{C}} h(k) e^{t \phi(k)}\,dk =
        \int_{\mathcal{C}'} h(k) e^{t \phi(k)}\,dk - 2\pi i
        (\text{Res}_h(k_1) + \text{Res}_h(k_3)) ,
      \end{equation}
      as shown by the deformation of $\mathcal{C}$ to $\mathcal{C}'$
      in Figure \ref{fig:LBBM_soln_descent}.  Evaluating the residues,
      we obtain
      $$
      \text{Res}_h(k_1) = \frac{1}{2\omega_0}
      \exp\left[i\left(\frac{\xi}{2\omega_0} - \omega_0 \right)t -
        \xi t \frac{\sqrt{4\omega_0^2-1}}{2\omega_0} \right], \quad
      \text{Res}_h(k_3) = - \text{Res}_h(k_1)^* .
      $$ 
      Then, the large $t$ asymptotics are
      \begin{equation}
        \begin{aligned}
          u( \xi t, t) &\sim \exp\left[-\frac{
                         \xi t\sqrt{4\omega_0^2-1}}{2\omega_0}\right]
                         \sin\left[
                         \left(\frac{\xi}{2\omega_0}-\omega_0 \right)
                         t \right] +  u_{s2}, \quad t \to
                         \infty, \label{eq:BBM_xi>1_soln_2}
        \end{aligned}
      \end{equation}
      where $u_{s2}$ is the saddle contribution given in
      \eqref{eq:BBM_xi>1_soln_1}. Because the exponential decay rate
      $|\phi(\rho_1)|$ of $u_{s2}$ is monotone increasing, in Region
      IV.b where $1<\xi<2\omega_0^2$, we have
      $|\phi(\rho_1)|<|\phi(\rho_1) |_{\xi=2\omega_0^2} < \xi
      \sqrt{4\omega_0^2-1} / (2\omega_0)$. Then, the solution
      \eqref{eq:BBM_xi>1_soln_2} is dominated by the saddle
      contribution $u_{s2}$.  Therefore, for large $t$, the leading
      order contribution for all $\omega_0$ and $\xi > 1$ is the same
      and is given by $u_{s2}$.
    \end{enumerate}
  \end{enumerate}

  The large $t$ solution given in eqs.~\eqref{eq:LBBM_soln_I},
  \eqref{eq:LBBM_soln_II}, \eqref{eq:BBM_soln_III},
  \eqref{eq:BBM_xi>1_soln_1}, and \eqref{eq:BBM_xi>1_soln_2} proves
  Theorem \ref{lbbm-t}.
  
  \begin{remark}
    At a qualitative level, the large $t$ solutions of the linear BBM
    and KdV wavemaker problems exhibit similar bifurcation diagrams in
    Figures \ref{fig:LKdV_soln_descent} and
    \ref{fig:LBBM_soln_descent}.  Consequently, the solution behavior
    is qualitatively similar, with traveling periodic waves in Region
    I, algebraically decaying waves in Region II, and exponential
    decay in Regions III and IV of the $\omega_0$-$\xi$ plane.
  \end{remark}

  \begin{remark}[Linear modulation theory] 
    The modulation equations \eqref{eq:Whitham_modulation} are
    applicable for general linear dispersive PDEs.  Following a path
    similar to that of the previous section, we now solve
    \eqref{eq:Whitham_modulation} with the associated boundary values
    \eqref{eq:Whitham_BC} where, in this case, $\omega=k/(1+k^2)$. For
    $\xi = x/t$, the self-similar solution is
    \begin{equation}
      \omega(\xi)= 
      \begin{cases}
        \omega_0, & 0\leq\xi\leq c_g(\omega_0)\\
        \frac{\sqrt{1-4\xi+\sqrt{1+8\xi}}}{2\sqrt{2}}, & c_g(\omega_0)<\xi\leq1
      \end{cases} , \quad 
      k(\xi)= 
      \begin{cases}
        k_0, & 0\leq\xi\leq c_g(\omega_0)\\
        \sqrt{\frac{-1-2\xi+\sqrt{1+8\xi}}{2\xi}}, & c_g(\omega_0)<\xi\leq1
      \end{cases},
    \end{equation}
    where $k_0 = k(\omega_0) > 0$ satisfies the linear dispersion
    relation $\omega_0=k_0/(1+k_0^2)$.  There are two choices
    for $k_0$.  Admissibility of the simple wave solution in the
    subcritical regime $0 < \omega_0 < 1/2$ implies the radiation
    condition $c_g(\omega_0) > 0$.  Then,
    $k(\omega_0)$ is uniquely determined ($k_0 = k_2$ in
    eq.~\eqref{eq:BBM_singularity}) and
    \begin{equation}
      \label{eq:BBM_grp_vel_wavemaker}
      c_g(\omega_0) = \frac{1}{2} \left( 1-4\omega_0^2 +
        \sqrt{1-4\omega_0^2} \right).
    \end{equation}
    Suitable integration of
    $\theta_x = k(x/t)$, $\theta_t = -\omega(x/t)$ determines the
    phase
    \begin{equation}
      \label{eq:phase_BBM_modulation}
      \theta(x,t) =
      \begin{cases}
        k_0 x - \omega_0 t, \quad 0 \le x \le c_g(\omega_0)t, \\
        \frac{2 \sqrt{2} x \left ( \sqrt{t(t+8x)}-t-2x \right
        )^{3/2}}{t\left ( \sqrt{1+8x/t}-1 \right )^2} + \theta_0, \quad
        c_g(\omega_0) t < x < t ,
      \end{cases}
    \end{equation}
    where $\theta_0 \in \mathbb{R}$.  The amplitude $a$ has the
    solution \eqref{eq:Whitham_a} with $c_g(\omega_0)$ given in
    \eqref{eq:BBM_grp_vel_wavemaker} and $F(\xi)$ an arbitrary, smooth
    function.  As before, these results agree with the asymptotic
    solutions \eqref{eq:LBBM_soln_I} and \eqref{eq:LBBM_soln_II} if we
    take
    \begin{equation}
      \label{eq:11}
      \theta_0 = 0, \quad F(\xi) = \frac{ \omega _0 \left(-4 \xi +\Xi-1\right)
        \xi^{1/4} \left(\Xi-1\right)^{3/2}
      }{2^{5/4} \sqrt{\pi } \left(\left(4 \xi -\Xi+1\right)
          \omega _0^2+\xi  \left(2 \xi - \Xi + 1\right)\right)
        \left(\Xi^2 -\Xi\right)^{1/2} 
        \left ( -2 \xi -1 +\Xi \right )^{1/2}},
    \end{equation}
    where $\Xi = \sqrt{8\xi + 1}$.
\end{remark}

\section{Comparison with conduit experiments}
\label{exp-s}

In the long-wavelength, small amplitude regime, the linear BBM
equation (equivalent to the linearized conduit equation
\cite{mao2023experimental}) describes the interfacial wave dynamics of
a particular type of viscous two-fluid core-annular flow
\cite{olson1986solitary, lowman2013dispersive} that models magma 
migration through the upper mantle \cite{barcilon_nonlinear_1986} 
and channelized water flow through glaciers \cite{stubblefield_solitary_2020}.
This flow system
comprises a pressure-driven, buoyant, Stokes core fluid within an
annulus of another, miscible, heavier Stokes fluid, with a small
core-to-annular fluid viscosity ratio. At the injection site, a piston
pump acts as a wavemaker by injecting a steady stream of core fluid
with an added, small amplitude, time-periodic injection to generate
interfacial waves (see \cite{mao2023experimental}). The system can be
realized in a simple laboratory setting and provides a quantitative
experimental platform for the study of linear and nonlinear wave
propagation \cite{olson1986solitary,scott1986observations,helfrich1990solitary,lowman_interactions_2014,maiden2016observation,maiden_solitonic_2018,anderson2019controlling,maiden2020solitary, mao2023experimental,mao2023creation}.  
The linear BBM equation was found 
to more accurately capture the dispersive wave properties
of this viscous core-annular flow across a wider range of frequencies
than the asymptotically equivalent long wave, linear KdV equation \cite{mao2023experimental}.
In this section, we investigate the linear
time-periodic wavemaker problem in this viscous core-annular system
and compare our experimental observations with the long-time
asymptotic solutions of the linear BBM equation, derived in section
\ref{bbm-asymp-s}.

Consider the linearized BBM equation initial-boundary value problem
\begin{subequations}
\begin{align}
  u_t+u_x - u_{xxt} &=0, \quad x > 0, \quad
                      t > 0, \label{eq:conduit-0} \\
  u(x,0) &= 0, \quad x > 0, \\
  u(0,t) &= \sin(-\omega_0 t), \quad t > 0.  \label{eq:conduit_BC}
\end{align}
\end{subequations}
The dependent variable $u$ represents small amplitude fluctuations of
the two-fluid interfacial cross-sectional area scaled by the steady,
constant cross-sectional area of the established, vertical, straight
conduit; $x$ is proportional to the vertical distance from the
boundary, and $t$ is proportional to time.  See
\cite{mao2023experimental} for the non-dimensionalization.  Here, we
will use the dimensional length $L$ and time $T$ scalings as fitting
parameters, which amounts to fitting the fluid viscosities and density 
difference, in order to compare with experiment.  Past work has 
demonstrated quantitative agreement between the nonlinear extension of 
the linearized BBM equation (the conduit equation)
and experiment with no fitting parameters 
\cite{lowman_interactions_2014,maiden_solitonic_2018,maiden2020solitary}. A linear
time-periodic sinusoidal time-series of unit amplitude and frequency
$\omega_0$ is generated at the boundary and propagates into the
constant background.

The experimental configuration is identical to that employed in
\cite{mao2023experimental}. The experiments are conducted within a
vertical square acrylic column that has the cross-section
$5~\text{cm} ~\times~ 5~\text{cm}$ and a height of $180$ cm. The core
fluid is composed of dyed, diluted glycerin which exhibits lower
density and significantly less viscosity than the reservoir fluid into
which it is injected. The core fluid rises buoyantly in
an external pure glycerin reservoir to establish a circularly
symmetric free interface between the two fluids.  A
computer-programmed flow rate pump is used to regulate the flow rate
$\widetilde{Q}$ (cm$^3$/s) of the core fluid, with the cross-sectional
area $\tilde{A}$ (cm$^2$) at the injection site satisfying the
Hagen-Poiseulle pipe flow law 
\begin{equation}
    \label{eq:poiseulle}
    \widetilde{Q} \propto \tilde{A}^2,
\end{equation} 
with a constant of proportionality involving the ratio of the fluid density difference to the interior viscosity. The
imposition of the time-periodic boundary data in \eqref{eq:conduit_BC}
is realized by the injection flow rate
\begin{equation}
  \label{eq:injection_time_series}
  \widetilde{Q}(\tilde{t}) =
  \widetilde{Q}_0(1+2a\sin(-\widetilde{\omega}_0 \tilde{t})  ), \quad
  a \ll1, \quad  \tilde{t}\ge0,
\end{equation}
where $\widetilde{Q}_0$ (cm$^3$/s) is the background flow rate, $a$ is
the (nondimensional) amplitude, $\widetilde{\omega}_0$ (rad/s) denotes
the dimensional angular input frequency and $\tilde{t}$ (s) is time. A high-resolution camera with a 3 Hz sampling rate is used for data acquisition above the injection boundary to capture the evolving wavetrains. The conduit edges are extracted from camera images  using edge detection and smoothing algorithms.  At each vertical image pixel, the conduit diameter is calculated as the number of horizontal pixels between the two conduit edges. The diameter is then converted to the dimensionless conduit cross-sectional area $A$ by dividing by the number of horizontal pixels in the undisturbed conduit prior to commencing time-periodic injection and squaring the result.

For experimental comparison with theoretical predictions, we
appeal to the nondimensionalization parameters $L$ (cm) as the
vertical length scale, $U$ (cm/s) as the vertical velocity scale and
$T=L/U$ (s) as the time scale. For each experiment presented below,
approximately 15 sets of linear periodic traveling waves were first
generated for calibration, and the dimensional angular frequency and wavenumber
($\Tilde{\omega},\Tilde{k}$) were measured. The
scales $L$ and $U$ were determined by fitting
the dimensional quantities ($\Tilde{\omega},\Tilde{k}$) to the BBM
linear dispersion relation
\begin{equation}
  \label{eq:dimensional_bbm_dispersion}
  f(\Tilde{k}) = \frac{U\Tilde{k}}{1+L^2\Tilde{k}^2}.
\end{equation}
The scales associated with each displayed data set are reported. We scale the dimensional space and time as
$x=\tilde{x}/L$, $t=\tilde{t}/T$, respectively. Upon linearization of
the cross-sectional area $A$ on the uniform background
\begin{equation}
    A(x,t) = 1 + a u(x,t), \quad a\ll1, \label{eq:conduit-to-BBM}
\end{equation}
and using the boundary condition \eqref{eq:injection_time_series} and the Hagen-Poiseulle law \eqref{eq:poiseulle} as 
\begin{equation}
    \begin{split}
        A(0,t) &= 1 + a u(0,t) \\
        &= \sqrt{\frac{\tilde{Q}(T t)}{\tilde{Q}_0}} = \sqrt{1 + 2 a \sin(-\tilde{\omega}_0 T t)} \\
        &\sim 1 + a \sin(-\omega_0 t), \quad a \to 0,
    \end{split} 
\end{equation}
where the dimensionless input frequency is $\omega_0 = \tilde{\omega}_0 T$,
we obtain the initial, boundary value problem \eqref{eq:conduit-0}, which
we compare with experiment.

\begin{figure}[tb!]
    \centering
     \begin{minipage}{.45\textwidth}
        \centering
        \sidesubfloat[]{\includegraphics[width=\textwidth]{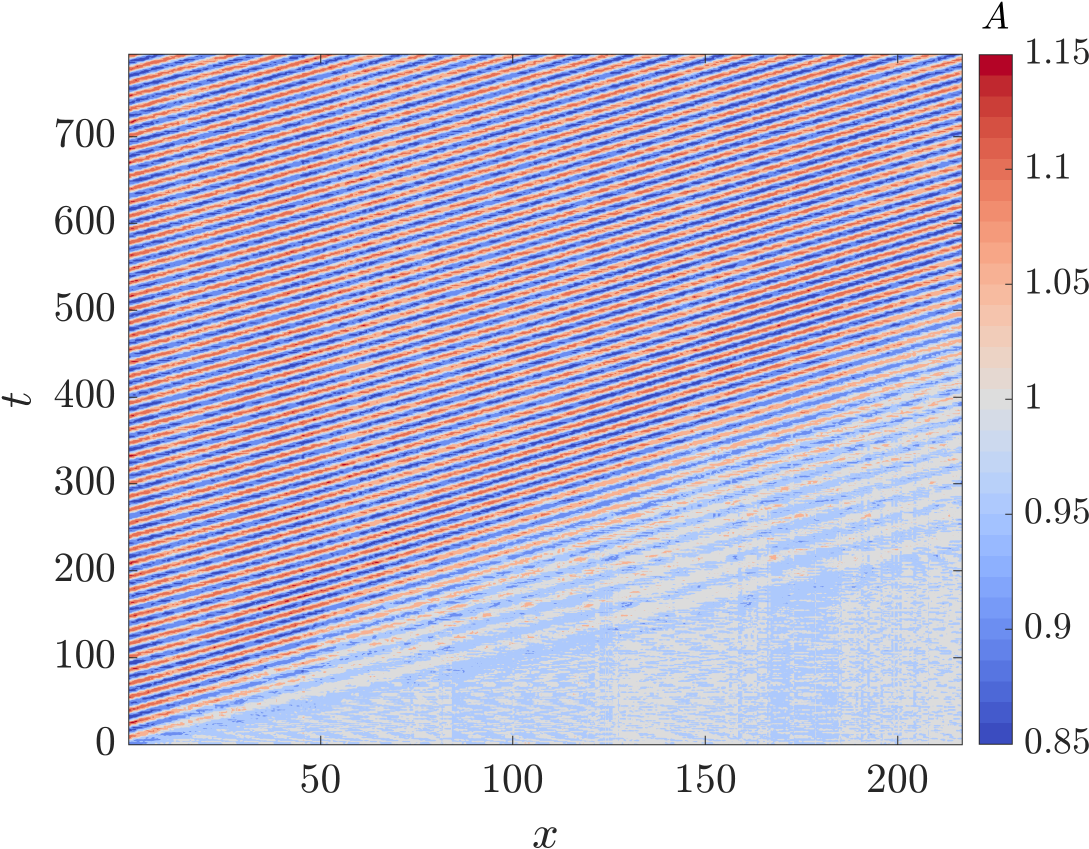} \label{fig:exp_contour_sub}}
    \end{minipage}
    \hspace{0.05\textwidth}  
    \begin{minipage}{.45\textwidth}
        \centering
        \sidesubfloat[]{\includegraphics[width=0.95\textwidth]{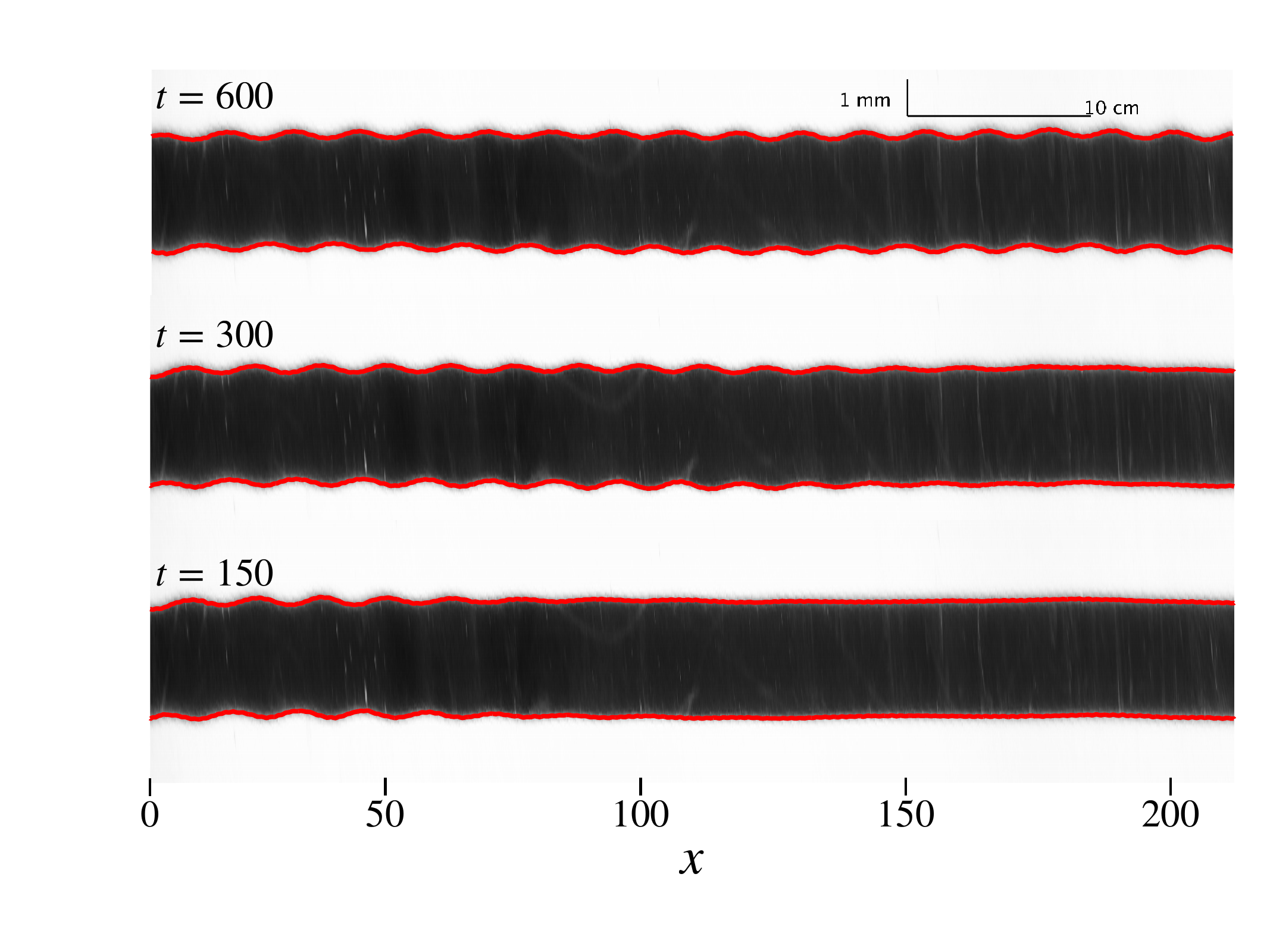} \label{fig:exp_pic_sub} }
    \end{minipage}
    \\
    \begin{minipage}{.45\textwidth}
        \centering
        \sidesubfloat[]{\includegraphics[width=\textwidth]{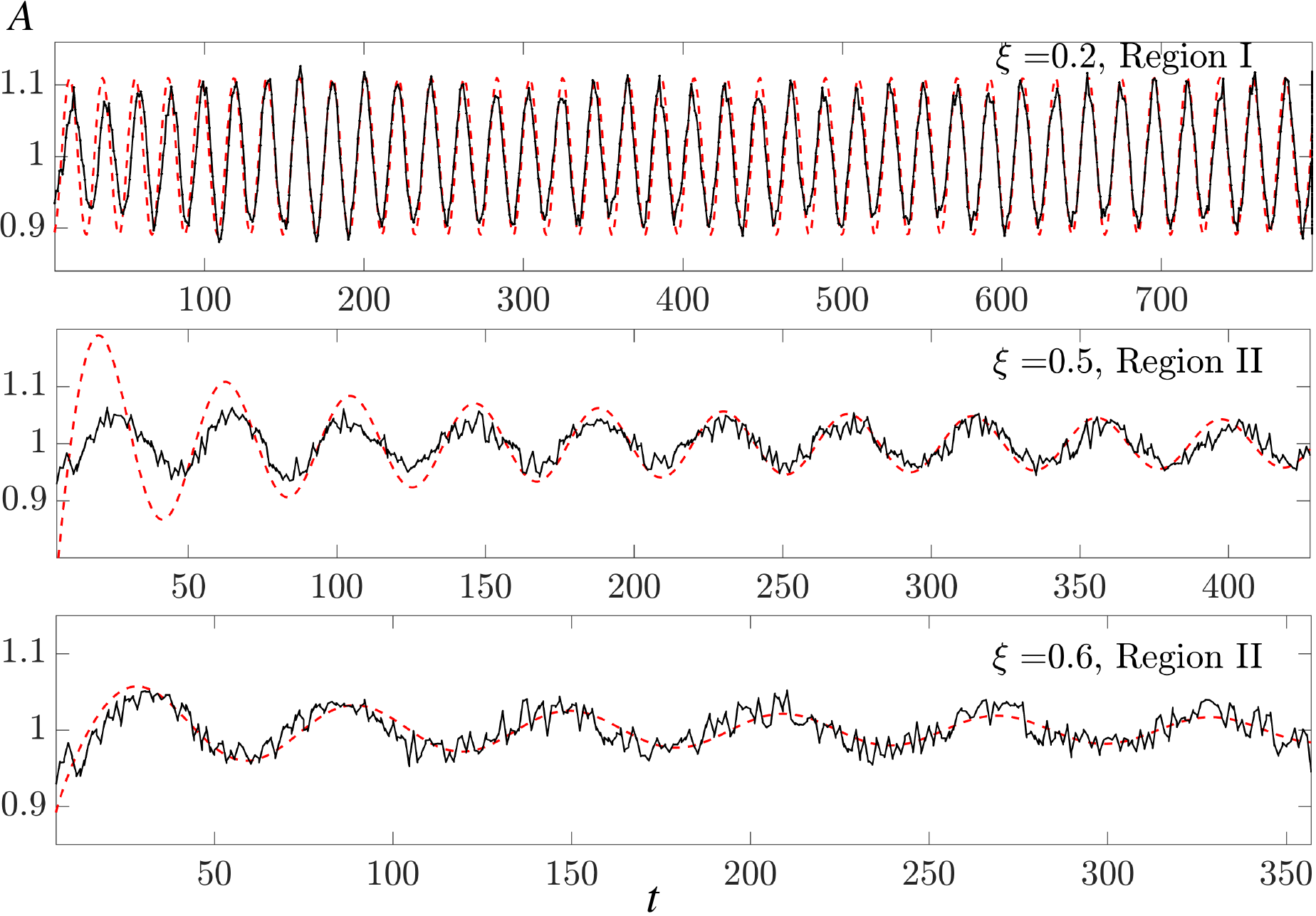} \label{fig:exp_compare_sub} }
    \end{minipage}
    \caption{Viscous core-annular flow experiment in the subcritical regime.  (a) Space-time contour of normalized cross-sectional area $A$ for time-periodic injection of core fluid into a uniform background.  The
      measured periodic wave amplitude and angular frequency are
      $(a,\omega)=(0.11 \pm 0.02,0.41 \pm 0.02)$, with the
      nondimensionalization scales
      $(L,U)=(0.28\pm 0.01\, \text{cm}, 0.96\pm 0.02\,
      \text{cm/s})$. (b) $90^\circ$ clockwise rotated, rescaled
      experimental time-lapse. The red curves are obtained from edge detection. (c) Comparison of the experimental data from (a)
      (black solid) with the long-time asymptotic solutions in Regions I 
      \eqref{eq:LBBM_soln_I} and II \eqref{eq:LBBM_soln_II} (red
      dashed) at different, fixed $\xi$.}
      \label{fig:exp_sub}
\end{figure}

\begin{figure}[tb!]
  \centering
  \includegraphics[width=0.45\linewidth]{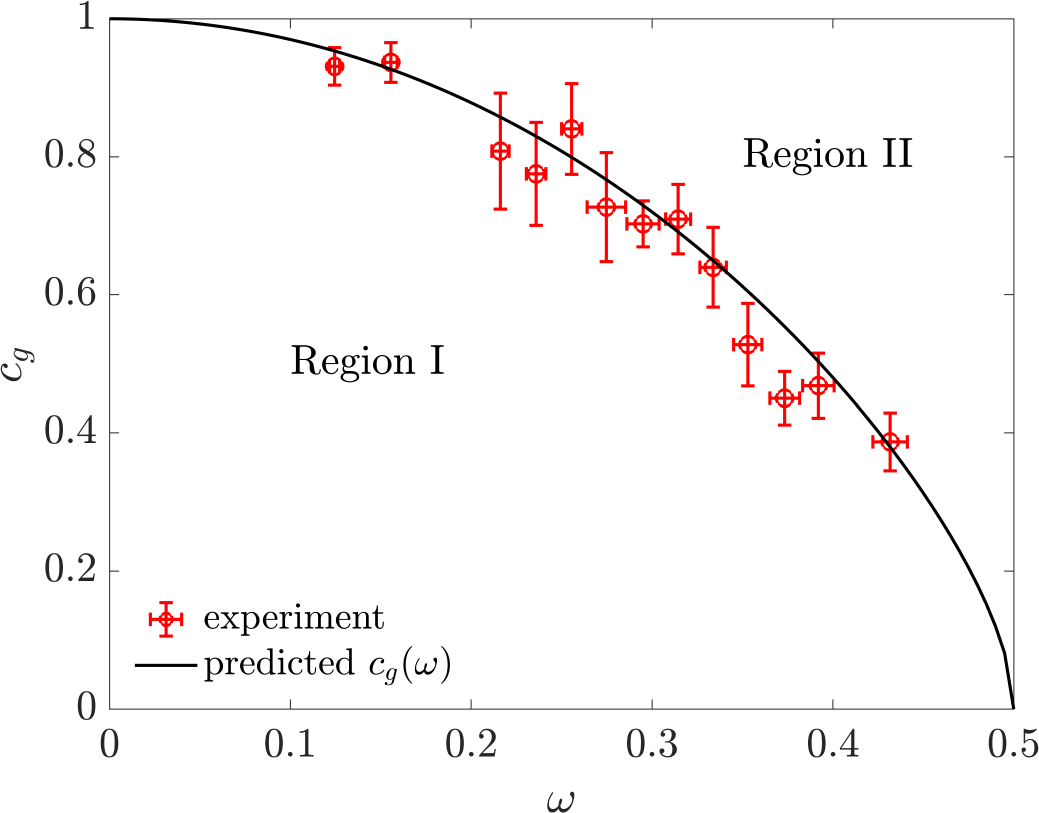}
  \caption{Comparison of the experimentally measured threshold between
    Regions I (periodic traveling waves) and II (algebraically
    decaying waves) at varied angular frequencies $\omega$ (red dots)
    in the subcritical regime with the group velocity
    \eqref{eq:BBM_grp_vel_wavemaker} (black curve), predicted by the
    long-time asymptotic analysis of the linear BBM wavemaker
    problem. The nondimensionalization scales for all thirteen
    measurements are
    $(L,U)=(0.31\pm0.01\,\text{cm}, 1.14\pm 0.02\,\text{cm/s})$. }
     \label{fig:exp_vg}
\end{figure}

Previous viscous core-annular flow experiments examined the accuracy of the linear dispersion relation \eqref{eq:dimensional_bbm_dispersion}, the existence of the critical frequency $\omega_{cr}$ separating sub and supercritical regimes, and the sinusoidal spatial profiles of linear waves, in addition to their nonlinear counterparts \cite{mao2023experimental}.  The BBM dispersion relation \eqref{eq:dimensional_bbm_dispersion} was shown to be an accurate model of observed linear waves, with discrepancies due to a weak recirculating flow in the outer, annular fluid.  In what follows, we compare new experiments with asymptotic predictions of i) the spatial structure of waves along rays $x/t = $ constant in the sub and supercritical regimes and ii) the subcritical group velocity transition between Regions I and II in Fig.~\ref{fig:LBBM_soln_descent}.

The space-time contour plot in Figure \ref{fig:exp_contour_sub} depicts the dimensionless cross-sectional area $A$ for the wavemaker problem in the subcritical regime $0 < \omega_0<\omega_{cr}=1/2$. The color scale denotes wave amplitude with the red-blue region
representing waves of substantial amplitude, whereas the lighter-colored,
blue region corresponds to the small amplitude, noisy background.
The method for measuring the parameters of the periodic traveling wave is described in \cite{mao2023experimental}. Measurement of this experimental dataset yields the periodic wave amplitude, angular frequency and wavenumber $(a,\omega,k)=(0.11\pm0.02, 0.41\pm0.02, 0.52\pm0.02)$, respectively.  
Figure \ref{fig:exp_pic_sub} displays example time-lapse images captured by a high-resolution camera during the experiment. At $t=150$ and $t=300$, the interfacial wave is observed as it enters the uniform conduit, resulting in decay at the leading edge. Subsequently, as the wave propagates, for example at $t=600$, it attains local periodicity with small amplitude relative to the average cross-sectional diameter or area. 
In Figure \ref{fig:exp_compare_sub}, we compare the experimentally
extracted traveling wave data from Fig.~\ref{fig:exp_contour_sub} along rays $x/t = $ constant with the
long-time asymptotic solutions of the linear BBM wavemaker
problem. The comparison utilizes the transformation
\eqref{eq:conduit-to-BBM} and the solutions (\ref{eq:LBBM_soln_I},
\ref{eq:LBBM_soln_II}) in Regions I and II of Figure
\ref{fig:LBBM_soln_descent}, with $a=0.11$,
$\omega_0=0.41$, and the previously fitted scaling parameters $L$ and $T$. However, the experiment identifies $x=0$ with the
starting point of the cropped image and $t=0$ with the initiation of
imaging, which differ slightly from the location of the injection site
and the initiation of the periodic time-series for the injection rate.
Consequently, a shift is applied to $x$ and $t$ everywhere they appear
in the solution formulas.  For example,
\begin{equation}
  \xi = \frac{x-x_0}{t-t_0} .
  \label{eq:exp_Xi}
\end{equation}
The experiment in Figure \ref{fig:exp_compare_sub} is subject to the
shifts $x_0=6$ $(\approx 2\,\text{cm})$ and $t_0=6$ $(\approx 2\, \text{s})$, which are selected by comparing the experimental data with the theoretical solutions at various $\xi$.  
The shift in $x$ is consistent with the experimental setup, in which case the injection site is, while not shown, on the left of the experimental images in Figure \ref{fig:exp_pic_sub}. The shift $t_0$ is the best fit time shift. 
In Figure \ref{fig:exp_compare_sub}, good
agreement is achieved for the comparison between the experiment and
the long-time asymptotic solutions in Regions I ($0 < \xi < c_g$) and II ($c_g < \xi < 1$). Region I
contains the solutions asymptoting to periodic traveling waves, and
the solutions in Region II asymptote to algebraically time-decaying
waves.  According to the linear BBM theoretical prediction, the
corresponding threshold between Regions I and II, i.e. the group
velocity, at $\omega_0=0.41$ is predicted to be
$c_g = \omega'(k_2) = 0.45$ (cf.~Equation
\eqref{eq:BBM_kappa_pos_1}).  For small $t$, the asymptotic solution
deviates somewhat from experiment.  However, this is expected since
the formula represents the large $t$ asymptotics of the solution.  In
particular, the observed oscillation amplitude, period and phase for
sufficiently large $t$ are very close to experiment.  Note that the
small, noisy oscillations correspond to the resolution of our imaging
system.  In particular, we
observe diameter fluctuations on the order of four pixels, which, for
a 30 pixel diameter established conduit corresponds to
$(4/30)^2 \approx 0.02$ fluctuations in the nondimensional
cross-sectional area $A$.  The reason for this resolution is the 200:1
vertical to horizontal aspect ratio of the imaging system used to capture the data in
Figure \ref{fig:exp_contour_sub}.

As shown in Figure \ref{fig:exp_sub}, wave propagation is periodic when $\xi$ is less than a
threshold $\xi_g$ and decays when $\xi > \xi_g$.  To obtain a value for this threshold from experiment, we employ the
following method. Time slices are extracted from the dataset for
$t \in [30,150]$. In each time slice, the rightmost
location $x_g(t)$ is defined to be the last spatial peak with a value greater than
$80\%$ of the averaged periodic wave amplitude at the leading edge
where the wave starts to decay.  The fitted slope of $x_g(t)$, $\xi_g$, is a reliable
approximation of the threshold between Regions I and II. This
measurement method is applied to thirteen distinct datasets of linear
waves generated from varied input frequencies in the subcritical
regime.  Figure \ref{fig:exp_vg} compares the measurement of $\xi_g$
with the group velocity $c_g(\omega)$ \eqref{eq:BBM_grp_vel_wavemaker}
at different input angular frequencies $\omega$.  The angular
frequencies are measured and averaged over the periodic wave
region. The error bars presented in Figure \ref{fig:exp_vg} represent
the uncertainties resulting from the linear fit of $x_g(t)$ and the
variation in the measured wave frequency.  Figure \ref{fig:exp_vg}
verifies the threshold $\xi_g$ as the linear BBM group velocity
prediction experimentally.

\begin{figure}
    \centering
    
    \begin{minipage}{.45\textwidth}
        \centering
        \sidesubfloat[]{\includegraphics[height=0.7\textwidth]{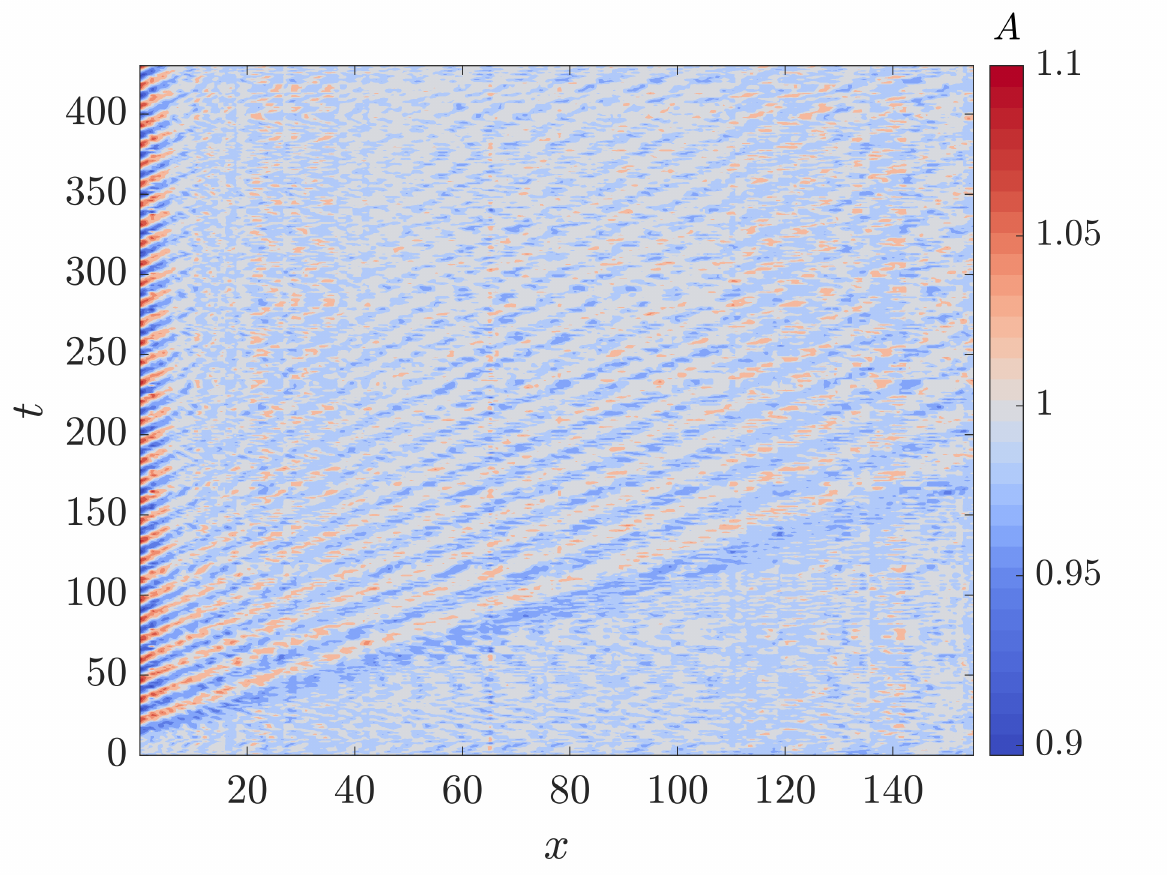} \label{fig:exp_contour_super}}
    \end{minipage}
    \hspace{0.1in}
   \begin{minipage}{.45\textwidth}
        \centering
        \sidesubfloat[]{\includegraphics[height=0.7\textwidth]{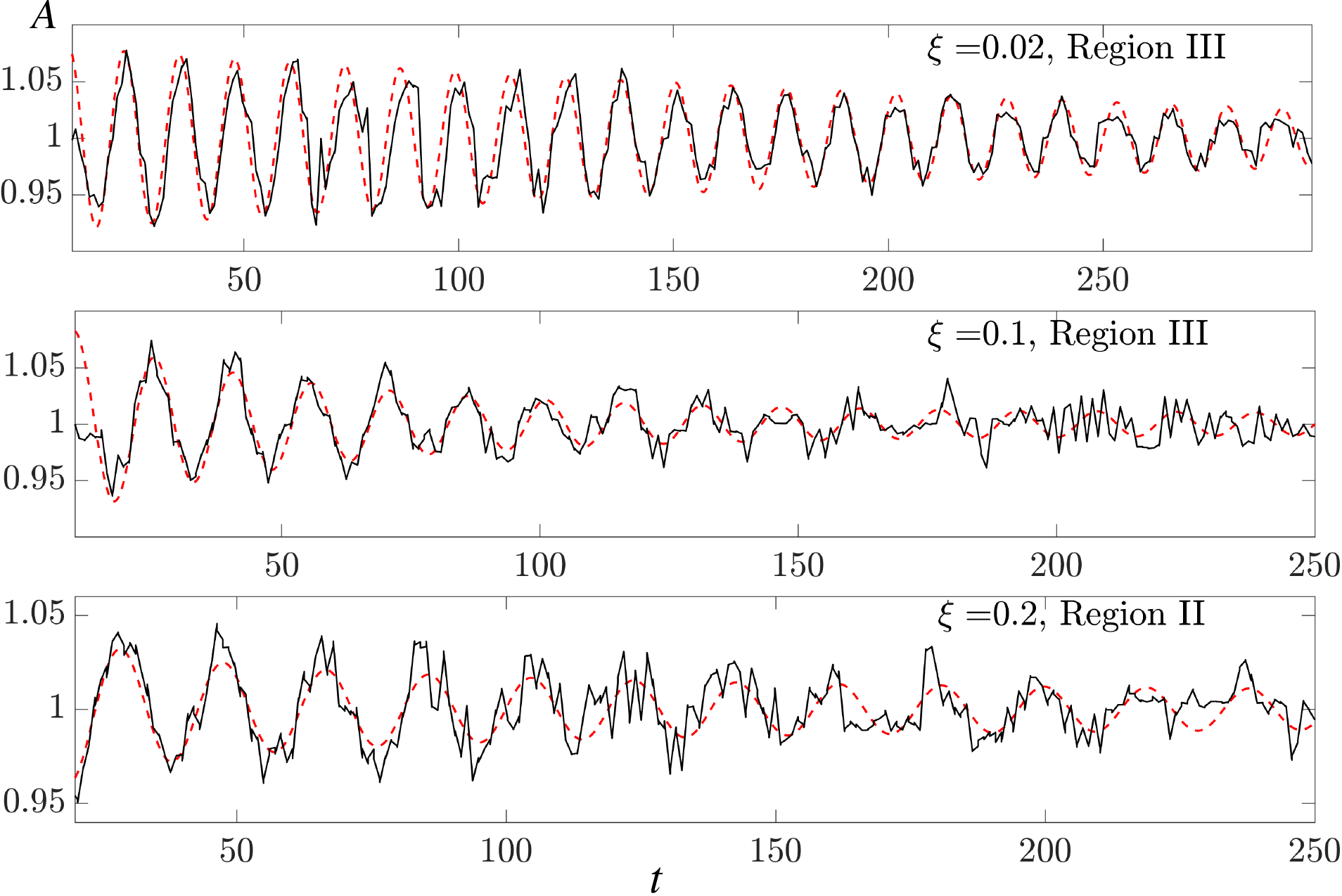} \label{fig:exp_conduit_super} }
    \end{minipage}
    \caption{Viscous core-annular flow experiment in the supercritical regime. 
   (a) Space-time contour of normalized cross-sectional area $A$ for time-periodic injection of core fluid into a uniform background. The measured angular frequency near the boundary is $\omega=0.51\pm0.02$, and the nondimensionalization scales are $(L,U)=(0.32\pm0.01\,\text{cm}, 1.04\pm0.02\,\text{cm/s})$. (b) Comparison of the experimental data in (a) (black solid) with the long-time asymptotic solutions in Regions II \eqref{eq:LBBM_soln_II} and III \eqref{eq:BBM_soln_III} (red dashed) at different, fixed $\xi$.}
    \label{fig:exp_super}
\end{figure}

Figure \ref{fig:exp_super} depicts an experiment in the supercritical regime with an input frequency
greater than the critical value, i.e., $\omega_0>1/2$ (see \cite{mao2023experimental} for an experimental study of the critical value). The wave
appears to be time-periodic near the boundary and rapidly decays as it
propagates. Near the boundary, the measured angular frequency
is $\omega=0.51\pm0.02$, very close to the critical value.
The space-time data in Fig.~\ref{fig:exp_contour_super} is compared with the long-time
asymptotic solutions of the linearized BBM wavemaker problem for fixed $\xi$ in Fig.~\ref{fig:exp_conduit_super}. Applying
the transformation \eqref{eq:conduit-to-BBM} and fixing $\xi$ as
defined in \eqref{eq:exp_Xi}, where the shift parameters are $x_0=0.3$ $(\approx 0.1\,\text{cm})$
and $t_0=-1$ (imaging started $\approx 0.3\,\text{s}$ ahead), the experiment is compared to the linear BBM long-time
asymptotic solutions \eqref{eq:LBBM_soln_II} in Region II and
\eqref{eq:BBM_soln_III} in Region III of Figure
\ref{fig:LBBM_soln_descent} using the mean wave parameters $(a=0.08,\omega_0=0.51)$.  The predicted threshold between
Regions II and III at $\omega_0=0.51$ is $\xi=0.12$ (cf.~Equation
\eqref{eq:BBM_kappa_pos_2}).  Figure
\ref{fig:exp_conduit_super} shows that the wave amplitude, frequency, and phase of the experimental result agree with the theoretical
predictions.  

\section{Conclusions}
\label{sec:conclusions}

This work presents the long-time, asymptotic analysis of
initial-boundary value problems for linear, scalar evolution equations
on the half-line with time-periodic boundary conditions. We obtain the
Dirichlet-to-Neumann (D-N) map for a general third-order model with
Dirichlet boundary data that is sufficiently smooth and asymptotically
periodic in time. The uniqueness of the D-N map is established
(Theorem \ref{dn-t}) if and only if the well-known radiation condition in
physics (Corollary \ref{rad-cond}) holds.  By employing the D-N map, we
obtain the leading-order asymptotic solution for fixed $x$ and
$t\to\infty$.  These results hinge upon the assumption of asymptotic
time-periodicity of the solution.

This assumption is proven in the specific case of the wavemaker
problems for the linear KdV and linear BBM equations with zero initial
condition and sinusoidal boundary condition with frequency $\omega_0$
(Theorems \ref{lkdv-t} and \ref{lbbm-t}). Using the Fokas method and
steepest descent analysis, we obtain the leading-order, asymptotic
solution as $t\to\infty$ that is uniform in $x \ge 0$.  Crossing the
interval of frequencies $(-\omega_{cr},\omega_{cr})$ identifies a
bifurcation in solution behavior.  When $|\omega_0| < \omega_{cr}$,
the solution is asymptotically the propagating wave
$u(x,t) \sim \sin(k_0x - \omega_0 t)$, $t \to \infty$ for
$x \in [0,c_g(\omega_0)t)$ where $c_g(\omega_0) = \omega'(k_0)$ is the
group velocity and $k_0 = k_0(\omega_0)$ is the unique branch of the
dispersion relation $\omega_0 = \omega(k_0)$ such that
$c_g(\omega_0) > 0$ and $\mathrm{Im}\,k_0 = 0$.  This choice of $k_0$
corresponds to the radiation condition.  When
$|\omega_0| > \omega_{cr}$, the asymptotic solution is the
time-periodic, spatially decaying, causal solution
$u(x,t) \sim \exp(-\mathrm{Im}(k_0)x)\sin( \mathrm{Re}(k_0)x -
\omega_0t)$ for $x = \mathcal{O}(1)$, $t \to \infty$ where
$k_0 = k_0(\omega_0)$ is the unique branch of $\omega_0 = \omega(k_0)$
such that $\mathrm{Im}\, k_0 > 0$ that again satisfies the radiation condition.

We observe quantitative agreement between viscous core-annular fluid
wavemaker experiments and the linearized BBM predictions.
In particular, we observe the bifurcation in solution behavior at the
group velocity and the convergence of spatial profiles to propagating
waves, waves that temporally decay $\mathcal{O}(t^{-1/2})$, and waves
that exponentially decay in space
$\mathcal{O}(\exp(-\mathrm{Im}(k_0)x))$.

We conclude this section by noting that the Fokas method has also
found application in the rigorous study of initial-boundary value
problems for nonlinear evolution equations that are not necessarily
completely integrable. In those non-integrable but still physically
relevant cases, the linear solutions obtained via the Fokas method
provide the starting point for developing the corresponding
well-posedness theory. This is done by treating the nonlinear problem
as a perturbation of its forced linear counterpart and invoking the
contraction mapping theorem, 
e.g. see the works \cite{fokas2016korteweg, fokas2017nonlinear,
  himonas2020well} and the references therein.  Additionally, there
are open questions regarding the long-time analysis of the nonlinear
wavemaker problem.  Whitham modulation theory is a promising analysis
tool for integrable and non-integrable equations although it presently
lacks a rigorous foundation.  As the analogue of the inverse scattering transform for boundary value problems, the Fokas method could also be used
for the rigorous analysis of the nonlinear wavemaker problem in
integrable systems \cite{fokas_linearization_1996,monvel_mkdv_2004,fokas_explicit_2010,deconinck_implementation_2021}.

\bibliographystyle{amsplain}
\bibliography{reference}

\end{document}